\documentclass[12pt,a4paper]{article}
\newcommand{\com}[1]{}
\newcommand{\lb}[1]{\label{#1}}\newcommand{\reff}{\ref}\newcommand{\rf}[1]{$(\reff{#1})$}
\newcommand{\cit}[1]{\cite{#1}}\newcommand{\bibi}[1]{\bibitem{#1}} 
\usepackage{amssymb,enumerate,amsmath,makeidx}
\usepackage[latin1]{inputenc}
\usepackage[T1]{fontenc}
\newcommand{\eq}[2]{\begin{equation}#2\,\,\, \,\lb{#1}\end{equation}}
\def\sep{\hspace{-2mm}.}
\def\noi{\noindent}
\def\med{\medskip}
\def\sm{\smallskip}
\def\np{\newpage}
\def\pr{\noi{\sl Proof\;\;\sep }\;\;}
\def\ep{\hfill$\square${\,}\medskip}
\def\eep{\medskip}
\def\tq{\;;\;}
\def\ssi{\Leftrightarrow}
\def\bul{\noindent$\triangleright$\;\;}

\newcommand{\scal}[2]{\langle#1\,|\,#2\rangle}
\renewcommand{\sec}[2]{\section{#2}{\lb{#1}}}
\newcommand{\sub}[2]{\subsection{#2}{\lb{#1}}}
\newcommand{\ssub}[2]{\subsubsection{#2.}{\lb{#1}}}

\newtheorem{lem}{{\bf Lemma}}[section]
\newtheorem{theor}[lem]{{\bf Theorem}}
\newtheorem{deff}[lem]{{\bf Definition}}

\newtheorem{conj}[lem]{{\bf Conjecture}}
\newtheorem{rem}[lem]{{\bf Remark}}
\newtheorem{coro}[lem]{{\bf Corollary}}
\newtheorem{exemple}[lem]{{\bf Example}}
\newtheorem{question}{{\bf Question}}
\newtheorem{exercice}[lem]{{\bf Exercise}}
\newtheorem{propo}[lem]{{\bf Proposition}}
\newcommand{\df}[2]{\begin{deff}{\lb{#1}}{\bf\sep}{\rm#2}\end{deff}}

\newcommand{\exo}[2]{\begin{exercice}\lb{#1}{\bf\sep}\rm#2\end{exercice}}
\newcommand{\ex}[2]{\begin{exemple}\lb{#1}{\bf\sep}\sl#2\end{exemple}}

\newcommand{\lm}[2]{\begin{lem}\lb{#1}{\bf\sep}\sl#2\end{lem}}
\renewcommand{\th}[2]{\begin{theor}{\lb{#1}}{\bf\sep}{\sl#2}\end{theor}}
\newcommand{\Lm}[2]{\begin{lem}\lb{#1}{\bf\sep}\sl#2\end{lem}}

\newcommand{\cor}[2]{\begin{coro}{\lb{#1}}{\bf\sep}{\sl#2}\end{coro}}
\newcommand{\prop}[2]{\begin{propo}{\lb{#1}}{\bf\sep}{\sl#2}\end{propo}}

\def\as{\asymp}

\def\a{\alpha}
\def\b{\beta}
\renewcommand{\d}{\delta}
\def\eps{\varepsilon}
\def\e{\varepsilon}
\def\f{\varphi}
\def\g{\gamma}
\def\k{\mathbb{K}}
\renewcommand{\l}{\lambda}
\def\q{\theta}

\def\sig{\sigma}

\def\w{\omega}
\def\A{\mathcal{A}}
\def\C{\mathbb{C}}
\def\D{\Delta}

\def\F{\mathbb{F}}
\def\H{\mathbb{H}}

\def\N{\mathbb{N}}
\def\O{\mathcal{O}}
\def\OO{\mathrm{O}}
\def\0{\{\id\}}
\def\P{\mathcal{P}}
\def\Q{\mathbb{Q}}
\def\R{\mathbb{R}}
\def\RP{\mathbb{R}\mathbf{P}}
\def\S{\mathbb{S}}
\def\sy{\mathcal{S}}
\def\T{\mathcal{T}}
\newcommand{\GA}{\mathop{\mathrm{GA}}\nolimits}
\newcommand{\GL}{\mathop{\mathrm{GL}}\nolimits}
\newcommand{\SL}{\mathop{\mathrm{SL}}\nolimits}
\newcommand{\PSL}{\mathop{\mathrm{PSL}}\nolimits}
\def\SO{\mathrm{SO}}
\def\Z{\mathbb{Z}}
\def\W{\Omega}
\renewcommand\bar[1]{\overline{#1}}
\renewcommand\le{\leqslant}
\renewcommand\ge{\geqslant}

\newcommand{\be}{\begin{enumerate}[{\bf a.}]}
\newcommand{\bee}{\begin{enumerate}[{\bf(i)}]}
\newcommand{\ee}{\end{enumerate}}\newcommand{\eee}{\end{enumerate}}
\renewcommand{\~}{\widetilde}
\renewcommand{\hat}{\widehat}
\def\modd{\!\!\mod}
\def\ta{{\rm^t}}

\def\card{{\rm card\,}}
\def\md{{\rm Med\,}}
\def\id{{\bf id}}
\def\bij{{\rm Bij\,}}
\def\aff{{\rm Aff\,}}
\def\tr{{\rm tr\,}}
\def\im{{\rm Im\,}}
\def\re{{\rm Re\,}}
\def\mat{{\rm Mat\,}}

\def\fix{{\rm Fix\,}}
\def\isom{{\rm Isom\,}}
\def\isomp{{\rm Isom^+\,}}
\def\isomm{{\rm Isom^-\,}}

\def\diag{{\rm diag}}
\def\gaf{{\sc gaf}}
\def\gag{{\sc gag}}

\usepackage{esvect}\renewcommand{\v}[1]{\vv{#1}} 			

\title{
{\Large Fixating Group Actions}
}
\author{\normalsize
Guido Ahumada, Bernard Brighi, Nicolas Chevallier, and Augustin Fruchard}
\date{\normalsize January 24, 2019}
\makeindex
\begin{document}
\maketitle
\begin{abstract}
A group of  bijections $ G $ acting on a set $ X $  is said
{\sl with fixed points} (abbreviated as
\gaf\;from the french ``groupe \`a points fixes'')
if any element of $ G $ has at least one fixed point in $X$.
The $G$ group is said {\sl with a common fixed point} (abbreviated as \gag\;
for ``groupe à point fixe global'')
if there is $ x \in X $ fixed by all elements of $ G $.
The group $ G $ is said {\sl fixating} if any subgroup of $ G $
which is a \gaf\;is automatically a \gag.
The article explores which groups are fixating.
The situation depends on the assumptions made on the
group of bijections and on the support set $ X $.
For example the group of
isometries of the Euclidean space $ \R^n $ is fixating for $n\le3$
but not for $n\ge4$. The case of isometries of 
elliptic and hyperbolic spaces is also considered, as well as that of
isometries of some discrete sets. As we will see, the situation depends
upon the fact whether the so-called median inequality is satisfied
in the ambient space or not.
Besides, some of our constructions of nonfixating groups 
rely on the existence of free subgroups of certain linear groups.
\end{abstract}

\noi Key-words: Group of bijections, fixed point, isometry,
median inequality, tree.
\med

\noi MSC2010 Classification: 	 	51M04, 51M09, 57M60, 20F65.
\np
\tableofcontents
\np
%
%
%
%
\sec{1.}{Introduction}
It is easy to find a group of bijective transformations, each having a
fixed point but without a common fixed point: for instance,
in the symmetric group of $\{1,2,3,4,5\}$, the subgroup generated by
the cycle $(123)$ and the double transposition $(12)(45)$,
or the group of rotations of the two dimensional sphere, or the 
group of homeomorphisms of the unit disk.
Such a group will be called {\sl eccentric}\index{eccentric group}.
In addition to the existence of a fixed point for each bijection of
the group, which information would imply the existence of a
common fixed point?
This additional information can be the invariance of a geometric
structure, the commutativity of the group or another algebraic
property, the uniqueness of the fixed point for each nontrivial bijection,
or a combination of the former informations.

The invariance of a geometric structure can often be stated as follows:
The group of bijections is a subgroup of a larger group $G$. The fact
that this information implies the existence of a common fixed point
can be seen as a property of the larger group $G$.
We shall say that a group  $G$ of bijections of a set $X$ is 
{\sl fixating}\index{fixating group} if it contains no eccentric
subgroup. 

To our knowledge, the notion of fixating group has not been
considered in earlier works. Notice that this
notion is not intrinsic to the group but depends on the action as a 
group of bijections. We will skip the word ``action'' for short.
The aim of this paper is to explore which groups of bijections
are fixating and, to a lesser extent, to find some sufficient
conditions for a group of bijections to have a common fixed point.

About fixating groups, we shall see that many things can happen,
depending on the nature and the dimension of the set $X$, and on
the nature of the bijections. We have paid particular attention to
isometry groups of some classical spaces.

When $X$ is a metric space, $\isom X$ denotes the group of
isometries of $X$. If moreover $X$ can be oriented, $\isomp X$
denotes the subgroup of isometries preserving the orientation. 
Given an integer $n\ge1$ we denote by
$\R^n$ the $n$-dimensional Euclidean space, 
$\Z^n$ the lattice of integer points in $\R^n$,
$\H_n$ the  $n$-dimensional hyperbolic space,
$\S_n$ the $n$-dimensional sphere, and
$\RP_n$ the $n$-dimensional projective space
(we prefer to use the exponent $n$ for cartesian products only).
Our results about these spaces are the following.

\bul
The groups $\isom\R^n$ and $\isomp\R^n$ are fixating if and only if $n\le3$.

\bul
The groups $\isom\H_n$ and $\isomp\H_n$ are fixating if $n\le3$, 
and nonfixating if $n\ge5$. When $n=4$, the question is open.

One could think that, for each family of spaces $\F_n=\R^n$, $\H_n$,
$\S_n$ or $\Z^n$, there exists a critical integer $n_0$ such that the
isometry group of $\F_n$ is fixating if and only if $n\le n_0$.
This holds for $\R^n$ (with $n_0=3$) and $\H_n$
(with $n_0=3$ or $4$) but neither for $\S_n$ nor for $\Z^n$: 

\bul
The group $\isomp\S_n$ is fixating if and only if $n=1$ or $3$. 

\bul
The group $\isom\Z^n$ is fixating for all $n\ge1$,  whether $\Z^n$ is
equipped with the Euclidean norm or the $L^1$ norm.

An important ingredient for the existence of global fixed points
is the {\sl median inequality}, see formula~\rf6 at the beginning
of Section~\reff{2.}. This inequality holds in the Euclidean and the
hyperbolic spaces but neither in the spheres nor in the projective spaces.

This inequality has been introduced by F.~Bruhat and J.~Tits
in~\cite{bt}\index{Bruhat, Fran\c cois}\index{Tits, Jacques} and their
result is as follows:
If a metric space $(E,d)$ is complete and if the median inequality holds,
then each isometry group with a bounded orbit admits a global fixed point, 
see Theorem~\reff{t6.3} and Corollary~\reff{c4.4}.

Moreover, we show that an isometry group $G$ in a complete metric space
satisfying the median inequality is fixating provided that it has
a normal fixating subgroup $H$ such that $G/H$
is cyclic, see Corollary~\reff{c4.5}. This is the reason why $\isom\R^n$
and $\isomp\R^n$ are fixating for the same dimension $n$; the same holds
for $\isom\H_n$ and $\isomp\H_n$. The median inequality is crucial for this
result: We shall see that $\isomp\S_3$ is fixating while $\isom\S_3$ is not.
Indeed:

\bul
The group $\isom\S_n$ is fixating if and only if $n=1$.

About the projective space we obtain:

\bul
The group $\isom\RP_n$ is fixating if and only if $n=1$;

\bul
The group $\isomp\RP_n$ is fixating if $n=1$ and nonfixating if $n=2$
or $n\ge4$. We do not know whether $\isomp\RP_3$ is fixating or not.
\sm

Another important ingredient in our study is the existence of some free
subgroups in certain linear groups. It is used to construct eccentric groups.
A very general result by A. Borel~\cite{bo} ensures this existence,
although some elementary results lead to explicit examples that are
enough for our need. 

The article is organized as follows. Section \reff{1bis.} is devoted
to  notations and some preliminary results. Section \ref{2bis.} is
about affine groups and Section \reff{2.} about isometry groups.
Isometry groups of classical spaces are studied in detail in
Section \reff{sec:globalisant}: first in Euclidean spaces, then in
hyperbolic spaces, and at last in elliptic spaces. Section
\reff{5.} is devoted to results in discrete spaces. We first deal with 
permutation groups, then with $\isom\Z^n$, and at last with some
graphs. Concerning graphs, we first present a result by
J.-P.~Serre\index{Serre, Jean-Pierre} about
fixed points of finitely generated isometry groups of trees~\cite s,
then we extend this result to a family of colored graphs.
Section~\reff {7.} provides exercises; the solutions are given in
Appendix~\reff{8.4}. Appendix~\reff{8.1} is a short introduction to
hyperbolic geometry and Appendices~\reff{8.2} and~\reff{8.5}
give details of some proofs.

We leave several questions open. It may happen that some
answers already exist in the literature, or that
the reader solves some of them.
In that case we would be grateful to the reader who will inform us!   
%
%
%
\sec{1bis.}{Preliminaries}
\sub{2.1}{Notations}
{\noi \bf Intervals.} We denot by $\,]a,b[\,$ 
the open interval from $a$ to $b$, in order to avoid confusion
with the pair $(a,b)$.
\med

{\noi \bf Groups.}  Given a group $G$, we denote $H\le G$ if
$H$ is a subgroup of $G$ and $H \unlhd G$ if $H$
is a normal subgroup of $G$.

If $ A $ is a subset of a group $ G $,
the subgroup of $ G $ generated by $ A $ is denoted $\langle A\rangle$.
When $ A $ contains a small number of elements, we skip the braces.
Thus we write
$$
\langle f,g\rangle=\bigcap_{f,g\in H\le G}H
=\big\{f^{i_1}g^{j_1}\dots f^{i_n}g^{j_n}\tq
n\in\N,\;i_k,j_k\in\Z\,\big\}.
$$
We call {\sl cyclic} a group, finite or not, generated by a single element.

The {\sl commutator} of $f$ and $g$ is $[f, g] = f^{- 1} g^{- 1} fg $.

For a set $ X $, $\bij X$ denotes the group of bijections of $X$.
\med

{\noi\bf Metric spaces.}
In a metric space $ (X, d) $, an {\sl isometry}\index{isometry}
$f:X\to X$ is a {\sl bijection} such that $d(f(x),f(y))=d(x,y)$ for
all $x,y\in X$.
As already said in the introduction, $\isom X$ is the group of
isometries of $X$ and $\isomp X$ is the subgroup of those preserving
the orientation if $X$ is orientable.
 
For an element $ x \in X $ and a real number $ r> 0 $, 
$ B (x, r) = \{y \in X \tq d (x, y) <r \} $ denotes the open ball of
center $x$ and radius $r$ and $B'(x,r)=\{y\in X\tq d(x, y)\le r\}$
is the corresponding closed ball.
Given two points $ a $ and $ b $ of $ X $, $ \md (a, b) $ denotes the
mediator\index{mediator} of $ a $ and $ b $:
\[
\md(a,b)=\{c\in X \tq d(a,c)=d(b,c)\}.
\]
{\noi \bf Affine and Euclidean Spaces.}
Given an affine bijection $f:\R^n\to\R^n$,
$\v f$ is the {\sl associated linear map}, defined by
$\v f(x)=f(x)-f(\v0)$.
The mapping $ f \mapsto \v f $ is a morphism of groups.
In particular we have $[\v f,\v g]=\v{[f,g]}$.

Similarly, if $F$ is an affine subspace of $\R^n$ then $\v F$ denotes
the associated vector subspace.

For a subset $ A \subset \R^n $, we denote by
$ \aff A $ the affine space generated by $ A $, i.e. the intersection of all
affine subspaces of $ \R^n $ containing $ A $.
If $ A = \{a_1, \dots, a_n \} $,
we shall write $ \aff (a_1, \dots, a_n) $ instead of 
$ \aff \big (\{a_1, \dots, a_n \} \big) $.
\med

{\noi \bf Hyperbolic Spaces.}
We will use the model of the upper half-space with the Poincar\'e metric, 
see Appendix \ref {8.1}.
\sub {2.2} {{\small GAF}, {\small GAG} and fixating group}
A {\sl group of bijections} $ (X, G) $ is the data of a set $ X $ and
a subgroup $ G $ of $ \bij X $.
Given a bijection $ g: X \to X $, its set of fixed points is
$$
\fix g = \{x \in X \tq g (x) = x \}.
$$
A group of bijections $ (X, G) $ is called 
{\sl with fixed points} (abbreviated as \gaf\index{g@\gaf})
if $ \fix g $ is nonempty for all $ g \in G $.
We say that $ (X, G) $ is a {\sl group with a common fixed point}
(\gag\index{g@\gag}) if
$$
\fix G: = \bigcap_{g \in G} \fix g \ne \emptyset.
$$
A \gaf\, which is not a \gag\, is called
{\sl eccentric}\index{eccentric group}.
With the above vocabulary we say that $ (X, G) $ is
{\sl fixating}\index{fixating group}
if for any subgroup $H\le G$ we have
$$
(X, H)\mbox{\;\gaf\;}\ssi\;(X,H)\;\mbox{\gag,}
$$
i.e. if it does not contain any eccentric subgroup.
Observe that any subgroup of a fixating group is fixating.
We shall omit the set $ X $ when the context is clear.
In the same way we will say that an action $ \rho $  of a group $ G $ 
on a set $ X $ is fixating if $ (X, \rho (G)) $ is a  fixating group
of bijections.
We will use the following result several times.
\prop {r2.1} {
Let $ f, g $ be two bijections on a set $ X $.
\be
\item\lb{i}
If $ f $ and $ g $ commute then $ g (\fix f) = \fix f $.
\item\lb{ii}
If $ \fix f $ is a singleton $ \{x_0 \} $ then
$ x_0 \in \fix g $ for any $ g $ commuting with $ f $.
\item\lb{iii}
Let $ G $ be a group of bijections on $ X $ and let $ H \unlhd G $.
For all $ g \in G $ we have $ g (\fix H) = \fix H $.
\ee
}
\pr
{\bf\reff{i}}.
If $ x \in \fix f $ then $ g (x) = g (f (x)) = f (g (x)) $ hence $ g (x) \in \fix f $.
This proves $ g (\fix f) \subseteq \fix f $.
Besides, $ g^{- 1} $ also commutes with $ f $, so
$ g^{- 1} (\fix f) \subseteq \fix f $. We then obtain
$ \fix f = g (g^{- 1} (\fix f)) \subseteq g (\fix f) $,
hence the equality follows.
\medskip

{\noi\bf\reff{ii}}.
Results from item {\bf a}.
\medskip

{\noi\bf\reff{iii}}.
Let $ g \in G $. For all $x\in\fix H$ and all $h\in H$,
we have $ k = g^{- 1} hg \in H $, hence $ x \in \fix k $.
Therefore we have $ g (x) = g (k (x)) = h (g (x)) $,
hence $g(x)\in\fix h$.
This shows that $g(\fix H)\subseteq\fix H$ for all $g\in G$.
Applied to $g^{-1}$, this gives 
$\fix H=g\big(g^{-1}(\fix H)\big)\subseteq g(\fix H)$.
\ep

We immediately deduce from item~\reff{ii} a first sufficient condition
for a group of bijections to be a \gag.
\prop{c2.2}{
Let $ G $ be a group of bijections on $ X $. If $ G $ is Abelian and
if there exists $ f_0 \in G $ having a single fixed point
then $ G $ is a \gag.
}
We will see in Section \reff {2bis.} that each of the words ``Abelian''
and ``unique'' is needed.
\med

Let us end this section with remarks of an algebraic nature.
\sm\\
1.
The concept ``fixating'' is compatible
with the product: If $ (X_1, G_1) $ and $ (X_2, G_2) $ are two fixating 
groups of bijections, it is easily checked that the action
of the product group $ G_1 \times G_2 $ on $ X_1 \times X_2 $ is
fixating.
\sm\\
2.
The notion ``fixating'' however is not compatible with
the induction of Frobenius~\cite {fr}\index{Frobenius, Ferdinand Georg}.
More precisely, let $G$ be a group, $H$ a subgroup of $G$, and $R$
a system of class representatives modulo $ H $.
An action of $ H $ on a set $ Y $ induces an action of $ G $
on $ X = R \times Y $ defined by $ g (r, y) = (r ', h (y)) $ 
where $ r' \in R $ and
$ h \in H $ are uniquely determined by $ gr = r'h $.

Section~\reff {5.2} on permutation groups provides an example
where the action $ (Y, H) $ is fixating but the induced action is not.
We consider the group $ G $ of permutations of $ \{1,2,3,4,5 \} $ and
$ H $ the subgroup of permutations fixing $ 5 $, seen as acting on
$ Y = \{1,2,3,4 \} $. As system of representatives, we choose 
the transposition $ r_i = (i \, 5) $
for $ 1 \le4 $ and $ r_5 = \id $.

On the one hand $(Y,H)$ is fixating by Theorem~\reff{t6.1},
on the other hand $ (X, G) $ is not:
Let $ K $ be the subgroup of $ G $ generated by the
permutations $ (123) $ and $ (12) (45) $. One finds

$$
K=\{\id,(123),(132),(12)(45),(13)(45),(23)(45)\},
$$
$\fix(12)(45)=\{(r_3,5)\}$, $\fix(13)(45)=\{(r_2,5)\}$,
$\fix(23)(45)=\{(r_1,5)\}$, and
$\fix(123)=\fix(132)=\{(r_4,4),(r_4,5),(r_5,4),(r_5,5)\}$,
hence $K$ is an eccentric subgroup of $(X,G)$. 
\sm\\
3.
The notion ``fixating'' strongly depends on the set on which
the group acts. In Exercise~\reff{e1} we introduce an
intrinsic notion of globalization:
A group is said to be {\sl superfixating} if, for any set $ X $ and
any morphism $ \rho: G \to \bij X $, the pair $ (X, \rho (G)) $ is fixating.
This notion finally has a rather limited interest: The result of
Exercise~\reff{e1} is that a group is superfixating if and
only if it is cyclic (finite or not).
On the other hand, the additive group $ \Q $ is {\sl finitely superfixating}
in the following sense:
If $ X $ is a finite set and $ \rho: \Q \to \bij X $ a morphism then
$ (X, \rho (\Q)) $ is fixating, cf. Exercise~\reff{e2}.
%
%
%
%
\sec{2bis.}{Groups of affine bijections} 
Among the geometric structures, the affine structure is the most 
fundamental one.
It is therefore natural to begin our study with groups of affine bijections
of $\R^n$. The results are rather negative except for dimension one.
\prop{l3.1}{
The group of affine bijections of $ \R $ is fixating.
}
\pr
Let $ H $ be a group of affine bijections of $\R$ which is a \gaf.
Since an affine bijection of  $\R $ different from the identity
has at most one fixed point, by Proposition~\reff {c2.2}
it is enough to prove that $ H $ is Abelian.
If $ f $ and $ g $ are in $ H $ then the commutator
$ [f, g] $, which also has a fixed point, is not a nontrivial translation.
Since the linear group on $\R$ is Abelian, 
this commutator can only be the identity.
\ep

The commutativity of the linear group and the uniqueness of the fixed points
are the two ingredients of the previous proposition.
Both ingredients are characteristics of the dimension one.
In larger dimension, adding only one of the hypotheses --- commutativity
or uniqueness of the fixed points --- is not enough to prove that a \gaf\,
is a \gag\, as shown in Examples~\reff{e2.2} and~\reff{e2.4}
below.
\ex{e2.2}{
Let $a\in\R\setminus\Q$ and let $f,g$ be the affine transvections 
of $\R^2$ given by
$$
f(x,y)=(x+y+1,y)~\mbox{ and }~g(x,y)=(x+ay,y).
$$
Then the group $G=\langle f,g\rangle$
is Abelian and eccentric.
}
\pr
It is straightforward to prove that $G$ is Abelian.
For $(m,n)\in\Z^2\setminus\{(0,0)\}$ we have
$f^mg^n(x,y)=\big(x+(m+na)y+m,y\big)$ hence
$\fix(f^mg^n)$ is nonempty: It is the straight line $D_{m,n}$ of equation
$y=\tfrac{-m}{m+na}$. It follows that $G$ is a \gaf.
We also have $D_{m,n}\cap D_{m',n'}=\emptyset$ provided that $mn'\ne m'n$
hence $G$ is not a \gag.
\ep

{\noi\sl Remarks.}
\sm\\
1.
If $ \R^2 $ is equipped with the discrete distance given by $ d (a, b) = 1 $ if
$ a \ne b $ and $ d (a, a) = 0 $, then the group $ G $ of the previous example
is an eccentric Abelian group of isometries.
Example~\reff {e4.3} of the next section is an example of an
eccentric and Abelian  group of isometries in a Hilbert vector space
of infinite dimension.
On the other hand Theorem~\reff{t2.4} shows that there does not exist
any eccentric Abelian group of isometries of a Euclidean
or hyperbolic space of finite dimension.
\sm\\
2.
The above shows that the group of affine bijections of
$\R^2$ is nonfixating.
For $n\ge3$, extending the previous maps $f$ and $g$
by the identity on the last $n-2$ components
we obtain an eccentric group of affine bijections of $\R^n$.
As a consequence the group of affine bijections of
$\R^n$ (with $n\ge1$) is fixating if and only if $n=1$.
\sm\\
3. 
We have not explored whether there is an eccentric Abelian group of affine
maps acting on $\Q^2$. On the other hand Example 3.3 below is valid with
$\Q$ instead of $\R$.
\sm\\
4.
Given an integer $d\ge1$ and a field $\k$, one can also ask whether
the group $\GL(d,\k)$ acting on $\k^d\setminus\{\v0\}$ is fixating.
The answer is that it is fixating only in the obvious cases
$\GL(1,\k)$ and $\GL(2,\F_2)$, where $\F_2=\{0,1\}$,
cf. Exercises~\reff{e16} and~\reff{e17}.
\ex{e2.4}{
Let $b\in\R^2\setminus\{\v0\}$ and let $f,g\in\bij\R^2$ be the affine bijections
$$
f:x\mapsto\v f(x)\qquad\mbox{ and }\qquad g:x\mapsto\v g(x)+b,
$$
where $\v f$ and $\v g$ are the elements of\; $\SL(2,\R)$ with matrices
$$
\mat(\v f)=A
=\left(\begin{matrix}\hfill0&1\\
-1&3\end{matrix}\right)\qquad\qquad 
\mat(\v g)=B
=\left(\begin{matrix}-1&-1\\
\hfill5&\hfill4\end{matrix}\right).
$$
Then the group $G_1=\langle f,g\rangle$ is eccentric.
More precisely, every element of $G_1\setminus\{\id\}$ has a unique fixed
point but $G_1$ is not a \gag.
}
We recall that $\SL(2,\Z)$ is the set of $2$ by $2$ matrices
with integer coefficients and determinant~$1$.
The proof is based on the following two results.
\lm{lem:libreSL2}{
{\rm (\cite{n} Chapter VIII)}
Let $ G_0 $ be the subgroup of\; $ \SL (2, \Z) $
generated by the matrices $ A $ and $ B $
above.
Then $ G_0 $ is free and any matrix $ M \in G_0 \setminus \{I \} $ has
trace different from $2 $. 
}
The proof is in Appendix~\reff {8.2}.
We deduce that $ 1 $ is not an eigenvalue of $ M $
for any matrix $ M \in G_0 \setminus \{I \} $, since $ \det M = 1 $.
\Lm {lem:valeurpropre1} {
Let $ h $ be an affine bijection on $ \R^n $ such that $ 1 $ is
not an eigenvalue of $ \v h $. Then $ \fix h $ is a singleton.
}
\pr
We have $h(x)=x\ssi x-\v h(x)=h(\v0)\ssi x=(\id-\v h)^{-1}(h(\v0)).$
\ep

{\noi\sl Proof of Example~\reff {e2.4}}.
Let $ h \in G_1 \setminus \{\id \} $ and let $ \g_1, \g_1 '\dots \g_r, \g'_r \in \Z $,
all nonzero except possibly $ \g_1 $ and $ \g'_r $, such that
$ h = f^{\g_1} g^{\g'_1} \cdots g^{\g'_r} $.
We have $ \v h = \v f^{\g_1} \v g^{\g'_1} \cdots \v g^{\g'_r} \ne \id $, 
so  $ 1 $ is not an eigenvalue of $ \v h $  according to 
Lemma~\reff {lem:libreSL2}, hence $ h $ has a unique fixed point by
Lemma~\reff {lem:valeurpropre1}.
However there is no fixed point
common to all elements of $ G_1 $ since
the fixed points of $ f $ and $ g $ are distinct. 
\ep

To finish this section, let us recall two classical results on
the existence of global fixed points
for groups or sets of affine applications.
The first result is the Markov-Kakutani\index{Markov, Andre\"\i \ Andre\"\i evitch}\index{Kakutani, Shizuo}\index{Markov-Kakutani theorem}\index{Theorem!Markov-Kakutani}
theorem and the second is the 
Kakutani theorem,\index{Kakutani theorem}\index{Theorem!Kakutani} see \cite{mk,ds,ru,k}. These two results require a
compactness assumption and an additional hypothesis.
Exercise~\reff{e9} presents a finite version of these results.
It is due to R.~Antetomaso~\cite{a}.
\th{Markov-Kakutani}{
{\rm (Markov-Kakutani, see for example~\cite {ru}
Theorem 5.23, p.140)\index{Markov, Andre\"\i \ Andre\"\i evitch}\index{Kakutani, Shizuo}} \,
Let $K$ be a convex compact nonempty subset of a 
separated topological vector space $E$ and let $G$ be a set of
{\rm affine and continuous} maps that commute and leave
$ K $ stable. Then $ K \cap \fix G $ is nonempty.
}
\th{kakutani}{
{\rm(Kakutani, cf.~\cite{ru} Theorem 5.11, p.127)\index{Kakutani, Shizuo}}\,
Let $K$ be a compact convex subset of a locally convex
topological vector space $E$ and $G$ be an equicontinuous group of 
affine bijections leaving $K$ stable.
Then $K\cap\fix G$ is nonempty.
}
The literature contains a large number of recent works containing
supplements and extensions of these results.
Among these, below is a complement to Theorem~\reff{Markov-Kakutani}
due to Anzai and Ishikawa~{\rm \cite {ak}}\index{Anzai, Kazuo}\index{Ishikawa, Goo}:
Under the assumptions of Theorem~\reff{Markov-Kakutani},
with moreover $ E $ locally convex,
if $ G $ is a finitely generated group, $ G = \langle T_1, \dots, T_n \rangle $,
then for all $ \a_j \in \,] 0,1 [\, $ satisfying $ \sum_{j = 1}^n \a_j = 1 $ we have
$\fix\big(\sum_{j=1}^n\a_jT_j\big)=\fix G$.
\medskip

The theorem below can be deduced from Theorem~\reff{kakutani} but
we  give an independent proof.
\th{orbite-bornee}{
A group of affine bijections of $ \R^n $ having a bounded orbit
is a \gag.
}
\pr
Let $ G $ be such a group and let $ x $ be an element of $ \R^n $ whose orbit
$ \O_x = \{f (x) \tq f \in G \} $ is bounded.
Let $ K $ be the closed convex hull of $ \O_x $.
Like $ \O_x $, $ K $ is invariant by all maps $ f \in G $.
Let $ F = \aff K $, the affine subspace generated by $ K $.
The interior of $ K $ relative to $ F $ is nonempty hence $ \l_F (K)> 0 $,
where $ \l_{F} $ denotes the Lebesgue measure of $ F $
(in the case where $K$ is a singleton, $\l_F$ is the counting measure).
Since $ K $ is compact, we also have $ \l_F (K) <+ \infty $.
The maps $ f \in G $ are affine, so send the measure $ \l_F $
to a multiple of itself. Since they send $ K $ on $ K $,
they preserve $\l_F$.
The centroid of $K$ for the restriction of $ \l_{F} $ to $ K $,
defined by $ \frac1 {\l_F (K)} \int_Kxd \l_F (x) $, is therefore
fixed by all $ f \in G $.
\ep
%
%
%
%
\sec{2.}{Groups of isometries}
\sub{sec:mediane}{The median inequality}
The Bruhat-Tits\index{Bruhat, Fran\c cois}\index{Tits, Jacques} fixed point theorem~\cite{bt}
gives a sufficient condition for a group of isometries on a metric space
to be a \gag: It is enough that the space satisfies the median inequality
below and that the group has a bounded orbit.
\df{d4.1}{
We say that a metric space $ (X, d) $ {\sl satisfies the median
inequality}\index{median inequality} if
\eq6{
\forall x,y\in X\;\exists m\in X\;\forall z\in X\quad
d(z,m)^2\le\tfrac12\big(d(z,x)^2+d(z,y)^2\big)-\tfrac14\,d(x,y)^2.
}
}
It is easy to prove that the point $ m $ is unique and that
$d(x,m)=d(y,m)=\tfrac12\,d(x,y)$, cf.~Exercise~\reff{e3}.
We say that $m$ is {\sl the midpoint\index{midpoint} of $\{x,y\}$}
and we denote it $m(x,y)$.

When $X$ is a Euclidean  space, or more generally a pre-Hilbert space,
\rf6 is actually an equality called the {\sl parallelogram identity}
and $m$ is the usual midpoint of the segment $[x,y]$, cf.~Exercise~\reff{e4}.a.
Conversely, it is known that a normed vector space satisfying~\rf6 is
necessarily pre-Hilbert, cf.~Exercise~\reff{e4}.b.

A combinatorial tree with its usual distance does not satisfy~\rf6
(an edge has no midpoint)
but its realization as a real metric space satisfies it.
Trees are also the only graphs with this property.
Complete Riemannian manifolds which are simply connected with 
negative sectional curvature, especially hyperbolic spaces
with their usual distance, satisfy~\rf6, see~\cite {bt}.
This is a consequence of the comparison theorem of Rauch~\cite {ec},
cf. Exercise~\reff{e6}.

The median inequality makes it possible to associate a single center with any
bounded subset of a complete metric space. Let $(X,d)$ a
metric space and $A$ be a nonempty bounded subset of $X$.
For every $x\in X$, let
$$
r(x,A)=\inf\{r>0\tq A\subseteq B'(x,r)\}=\sup\{d(x,a)\tq a\in A\}.
$$
We define the {\sl radius of $A$} by
$$
r_A=\inf\{r(x,A)\tq x\in X\}.
$$
If there exists $ x_0 \in X $ such that $ r_A = r (x_0, A) $,
we will say that $ x_0 $ is a {\sl center\index{center} of $A$}.
In this case the closed ball $ B '(x_0, r_A) $ is a
{\sl ball circumscribed to $ A $}.
\lm{l4.2b}{ {\rm\cite{bt}}
If $ (X, d) $ is a complete metric space satisfying \rf6
then any nonempty bounded subset of $ X $ has a unique center.
}
\pr
Let $ A $ be a nonempty bounded subset of $ X $, let $ x, y \in X $, and
let $ m $ be the midpoint of $ \{x, y \}. $
Writing \rf6 for all $ a \in A $ we get
\begin{align*}
d(m,a)^2&\le\tfrac12\big(d(x,a)^2+d(y,a)^2\big)-\tfrac14d(x,y)^2\\
&\le\tfrac12\big(r(x,A)^2+r(y,A)^2\big)-\tfrac14\,d(x,y)^2,
\end{align*}
from which we successively deduce
\[
r_A^2\le r(m,A)^2\le\tfrac12\big(r(x,A)^2+r(y,A)^2\big)-\tfrac14\,d(x,y)^2,
\]
and
\eq5{
\tfrac12\,d(x,y)^2\le r(x,A)^2+r(y,A)^2-2r_A^2.
}
We deduce the uniqueness of a possible center:
If $ x $ and $ y $ are centers of $ A $ then $ r (x, A) = r (y, A) = r_A $,
and~\rf5 implies $ x = y $.

To prove the existence, let $ (x_n)_{n \in \N} $ be a sequence such that
$r(x_n,A)$ tends to $r_A$ as $n$ tends to $+\infty$. 
Taking $x = x_n$ and $y = x_{n+p}$ in~\rf5 we obtain
$$
\tfrac12\,d(x_n,x_{n+p})^2\le r(x_n,A)^2+r(x_{n+p},A)^2-2r_A^2\to0
$$
uniformly in $ p $ as $ n $ goes to infinity.
Thus the sequence $ (x_{n})_{n \in \N} $ is Cauchy,
therefore has a limit $ \ell $ verifying $ r (\ell, A) = r_A $,
hence $ \ell $ is a center of $ A $.
\ep

The Bruhat-Tits\index{Bruhat, Fran\c cois}\index{Tits, Jacques} fixed point theorem\index{Bruhat-Tits theorem}\index{Theorem!Bruhat-Tits} is stated
as follows.
\th{t6.3}{
{\rm\cite{bt}}\,
Let $G$ be a group of isometries of a complete metric space
$ (X, d) $ verifying the median inequality~\rf6.
If there is a nonempty bounded subset of $ X$ which is invariant by all
the elements of $ G $ then $ G $ is a \gag.
}
\pr 
Let $ A $ be a nonempty bounded subset of $ X $, invariant by any
$ g \in G $, and let $ a $ be the center of $ A $.
Then, for all $ g \in G $, $ g (B '(a, r_A)) = B' (g (a), r_A) $
is the ball circumscribed to $ g (A) = A $, therefore by uniqueness of the center,
$g(a) = a $. As a consequence $ a \in \fix G $ hence $ G $ is a \gag.
\ep

We immediately deduce the following result.
\cor{c4.4}{
Let $ G $ be a group of isometries of a Euclidean or hyperbolic space.
If $ G $ has a bounded orbit then $ G $ is a \gag.
}
Some of the following results will be used in Section~\reff{sec:globalisant}.
\df{d2}{
Let $(E, d)$ be a complete metric space verifying~\rf6.
A subset $C$  of $E$ is called {\sl half-convex}\index{half-convex} if
for all $x,y\in C$ the midpoint of $\{x,y\}$ is in~$C$.
}
It can be easily shown that a half-convex closed subset of a normed vector
space is convex in the usual sense.

The following proposition asserts the existence and uniqueness of an
``orthogonal  projection'' on the set of fixed points of a group of isometries.
We split it into three statements, each having its own interest.
We recall that $ d (x, A) = \inf_{a \in A} d (x, a) $.
\prop{p4.5}{
Let $ (E, d) $ be a complete metric space verifying~\rf6.
\be \item \lb {p4.5a}
If $ C \subset E $ is a closed half-convex part of $ E $ then,
for all $ x \in E $, there exists a unique $ y \in C $ such that
$ d (x, C) = d (x, y) $. This point is denoted $ y = \pi_Cx $.
\item \lb {p4.5b}
If $ g $ is an isometry of $ E $ then $ \fix g $ is closed and
half-convex.
\item \lb {p4.5c}
If $ G $ is a group of isometries of $ E $ then $ \fix G $ is closed and
half-convex.
\ee
}
\pr{\bf\reff{p4.5a}}.
Let $ \d = d (x, C) $.
By definition, for any $ \eps> 0 $, there exists $ y \in C $
such that $ d (x, y) <\d + \eps $. If $z\in C$ also satisfies this inequality,
since $ d \big (x, m (y, z) \big) \ge \d $, \rf6 then gives
\eq{2b}{
d(y,z)^2\le8\d\eps+4\eps^2.
}
For each $ n \in \N^* $, let $ y_n \in C $ be such that $ d (x, y_n) <\d + \tfrac1n $.
By~\rf{2b} the sequence $ (y_n)_{n \in \N^*} $ hereby defined is
 Cauchy, hence converges to a point $y\in C$ satisfying $d(x,y) = \d$;
this proves the existence.
Inequality~\rf{2b} also proves the uniqueness.
\medskip 

{\noi\bf\reff{p4.5b}} and {\noi\bf\reff{p4.5c}}.
The proof is straightforward.
\ep
\th{t2.6}{
Let $ (E, d) $ be a complete metric space satisfying~\rf6,
let $ G $ be a group of isometries of $ E $,
and let $ H \unlhd G $ be such that $ G / H $
is cyclic. If $ G $ is a \gaf\;and $ H $ a \gag\;then $ G $ is a \gag.
}
\pr
Let $ \eps \in G $ be such that $ \eps H $ generates $ G / H $.
Denote $ F = \fix H $.
Since $ H $ is normal in $ G $, we have $ g (F) = F $ for all $ g \in G $
according to Proposition~\reff {r2.1}.\reff {iii}, in particular $ \eps (F) = F $.

Let $ x \in \fix \eps $.
By uniqueness of the orthogonal projection and since $ \eps $ is an isometry,
one has $ \eps (\pi_Fx) = \pi_{\eps (F)} \eps (x) = \pi_Fx $, hence
$ \pi_Fx \in \fix \eps $. Since $ G = \langle \eps, H \rangle $,
we obtain $ \pi_Fx \in \fix G $.
\ep
\cor{c4.5}{
Let $ (E, d) $ be a complete metric space satisfying \rf6 and
let $ H \unlhd G \le\isom E $ be such that $ G / H $ is solvable and finite.
If $ H $ is fixating then $ G $ is fixating.
In particular a group of isometries of $ E $ is fixating provided that it
contains a fixating subgroup of index $ 2 $.
}
\pr
We first assume that the quotient $ G / H $ is cyclic.
Let $ G_1 \le G $ be a \gaf; then $ H_1 = G_1 \cap H $ is a \gaf,
therefore a \gag\;since $ H $ is fixating.
Besides, $ G_1 / H_1 $ is isomorphic to a subgroup of $ G / H $,
therefore cyclic, hence $ G_1 $ is a \gag \, according to Theorem~\reff {t2.6}.

Since finite Abelian groups are products of cyclic groups,
under the hypothesis that $ G / H $ is solvable and finite,
there exists a finite sequence
$ H = H_0 \unlhd \cdots \unlhd H_n = G $ such that
for each $ i = 1, \dots, n $
the quotient $ H_{i} / H_{i-1} $ is cyclic.
The result is then successively applied to the cyclic quotients.

For the last assertion, if $ H $ is a subgroup of $G$ of index $2$
then $ H $ is normal in $ G $ and $ G / H $ is cyclic of order $ 2 $.
\ep
\np
{\noi\sl Remarks.}
\sm\\
1.
Our proof is not valid if $G/H$ is only solvable.
The right notion in our context is $G/H$
{\sl polycyclic}~\cit{se}\index{polycyclic group},
that is, $G/H$ admits a finite sequence 
$\{e\}=H_0 \unlhd\cdots\unlhd H_n=G/H$
with $ H_ {i} / H_ {i-1} $ cyclic.
We do not know whether or not our result holds when $G/H$ is only assumed
solvable and of finite type.
\sm\\
2. 
The fact that the ambient space satisfies~\rf6 is essential:
We will see in Section~\reff {section:globalisantelliptique} that
$\isomp\S_3 $ is fixating while $\isom\S_3$ is not.
\sub{sec:resoluble}{Solvable\index{solvable group} subgroups:
The Euclidean and hyperbolic cases}
We fix an integer $ n \ge0 $. The notation $ \F_n $ will indicate either the 
Euclidean space $ \R^n $, or the hyperbolic space $ \H_n $.
We know that $ \F_n $ satisfies~\rf6,
see Exercise~\reff{e6} for the case of $\H_n$.
\th{t2.4}{
Let $ G $ be an {\rm Abelian} group of isometries of $ \F_n $.
If $ G $ is a \gaf\,then $ G $ is a \gag.
}
\pr
We proceed by induction on the dimension $n$.
For $ n = 0 $ the result is trivial.
Now let $ n \ge1 $ and assume that the property is true for all $ k <n $.

If $ G = \{\id \} $ we are done. Otherwise let $ f \in G \setminus \{\id \} $.
Then $ F = \fix f $ is a subspace (affine or hyperbolic) of
$ \F_n $, of some dimension $ k <n $.
Let $ g \in G $. As $ f $ and $ g $ commute,
we have $ g (F) = F $ according to Proposition \reff {r2.1}.{\reff {i}}.
So for all $ g \in G $ the restriction of $ g $ to $ F $,
denoted by $ g_{| F} $, is well defined from $ F $ to $ F $ and is an
isometry of $ F $ which is itself isometric to $ \F_k $.

By hypothesis, $ \fix g $ is nonempty.
Let $ x_g \in \fix g $ and set $ y_g = \pi_Fx_g $ 
given by Proposition~\reff{p4.5}.
Since $ g $ is an isometry and by uniqueness of the orthogonal projection,
we have
$g(y_g)=\Pi_{g (F)}g (x_g)=\Pi_Fx_g =y_g$.
Thus, for every $ g \in G $, $ g_{| F} $ has at least one fixed point $ y_g $.

Let $ G_F = \{g_{| F} \tq g \in G \} $.
Then $ G_F $ is a \gaf \, on a space of dimension $ k <n $, hence is a \gag \
by induction. Since $ \fix G_F = F \cap \fix G $, we deduce
that $ \fix G $ is nonempty hence $ G $ is a \gag.
\ep

The following example shows that the finite dimension is necessary.
\ex{e4.3}{
{\rm
Let $ E = \ell^2 (\N, \R) $, the space of square summable real sequences;
it is a Hilbert space. Let $ h_k $ be the symmetry of center $ 1 $ on
the $ k $-th component, i.e. the isometry of $ E $ defined by
$$
h_k(x_0,x_1,\dots)=(x_0,\dots,x_{k-1},2-x_{k},x_{k+1},\dots).
$$
Let $ G_n = \langle h_0, \dots, h_n \rangle $ and let $G = \bigcup_{n \in \N} G_n $.
It is immediate that $ G $ is Abelian.
Let $ (e_n)_{n \in \N} $ denote the canonical basis of $ E $ and
$s_n = \sum_{k = 0}^ne_k$.
We have $ s_n \in \fix f $ for all $ f \in G_n $, so $ G $ is a \gaf.
By contradiction, if $ G $ were a \gag \, and if
$ x = (x_0, x_1, \dots) \in \fix G $ then
for all $ n \in \N $ we would have $ x_n = 1 $ but the constant
sequence equal to 1 is not in $ E $ hence $ G $ is eccentric.
}}
Theorem~\reff {t2.4} can be generalized by changing the assumption
``Abelian'' into ``solvable''.
\th{t2.5bis}{
Let $ G $ be a {\em solvable} group  of isometries of $ \F_n $.
If $ G $ is a \gaf\;then $ G $ is a \gag.
}
\pr
Recall that a group $G$ is called {\sl solvable}\index{solvable group}
if there exists a finite and increasing sequence of subgroups
$\{\mathbf{e}_G\}=H_0\unlhd H_1\unlhd\dots\unlhd H_p=G$ 
(i.e. each normal in the next one) such that all the quotients
$ H_{k + 1} / H_k $ are Abelian.
The {\sl solvability index of $ G $} is the smallest integer $p\ge0$
with this property.
This integer is reached for example by taking the sequence of derived groups:
We put $ G_0 = G $ and, for $ k \ge0 $, $ G_{k + 1} = G'_k = [G_k, G_k] $,
the group generated by the commutators of $ G_k $.
Finally we choose $ H_k = G_{p-k} $.

The proof is by induction on the solvability index of $ G $.
The property is trivial for $ p = 0 $.
Assume that it is true for any  group of solvability index $ p-1 $ and
let us show it for $ G $.

Let $ G_1 = [G, G] $. Since $ G $ is a \gaf, $ G_1 $ is a
\gaf\;hence a \gag \;by induction hence $ F = \fix G_1 $ is nonempty.
Since $ G_1 $ is normal in $ G $ we have $ g (F) = F $ for
all $ g \in G $.

Now the group $ G / G_1 $ acts naturally on $ F $:
If $ \bar g = \{gh \tq h \in G_1 \} $ is an element of $ G / G_1 $ and if $ x \in F $,
then $ \bar g (x): = g (x) $ does not depend on the choice of the representative $ g \in \bar g $
since $ h (x) = x $ for all $ h \in G_1 $.

We claim that the pair $ (F, G / G_1) $ is a \gaf.
Indeed, since $ G $ is a \gaf, if $ g \in G $ and $ x \in \fix g $ then,
as in the proof of Theorem~\reff {t2.4}, by uniqueness of the
orthogonal  projection, the projection $ \pi_Fx $ is also in $ \fix g $.
As a consequence $ F \cap \fix g \ne \emptyset $, which gives $ \fix \bar g \ne \emptyset $
for all $ \bar g \in G / G_1 $ seen as isometry of $ F $.

Since $ G / G_1 $ is an Abelian group and $ F $ is a space (Euclidean or
hyperbolic) of finite dimension, $ G / G_1 $ is a \gag \, according to
Theorem~\reff{t2.4}. Any fixed global point $ x \in F $ of $ G / G_1 $ is then
fixed by any element of $ G $, so $ G $ is itself a \gag.
\ep
%
%
\sec{sec:globalisant}{Groups of isometries of the classical spaces}
In the whole Section~\reff{sec:globalisant}, $ n $ is a strictly positive integer.
\sub{sec:globalisanteuclidien}{The Euclidean\index{Euclidean space} case}
Recall that $\isom\R^n$ is the group of isometries of $\R^n$ equipped with the
usual Euclidean distance and $\isomp\R^n$ is the subgroup of those preserving
the orientation. It is known that the elements of $\isom \R^n$
are affine applications, cf.~Exercise~\reff {e4b}.
In this Section~\reff{sec:globalisanteuclidien} we prove the following result.
\th{t3.1}{
The group $ \isom \R^n $ is fixating if and only if $ n \le3 $.
}
We leave to the reader the pleasure to show
that $ \isom \R^2 $ is fixating, see Exercise~\reff{e0}.
We will show successively that $ \isom \R^3 $ is fixating,
then that $ \isomp \R^4 $ is nonfixating, which will imply
that $ \isom \R^n $ is not  fixating for $ n \ge4 $.
\medskip

{\noi\bf The case of dimension 3}
\medskip

Let us recall that the elements $ f \in \isomp \R^3 $ such that
$ \fix f \ne \emptyset $
are either the identity or the rotations around an axis $ \fix f $.
Those with an empty set of fixed points
are the translations and the screw displacements
(a {\sl screw displacement} is the Abelian product of a rotation
and a translation of nonzero vector parallel to the axis of the rotation).

The following lemma is a key step to prove that $\isom\R^3$ is
fixating. We will have this same step in the proof that the group of
isometries of the hyperbolic space of dimension $3$ is fixating,
see Lemma~\ref{l4.5}.
\lm{l3.4}{
If $ f, g \in \isomp \R^3 $ are such that $ \fix f \cap \fix g = \emptyset $,
then there is $ h \in \langle f, g \rangle $ such that $ \fix h = \emptyset $.
}
\pr
If $ \fix f $, $ \fix g $, or $ \fix (f^{- 1} g) $ is empty we are done.
Otherwise let $ a \in \fix (f^{- 1} g) $ and let $ b = f (a) = g (a) $.
We have $ b \ne a $ since $ \fix f \cap \fix g = \emptyset $.
So for all $ c \in \fix f $ we have
$ d(a, c) = d (f (a), f (c)) = d (b, c) $. Therefore $ \fix f $ is in
$ \md (a, b) $, the mediator plane of $ a $ and $ b $.
The same holds for $ \fix g $.
Since $ \fix f \cap \fix g = \emptyset $, $ \fix f $ and $ \fix g $ are
two parallel lines. It follows that $ \v {[f, g]} = \id $.
If $ f $ and $ g $ were commuting, we would have $ f (\fix g) = \fix g $
by Proposition \reff {r2.1}.{\reff {i}}, in contradiction with
$f\ne\id$. It follows that $ [f, g] $ is a nontrivial translation
hence $ \fix [f, g] = \emptyset $.
\ep
\prop{p3.6}{
The group $ \isom \R^3 $ is fixating.
}
\pr
According to Corollary~\reff {c4.5}, it suffices to show that
$ \isomp \R^3 $ is fixating.
Let $ G \le \isomp \R^3 $ be a \gaf. We must show that $ \fix G \ne \emptyset $.
Lemma~\reff {l3.4} already implies $ \fix f \cap \fix g \ne \emptyset $
for all $ f, g \in G $.

We call {\sl half-turn} a symmetry about a straight line, called its {\sl axis}. 
We will use the following fact: If the product of two half-turns is
a half-turn then their axes are orthogonal and secant.

If $ G \setminus \{\id \} $ contains only half-turns then,
either $ G = \{\id, f \} $ where $ f $ is a half-turn, or
$ G = \{\id, f_1, f_2, f_3 \} $ where $ f_1, f_2, f_3 $ are three half-turns of axes
orthogonal and pairwise secant, so secant all three in one
point, hence $ G $ is a \gag.

From now on we assume that there exists $ f \in G \setminus \{\id \} $ that is
not a half-turn. If $ \fix g = \fix f $ for all $ g \in G \setminus \{\id \} $
we are done: $ \fix G = \fix f \ne \emptyset $. Let us now assume that there
exists $ g \in G \setminus \{\id \} $ such that $ \fix g \ne \fix f $.
Then $ \fix f $ and $ \fix g $ are two straight lines crossing at some point
denoted $ \w $.
Set $ P = \aff (\fix f \cup \fix g) $, the affine plane containing
$ \fix f $ and $ \fix g $.

\begin{center}
\begin{picture}(200,110)(-30,-20)
\put(-30,-25){\line(1,0){240}}
\put(-30,-25){\line(0,1){120}}
\put(-28,-20){$P$}

\put(0,0){\line(3,2){130}}
\put(128,70){$\fix h$}
\put(-5,20){\line(6,-1){160}}
\put(147,0){$\fix g$}
\put(77,89){\line(1,-3){34}}
\put(49,80){$\fix f$}
\put(23,15.33){\circle*{3}}
\put(20,2){$b$}
\put(106,1.7){\circle*{3}}
\put(96,-8){$\w$}
\end{picture}
\end{center}

We will be done if we show that $\w\in\fix G$.
We proceed by contradiction. Assume that there exists $h\in G$
such that $ \w \notin \fix h $.
Since $ \fix f \cap \fix h $ and $ \fix g \cap \fix h $ are nonempty,
we have $ \fix h \subset P $.
Let $ a \in \R^3 $ be such that $ f (a) = b \in \fix g \cap \fix h $.
We have $ b \in P $ and $ b \notin \fix f $ but $ f $ is not a half-turn,
so $ a \notin P $.
Let $ \~g = f^{- 1} gf $. We have $ \~g (a) = a $ and $ \~g (\w) = \w $,
hence $ \fix \~g = \aff (a, \w) $, which is a straight line disjoint from
$ \fix h $, a contradiction.
\ep

{\noi\sl Remark.}
This result can also be proved by copying the proof of Lemma~\reff {l4.6}.
\medskip

{\noi\bf The case of higher dimensions}
\prop{p3.9}{
The group $ \isomp\R^4 $ is nonfixating.
As a consequence, the group $ \isom\R^4 $ is nonfixating.
}
\pr
We reproduce below the construction by Wagon~\cite w\index{Wagon}
of a free subgroup of rank 2 in $ \SO_4 $ whose action on the sphere
$ \S_3 $ is without fixed point.
Let $ \q \in \R $ be such that $ \cos \q $ is transcendent,
for example let $ \q = 1 $, and let $ \sig $ and $ \tau $ be 
the elements of $ \SO_4 $ of matrices
\eq{41}{
S_4=\left(\begin{matrix}\hfill\cos\q&-\sin\q&0&0\\
\sin\q&\cos\q&0&0\\
0&0&\cos\q&-\sin\q\\
0&0&\sin\q&\cos\q\end{matrix}\right)
\quad\mbox{ and }\quad
T_4=\left(\begin{matrix}\hfill\cos\q&0&0&-\sin\q\\
0&\cos\q&-\sin\q&0\\
0&\sin\q&\cos\q&0\\
\sin\q&0&0&\cos\q\end{matrix}\right).
}
\lm{l5.5}{ {\rm(\cite w Theorem 5.2, p.53)}
The subgroup $ G_0 $ of $ \SO_4 $ generated by
$ \sig $ and $ \tau $ is free. 
Moreover, $ 1 $ is not an eigenvalue
of any element of  $ G_0 \setminus \{\id \} $.
}
For the convenience of the reader, we have written the proof in
Appendix~\reff{8.5}.
\med

Now choose the affine rotations $ \sig $ and
$ \~\tau: x \mapsto \tau x + a $ with $ a \ne\v0 $, for example $ a = (1,0,0,0) $.
Let $G$ be the subgroup of $ \isomp \R^4 $
generated by $ \sig $ and $ \~\tau $.
Then $ G $ is free and
$ \fix g $ is a singleton for all $ g \in G $.
Since $ \fix \sig \cap \fix \tau = \emptyset $
we deduce that $ G $ is eccentric.
\ep

{\noi\sl Remark.}
The existence of free subgroups of $ \SL (2, \R) $ and $ \SO_4 $
whose elements, apart from identity, never admit
$ 1$ as eigenvalue, is the essential ingredient of
constructions of eccentric subgroups of affine applications or
of affine isometries (Example \reff {e2.4} and Proposition
\reff {p3.9}). We have used explicit examples of
such subgroups. These subgroups, although sometimes difficult to
exhibit, are not exceptional. Indeed it can be shown
thanks to the Baire theorem that, if $ G $ is a closed subgroup of
 $ \GL (n, \R) $, then the set of pairs of elements of
$G$ generating a free group is either empty or contains
a dense $G_{\delta}$, i.e. a
countable intersection of dense open subsets of  $ G \times G $.
The same result holds with the additional constraint on the eigenvalue $1$.
Let us add that A.~Borel\index{Borel, Armand} has proved a very general result
encompassing our constructions, see~\cite{bo}:
\smallskip\\
{\sl 
If $ G $ is a semi-simple linear algebraic group defined
on $ \R $ then the set of $n$-tuples of $G(\R)^n$ generating
a free group contains a dense $G_{\delta}$ subset.}
\prop{prop:Ingeq4}{
For any integer $ n \ge4 $, the group $ \isom\R^n $ is nonfixating.
}
\pr
For $ n \ge5 $,
let $ \sig_n $ and $ \tau_n $ be the elements of $ \isomp \R^n $ of matrices
$$
S_n = \left (\begin {matrix} \hfill S_4 & 0 \\0 & \id_{n-4} \end {matrix} \right)
\mbox { and }
T_n = \left (\begin {matrix} \hfill T_4 & 0 \\0 & \id_{n-4} \end {matrix} \right)
$$
respectively, where $S_4$ and $T_4$ are defined in~\rf{41}.
Let $G$ be the subgroup of
$ \isom \R^n $ generated by $ \sig_n $ and $ \~\tau_n: x \mapsto \tau_nx + a $
with $ a \in \R ^n \setminus \{\v0 \} $.
It is easy to check that $ G $ is free and eccentric.
\ep

{\noi\bf The case of non Euclidean norms.}
\medskip

Let us  endow $ \R^n $ with an arbitrary norm, 
and let $ G $ denote  the group of  isometries associated with this norm.
According to the Mazur-Ulam theorem\index{Mazur-Ulam theorem},\index{Theorem!Mazur-Ulam}\index{Mazur, Stanislaw}\index{Ulam, Stanislaw} the elements of $ G $ 
are affine maps, cf.~\cite{mu} or
Exercise~\reff {exo:mazurulam} and its solution.
The group $ \v G $ of the linear parts of the elements of
$ G $ is closed and bounded in the vector space of endomorphisms
of $ \R^n $, so $ \v G $ is compact.
By a classical argument we can construct a scalar product which is
invariant by the elements of $ \v G $. Therefore $ G $ can be seen as a subgroup of
$ \isom \R^n $ and $ G $ is thus fixating if $ n \le 3 $.
When $ n \ge4 $, Proposition \ref {prop:Ingeq4} does not apply
and indeed $ G $ may be fixating for some norms:
\prop{p5.6}{
Let $ N $ be a norm on $ \R^n $ and let $ G $ be the group of isometries
associated with $ N $. If the linear group associated with $ \v G $ is finite then $ G $
is fixating. In particular, if $ N $ is one of the usual $ N_p $ norms
with $ p \in [1, + \infty] \setminus \{2 \} $ then $ G $ is fixating.
}
\pr
Let $ H $ be a subgroup \gaf \, of $ G $.
Since the only translation in $ H $ is $ \id $, the morphism
$ \f: H \to \v G, \, f \mapsto \v f $ is injective, so $ H $ is finite.
We check that the point $ \w = \frac1 {| H |} \sum_{h \in H} h (\v0) $ is fixed by
all elements of $ H $ hence $ H $ is a \gag.

If $ N = N_p $ then $ \v G $ contains all the permutations of the axes,
so a scalar product invariant by the elements of $ \v G $ is
necessarily proportional to the usual Euclidean scalar product.
The Euclidean unit sphere touches the unit ball for $ N_p $
exactly at the vertices of
the hyper-octahedron $E=\{\pm\v e_1,\dots,\pm\v e_n\}$.
The set $E$ is invariant by $\v G$, so $\v G$ is finite.
\ep
\sub{sec:globalisanthyperbolique}{The hyperbolic\index{hyperbolic space} case}
We use the Poincar\'e model of the half-space:
$\H_n=\R_{>0}\times\R^{n-1}$ endowed with the Poincar\'e metric given by
$ds^2=\frac1{x_n^2}\sum_{i=1}^ndx_i^2$, cf. Appendix~\reff{8.1} for details.

%
\th{t6.2b}{
If $n\le3$ then $\isom\H_n$ is fixating.
If $n\ge5$ then $\isom\H_n$ is nonfixating.
}
We conjecture that the group $\isom\H_4$ is nonfixating.
\med

We will show first that $ \isom \H_n $ is nonfixating if $n\ge5$
and then show that $ \isom \H_2 $ and $ \isom \H_3 $ are fixating.
The case $n=1$ is obvious.
\med
\np
{\noi\bf The case of dimensions higher than 5}
%
\prop{p4.21}{
For any $n\ge5$, the group $\isom\H_n$ is nonfixating.
}
\pr
With each isometry of the Euclidean space of dimension $ n-1 $ we associate
an isometry of the hyperbolic space of dimension $ n $ in the following way: 
If $ f $ is an isometry of $ \R^{n-1} $ then
the map $ F: \R^{n-1} \times \R_{> 0} \rightarrow \R^{n-1} \times \R_{> 0}$
defined by $ F (x, t) = (f (x), t) $ is an isometry of $ \H_n $.
The mapping $ f \in \isom \R^{n-1} \mapsto F \in \isom \H_n $
is a morphism of groups and we have $ \fix F = \fix f \times \R_{> 0} $.
Therefore the image by this morphism of an eccentric subgroup of
$ \isom \R^{n-1} $ is an eccentric subgroup of $ \isom \H_n $.
Since, for $ n \ge5 $, $ \isom \R^{n-1} $ is nonfixating,
$ \isom \H_n $ is nonfixating.
\ep
\med

{\noi\bf The case of dimension 2}
\smallskip

Although this case can be deduced from the $3$-dimensional case
(cf. Exercise~\reff{e12}) we prefer to present proofs specific to dimension
$2$ because they are more elementary and they provide
a good introduction to hyperbolic geometry.
We use the model of the complex half-plane $ \H_2 = \{z \in \C \tq \im z> 0 \} $.
The positive isometries of $ \H_2 $ are the homographies
$$
h_{a,b,c,d}:z\mapsto\dfrac{az+b}{cz+d}~\mbox{ with }~a,b,c,d\in\R,\,ad-bc=1.
$$
The other isometries of $\H_2$, called {\sl negatives}, are the homographies
composed with the symmetry $z\mapsto-\bar z$. The group $\isomp\H_2$
of positive isometries is called the {\sl M\"obius group}\index{M\"obius group}.
We recall that the mapping
$$
\Phi:\isomp \H_2\to\PSL(2,\R),\;h_{a,b,c,d}\mapsto\{M,-M\}~\mbox{ with }~
M=\left(\begin{matrix}a&b\\c&d\end{matrix}\right)
$$
is an isomorphism of groups.
We call {\sl matrices associated with an isometry} $ h $ the elements of
$ \Phi (h) $. We define the {\sl trace} of an isometry $ h $ by
$ \tr h = | \tr M | $ where $ M \in \Phi (h) $.
We will use the following result, the proof of which is in
Appendix~\reff {sec:h2}.
\lm{l6.3}{
Let $h\in\isomp \H_2\setminus\{\id\}$.
\be\item
We have $\fix h\ne\emptyset$ if and only if $\tr h<2$.
Moreover, in that case, $\fix h$ is a singleton.
\item
In the case where $\fix h=\{i\}$, the associated matrices $M$ and $-M\in\Phi(h)$
are the orthogonal  matrices
$R\big(\tfrac\q2\big)=\left(
\begin{matrix}\cos\tfrac\q2&-\sin\tfrac\q2\smallskip\\
\sin\tfrac\q2&\cos\tfrac\q2\end{matrix}\right)$
and $R\big(\tfrac\q2+\pi\big)$. 
We then say that $h$ is a {\rm rotation of center $i$  and angle}
$\q\in\R/(2\pi\Z)$,
and we denote $h=r_\q$.
\item
In the general case, if $\fix h=\{z_0\}$ and if $\f\in\isomp \H_2$ is such that
$\f(i)=z_0$ then $\f^{-1}h\f$ fixes $i$ hence is a rotation $r_\q$, 
and the angle $\q$ does not depend on the choice of $\f$. We will say that
$h$ is a {\rm rotation of center $z_0$ and angle} $\q$.
\item\lb{d}
For every $x>0$ and every $\q\in\R/(2\pi\Z)$, the rotation of center $ix$ and
angle $\q$ is associated with the matrices
$\pm\left(\begin{matrix}
\cos\tfrac\q2&-x\sin\tfrac\q2\smallskip\\
x^{-1}\sin\tfrac\q2&\cos\tfrac\q2
\end{matrix}\right)$.
\ee}
An element $h\in\isomp \H_2$ is called {\sl elliptic} if $\tr h<2$,
{\sl parabolic} if $\tr h=2$, and {\sl hyperbolic} if $\tr h>2$.

We divided by $ 2 $ the angle in the matrices so that the angle defined
in the statement above corresponds to the usual notion of angle:
The rotation $ r_\q: z \mapsto \tfrac {cz-s} {sz + c} $ with
$ c = \cos \tfrac \q2 $ and
$ s = \sin \tfrac \q2 $ turns approximately the points very close to its
center $ i $ by an angle $ \q $ and not $ \tfrac \q2 $.

In addition we will use that two segments of $ \H_2 $ of the same
length are positively isometric:
If $ x, y, x ', y' \in \H_2 $ are such that $ d (x, y) = d (x ', y') $
then there exists $ f \in \isomp \H_2 $ such that $ f (x) = x '$ and $ f (y) = y' $;
see Section \ref {sec:h2} for a proof.

As a preliminary we show the following result.
\lm{l4.2}{
If $f,g\in \isomp \H_2$ are such that $\fix f\cap\fix g=\emptyset$
then there exists $h\in\langle f,g\rangle$ such that $\fix h=\emptyset$.
}
\pr
If $\fix f$ or $\fix g$ is empty we are done. Now we assume that 
$\fix f$ and $\fix g$ are nonempty.
Thus $f$ and $g$ are two rotations of centers 
$a$ and $b$ respectively, with $a\ne b$.
Below we prove that $\tr[f,g]>2$, yielding $\fix[f,g]=\emptyset$
by Lemma~\reff{l6.3}.a.

Let us start by sending $a$ and $g(a)$ into $i\R$:
Let $\f\in\isomp \H_2$ and $x\in\,]0,+\infty[\,\setminus\{1\}$
be such that $\f(a)=i$ and $\f(g(a))=ix$.
Set $\~f=\f f\f^{-1}$ and  $\~g=\f g\f^{-1}$. We have
$[f,g]=\f^{-1}[\~f,\~g]\f$, hence $\tr[f,g]=\tr[\~f,\~g]$.

Now we write the commutator of $\~f$ and $\~g$ in the form
$[\~f,\~g]=\~fh$, with $h=\~g\~f^{-1}\~g^{-1}$.
The isometry $\~f$ is the rotation of center $i$ and angle $\q\ne0$;
it is associated with the  matrices
$\pm\left(\begin{matrix}\cos t&-\sin t\\
\sin t&\cos t\end{matrix}\right)$
with $t=\tfrac\q2\not\equiv0\!\mod\pi$.
The isometry $h$ is conjugated to $\~f^{-1}$, hence has an angle
$-\q$, and fixes the point $ix$. By Lemma~\reff{l6.3}.\reff d,
it follows that $h$ is associated with the matrices
$\pm\left(\begin{matrix}
\cos t&x\sin t\\
-x^{-1}\sin t&\cos t
\end{matrix}\right)$.
As a consequence, the product $\~fh$ is associated with the matrices
$$
\pm\left(\begin{matrix}
\cos^2 t+x^{-1} \sin^2 t&(x-1)\cos t\sin t\\
(1-x^{-1})\cos t\sin t&x\sin^2 t+\cos^2 t
\end{matrix}\right).
$$
We then have 
$$
\hspace{38mm}
\tr[f,g]=\tr\~fh=2\cos^2 t+(x+x^{-1})\sin^2 t>2.
\hspace{40mm}\square
$$
The proof of Lemma~\ref{l4.2}  is based on the computation of a trace. 
A geometric proof is also available by an adaptation of the proof of
Lemma~\ref{l4.5} in the sequel. 
\prop{p4.2}{
The group $\isom \H_2$ is fixating.
}
\pr
By Corollary~\reff{c4.5}, it is enough to prove that $\isomp \H_2$
is fixating.
Let $G\le\isomp \H_2$ be a \gaf; we have to prove that $G$ is a \gag.
If $G=\{\id\}$ the result is obvious. We now assume $G\ne\{\id\}$.
By Lemma~\reff{l4.2}, all elements of $G\setminus\{\id\}$
are rotations of the same center hence $G$ is a \gag.
\ep

{\noi\sl Remarks.}
\sm\\
1.
Ironically, the M\"obius group acting on $\hat\R$ instead of $\H_2$ is
nonfixating in spite of a lower dimension, 
cf. Exercise~\reff{ex:homographie1}.
\sm\\
2.
As already said, one can also deduce Proposition~\reff{p4.2} from the
forthcoming Corollary~\reff{t4.4} and Exercise~\reff{e12}.
\bigskip

{\noi\bf The case of dimension 3}
\smallskip

We start with the hyperbolic analogue of Lemma~\ref{l3.4} on $\isomp\R^3$.
\lm{l4.5}{
Let $f,g\in\isomp \H_3$. If $\fix f\cap\fix g=\emptyset$ then there exists
$h\in\langle f,g\rangle$ such that $\fix h=\emptyset$.
}
\pr
The proof is done in several steps by proving that one of the 
isometries $f$, $g$, $f^{-1}g$, $gfg^{-1}f^{-1}$, or $gfgf^{-1}$ has no
fixed points.
 
Let $h=f^{-1}g$. If $\fix f$, $\fix g$, or $\fix h$ is empty we are done,
otherwise let $x_0\in\fix h$. We have $f(x_0)=g(x_0)$ and since
$\fix f\cap\fix g$ is empty one has $f(x_0)=g(x_0)\neq x_0$.
Let $P_0=\md(x_0,f(x_0))$; it is a hyperbolic plane by
Lemme~\ref{lem:med}. Since $f$ and $g$ are isometries,
$\fix f$ and $\fix g$ are included in $P_0$.

We claim that $ \fix g \cap f (\fix g) $ is empty.
Otherwise let $ x_1 \in \fix g $ be such that $ y_1 = f (x_1) \in \fix g $.
According to Lemma~\ref{lem:rotationH3}, the sets of fixed points of
elements of $ \isomp \H_3 $ are hyperbolic lines when they 
are nonempty.
Therefore, according to Lemma~\ref{lem:projdroite}
on the projections, there is a single point
$z_1=\pi_{\fix f}x_1$ realizing the distance from $x_1$ to $\fix f$.
The three points $x_1,y_1$, and $z_1$ are distinct
since $\fix f\cap\fix g$ is empty.

Let $\g$ denote the geodesic passing through $x_1$ and $z_1$.
Since $x_1$ and $\fix f$ are in the plane $P_0$,
$\g$ is a hyperbolic line of the plane $P_0$, orthogonal to
$\fix f$ at the point $z_1$ according to Lemma~\ref{lem:projdroite}.
Since $f(x_1)\in\fix g\subset P_0$, the plane $P_0$ is stable by $f$.
Therefore $f(\g)$ is a hyperbolic line of $P_0$.
It is also orthogonal to $\fix f$ at $z_1$ because $f$ preserves angles.
As a result, $f(\g)=\g$ and
the restriction of $f$ to $ \g $ is a symmetry of center $ z_1 $.
Since $ y_1 = f (x_1) \in f (\g) = \g $, $ y_1 $ is
the symmetric on $ \gamma $ of
$ x_1 $ with respect to $ z_1 $.
It follows that $ \g = \fix g $, and therefore $ z_1 $
belongs to $ \fix g $ and $ \fix f $, in contradiction with
$ \fix f \cap \fix g $ empty.

Let $ \~g = fgf^{- 1} $; it is a hyperbolic rotation of the same angle $ \q $
as $ g $ up to the sign (the angle is only defined up to the sign).
We have $ \fix \~g = f (\fix g) $ so
$ \fix \~g \, \cap \, \fix g = \emptyset $ from above.
As before, if $ \fix (g \~g^{- 1}) $ is empty we are done;
otherwise, with $ x_1 \in \fix (g \~g^{- 1}) $, 
the plane $ P = \md (x_1, g (x_1)) $
contains $ \fix g $ and $ \fix \~g $.

We will show that $ \fix (g \~g) $ or $\fix(g\~g^{- 1})$ is empty.

Let $ \sig_P $ be the reflection about the plane $P$.
According to Lemma~\ref{lem:rotationsymetrie}, there is a hyperbolic plane
$ S $ containing the line $ \fix g $ and making an angle $ \pm \tfrac \q2 $ with
$ P $ such that $ g = \sig_S \sig_P $.

Since $ \fix \~g $ is also included in $P$, we can find two hyperbolic 
planes $\~S$ and $\~S'$ containing $\fix\~g$ and making the angles
$\pm\tfrac\q2$ with $P$. According to Lemma~\ref{lem:rotationsymetrie},
$\~g$ or $\~g^{- 1} = \sig_P\sig_{\~\S}$.
It can be assumed without loss of generality that
$\~g= \sig_P \sig_{\~S}$; then we have $\~g^{- 1} = \sig_P\sig_{\~S'}$.

We claim that at least one of the intersections
$ S \cap\~S$ or $S\cap\~S'$ is empty.
Indeed, conjugating by an element of $ \isom \H_3 $ if necessary,
we can assume that $ \fix g $ is a vertical half-line of endpoint
$ a $ in the horizontal plane $ \partial \H_3 $;
$ P $ and $ S $ are then vertical half-planes containing the half-line
$ \fix g$. The boundary of $ P $  is an affine straight line $ \Delta $
of the horizontal plane that contains $ a $
as well as the endpoints $ b $ and $ c $ of the hyperbolic line $ \fix \~ g $
(it may occur
that the hyperbolic line $ \fix\~g$ has only one endpoint;
it is then vertical and the sequel becomes simpler).
Since the hyperbolic lines $ \fix g $ and $ \fix \~ g $ do not cross,
the points $ b $ and $ c $ are on the same side of $ a $ on the line $ \Delta $,
so one of the points, say $ b $, is between $ a $ and $ c $.
Let $ V $ be the vertical plane containing $ b $ and parallel to the plane $ S $.
Since the angles of the hyperbolic planes $\~S $ and $\~S' $ with the plane
$ P $ are $ \pm \tfrac \q2 $, one of the planes $\~S$ or $\~S'$
is tangent to $ V $.
Suppose it is $\~S$;
if this is not the case simply replace $ \~ g $ by $ \~ g^{- 1} $.
The plane $\~S $ does not cross $ S $ because it is located on one side of $ V $ 
and its closure contains $ c $ while $ S $ contains $ a $.
We have $g\~g=\sig_S\sig_P\sig_P\sig_{\~S}=\sig_S\sig_{\~S}$.

We conclude by showing that $ \fix (\sig_S \sig_{\~S}) $ is empty.
Assume by contradiction that there exists
$ x \in \fix (\sig_S \sig_{\~S}) $ hence
$ \sig_S (x) = \sig_{\~ S} (x) $. Since
$ \fix \sig_{\~S} \cap \fix \sig_S = S\cap\~S = \emptyset$,
we have $ \sig_S(x)\ne x$
hence the hyperbolic planes $S$ and $\~S$ are both contained in
the plane $\md(x,\sig_S(x)) $, so coincide, in contradiction with
$S\cap\~S = \emptyset$.
\ep
\lm{l4.6}{
Let $G\le\isomp\H_3$ be such that $\fix f\cap\fix g\ne\emptyset$ for all
$f,g\in G$. Then $G$ is a \gag.
}
\pr
If $ \fix f = \fix g $ for all $ f, g \in G \setminus \{\id \} $ we are done.
Otherwise let $ f, g \in G \setminus \{\id \} $ with $ \fix f \ne \fix g $,
let $ \Pi $ be the hyperbolic plane containing the
hyperbolic  lines $ \fix f $ and $ \fix g $, and let $ \omega \in \Pi $ be
the point of intersection of these lines.
It remains to prove that $ \omega \in \fix G$.
Assume by contradiction that there exists $ h \in G $ 
such that $ \omega \notin \fix h $.
Let $ \a $ be such that $ \fix f \cap \fix h = \{\a \} $ and $ \beta $
such that $ \fix g \cap \fix h = \{\beta \} $.
Therefore we have $ \a \ne \beta $ and $ \a, \beta \in \Pi \cap \fix h $
hence $ \fix h \subset \Pi $.
By the way, we have shown
\eq7{
\forall k\in G\quad(\omega\notin\fix k\Rightarrow\fix k\subset\Pi).
}
With the same $ h $, let
$ \delta \in \fix (fh) \setminus \{\a \} $.
We have $ \delta \notin \fix h $ (otherwise $ \delta \in \fix f \cap \fix h $,
but $ \delta \ne \a $) hence $ \fix h \subset \md (\delta, h (\delta)) $.
Similarly
$ \fix f = \fix f^{- 1} \subset \md (\delta, f^{- 1} (\delta)) = \md (\delta, h (\delta)) $,
which yields $ \md (\delta, h (\delta)) = \Pi $. It follows that
$ \delta \notin \Pi $, hence $ \fix (fh) \not \subset \Pi $.
Now the contraposition of~\rf7 implies
$ \omega \in \fix (fh) $ but this implies $ \omega \in \fix h $,
a contradiction.
\ep
\cor{t4.4}{
The group $\isom\H_3$ is fixating.
}
\pr
Since $ \isomp \H_3 $ is a subgroup of index $ 2 $ of $ \isom \H_3 $,
by Corollary~\reff {c4.5} it suffices to prove that $ \isomp \H_3 $
is fixating.
Let $ G \le \isomp \H_3 $ be a \gaf. According to Lemma~\reff{l4.5}, 
we have
$ \fix f \cap \fix g \ne \emptyset $ for all $ f, g \in G $ hence $ G $
is a \gag\, by Lemma~\reff{l4.6}.
\ep
\sub{section:globalisantelliptique}{The spherical\index{spherical space} case}
The sphere $ \S_n $ is endowed with the spherical distance 
$ d (x, y) = \arccos\scal xy $, where $\scal\;{\;\,}$ is the usual scalar product.
Since the function arccos is bijective, the isometries of $ \S_n $ for
the spherical distance coincide with the isometries of $ \S_n $ for
the distance induced by the Euclidean distance of $ \R^{n + 1} $.
It is known that an isometry of $ \S_n $ is the restriction to $ \S_n $
of a vector isometry of $ \R^{n + 1} $, see~\cite b, Chap.18 or
Exercise~\reff{e5.3}.
More precisely, if $ \OO_{n + 1} $ denotes the group of vector isometries
of $ \R^{n + 1} $ then we have a group isomorphism
$$
\f:\OO_{n+1}\to\isom\S_n
$$
which maps any element $f\in\OO_{n + 1}$ to its restriction to $\S_n$.
We have $ \fix f = \{\v0 \} $ if and only if $ \fix \f (f) = \emptyset $.
By abuse of language, we will say that a subgroup $ G $ of $ \OO_{n + 1} $
is a \gaf, resp. a \gag, resp. eccentric, if its image $ \f (G) $ is
a \gaf, resp. a \gag, resp. an eccentric subgroup of $ \isom \S_n $.
In other words we consider $\OO_{n+1}$ as acting on
$\R^{n+1}\setminus\{\v0\}$.
\th{t5.7}{
The group $\isom\S_n$ is fixating if and only if $n=1$.
\\
The group $\isomp\S_n$ is fixating if and only if $n=1$ or $n=3$.
}
Observe that Corollary~\reff{c4.5} does not apply here because 
the median inequality does not hold in elliptic spaces.
The proof of Theorem~\reff{t5.7} is split into several parts.

For $ n = 1 $, the only \gaf\, of $ \isom \S_1 $ are the trivial
group $\{\id\}$ and the groups $ \{\id, s \} $ where $ s $ is any fixed
reflection along a straight line, which are obviously \gag.
\prop{p5.7}{
If $n\ge2$ then $\isom\S_n$ is nonfixating.
}
\pr
We use the framework of $ \OO_{n + 1} $.
Let us start with $ n = 2$. Let $ (\v i, \v j, \v k) $ be 
an orthonormal basis of $ \R^3 $. Let $ f, g, h \in \OO_3 $
be the vector isometries of matrices respectively 
$ \mat f = \diag (1, -1, -1) $, $ \mat g = \diag (-1,1, -1) $
and $ \mat h = \diag (-1, -1,1) $.
The set $ G = \{\id, f, g, h \} $ is a group isomorphic to the Klein group
$ (\Z / 2 \Z)^2 $: We have $ f^2 = g^2 = h^2 = \id $,
$ fg = gf = h $, $ gh = hg = f $ and $ hf = fh = g $.
We therefore have $ \fix k \ne \{\v0 \} $ for all
$ k \in G $ but $ \fix G = \{\v0 \} $.

For $ n \ge3 $, we complete the matrices of $ f $ and $ g $ by $ -1 $
and the matrix of $ h $ by $ 1 $: In an orthonormal basis, one chooses
$ \mat f = \diag (1,-1,\dots,-1) $, $ \mat g = \diag (-1,1, -1,\dots,- 1)$,
and $ \mat h = \diag (-1, -1,1,1, \dots, 1) $.
We check that $ \{\id, f, g, h \} $ is still isomorphic to the Klein group
and that it is eccentric.
\ep

Let us now treat the case of $\isomp\S_n$.
The isomorphism $ \f: \OO_{n + 1} \to \isom \S_n $ induces by restriction an
isomorphism from $ \SO_{n + 1} $ into $ \isomp \S_n $ which will 
also be denoted by~$\f$.
\prop{p5.8}{
If $n$ is even then $\isomp\S_n$ is eccentric.
}
\pr
Let $ g \in \isomp \S_n $ and let $ f = \f^{- 1} (g) \in \SO_{n + 1} $.
All the eigenvalues of $ f $ are of modulus 1 and their product is
equal to 1. Moreover, if $ \l $ is an eigenvalue of $ f $ then $ \bar \l $ too,
with the same multiplicity.
Since $ n + 1 $ is odd, we deduce that 1 is an eigenvalue of $ f $
hence $ \fix f \ne \{\v0 \} $. It follows that $ \isomp \S_n $ is a
\gaf\;but is obviously not a \gag.
\ep

It remains to treat the case when $n$ is odd.
For $ n = 1 $, $ \isomp \S_1 $ is obviously fixating.
We now treat the case $n\ge5$ and we will end up with $n=3$.
\prop{p5.10}{
If $n$ is odd and $n\ge5$ then $\isomp\S_n$ is nonfixating.
}
\pr
%
For $ n = 5 $, let $ f, g, h \in \SO_6 $ be the isometries of 
block diagonal matrices $ \mat f = \diag (I_2, -I_2, -I_2) $,
$ \mat g = \diag (-I_2, I_2, -I_2) $ and $ \mat h = \diag (-I_2, -I_2, I_2) $
 respectively.
Similarly to the proof of Proposition~\reff{p5.7}, 
we verify that $ \{\id, f, g, h \} $ is isomorphic to
the Klein group and eccentric.

For $ n $ odd, $ n \ge7 $, we expand the matrices of $ f $ and $ g $ by
$-I_2$ and the matrix of $h$ by $I_2$ as in the proof
of Proposition~\reff{p5.7}. The group $\{\id,f,g,h\}$ is still
isomorphic to the Klein group and eccentric.
\ep
\prop{p5.9}{
The group $\isomp\S_3$ is fixating.
}
\pr
By contradiction, let $G\le\SO_4$ be an eccentric group.
An element of $G\setminus\{\id\}$ has a plane
(i.e. a subspace of dimension 2) of fixed points.
\med\\
{\noi\sl Step 1. \,One has $\fix f\cap\fix g\ne\{\v0\}$
for all $f,g\in G$.
}
\sm

Let $ \v u \in \fix (f^{- 1} g) \setminus \{\v0 \} $ (which is nonempty, as
$G$ is a \gaf). If $ \v u \in \fix f \cap \fix g $ we are done.
Otherwise we have $ f (\v u) = g (\v u) \ne \v u $.
Then $ \fix f $ and $ \fix g $ are two planes included in
the hyperplane (i.e. three dimensional subspace)
$ \md (\v u, f (\v u)) $ hence intersect each other.
\eep

Now set $ f_0 \in G \setminus \{\id \} $.
Since $ \fix G = \{\v0 \} $, there exists $ g_0 \in G \setminus \{\id \} $ such that
$ \fix f_0 \ne \fix g_0$. From Step 1, $ \fix f_0 \cap \fix  g_0$
is a straight line denoted by $D$ and $ \fix f_0 + \fix g_0 $ is a
hyperplane denoted by $H$.
\med\\
{\noi\sl Step 2. \,
For all $f\in G\setminus\{\id\}$, one has $\fix f\subset H$.
}
\sm

Indeed, let $ h_0 \in G $ be such that $ \fix h_0 $ does not contain $ D $.
Such an element $ h_0 $ exists since $ \fix G = \{\v0 \} $.
We have $ \fix f_0 \ne \fix g_0 \ne \fix h_0 $.
From Step 1, $ \fix f_0 \cap \fix h_0 $ is a straight line denoted by $D'$
and $\fix g_0\cap\fix h_0$ is a straight line denoted by $D''$.
We have $ D'\ne D''$ (otherwise $D' = D'' = D$, in contradiction with
$D\not\subset \fix h_0 $) and $ \fix h_0 $ is two dimensional
hence $ \fix h_0 = D '+ D'' \subset H $.
Also notice that the three straight lines $ D $, $ D '$,
and $D'$ are not coplanar (otherwise we would have
$\fix f_0=\fix g_0=\fix h_0$)
hence $H = D + D'+ D''$, so that $ f_0 $, $ g_0 $ and $ h_0 $ play
symmetrical roles.
Now let $ f \in G \setminus \{\id \} $. Therefore $ \fix f $, which
is a plane, cannot contain at the same time $ D $, $ D' $ and $ D''$
hence, from the above, $ \fix f $ is included in
one of the three subspaces $ \fix f_0 + \fix g_0 $,
$ \fix g_0 + \fix h_0 $ or $ \fix h_0 + \fix f_0 $,
all of which in fact coincide with $H$.
\eep

Let $\D$ denote the line orthogonal to $H$: $\D=H^\perp$.
\med\\
{\noi\sl Step 3. \,
One has $f(\D)=\D$ for all $f\in G$.
}
\sm

Otherwise let $ f \in G $ be such that $ f (\D) \ne \D $ and let
$ g \in G \setminus \{\id \} $.
Let us show that $ (\fix g)^\perp = \D + f (\D) $.
We have $ \D \subset (\fix h)^\perp $ for all 
$ h \in G \setminus \{\id \} $ and,
since $ f $ is an isometry, $ f (\D) \subset (f (\fix h))^\perp $.
For $ h = f^{- 1} gf $, this gives 
$ f (\D) \subset (f (\fix h))^\perp = (\fix g)^\perp $.
Therefore $ \D + f (\D) \subset (\fix g)^\perp $, which have the same 
dimension hence are equal.
Since this holds for all $ g \in G \setminus \{\id \} $,
we get $ \fix G = (\D + f (\D))^\perp \ne \{\v0 \} $, a contradiction.
\eep

To sum up we found a hyperplane $ H $ such that 
$ \fix f \subset H $ and $ f (H^\perp) = H^\perp $
for all $ f \in G \setminus \{\id \} $.
It follows that the restriction of $ f $ to $ H^\perp $ is $ - \id $,
thus the only possible eigenvalues of $f$ are $1$ and $-1$, hence $f^2=\id$
for all $f\in G$, so $G$ is Abelian (we have $\id = f^2g^2 = (fg)^2$
hence $fg = gf$ for all $f, g \in G$).

Since $ G $ is Abelian, the elements of $G$ are diagonalizable in a
common basis, denoted by $\mathcal B=(\v i,\v j,\v k,\v\ell)$,
with $1$ and $-1$ as double eigenvalues
(these are positive isometries of $ \R^4 $).
Consider the endomorphisms $f_1$, $f_2$ and $f_3$ whose matrices
in the basis $\mathcal B$ are respectively
$ \diag (1,1, -1, -1) $, $ \diag (1, -1,1, -1) $, and $ \diag (1, -1, -1,1) $.
Thus $G$ is a subgroup of 
$G_0=\{\pm\id,\pm f_1,\pm f_2,\pm f_3\}$.
\med

The list of the sixteen subgroups of $G_0$ splits into
\med

\bul eleven \gag:  
$\0$, $\{\id,f_n\}$ with $1\le n\le3$, $\{\id,-f_n\}$ with $1\le n\le3$,
\noi
$\{\id,f_1,f_2,f_3\}$ which fixes $\v i$, 
$\{\id,f_1,-f_2,-f_3\}$ which fixes $\v j$,
\noi
$\{\id,-f_1,f_2,-f_3\}$ which fixes $\v k$,
$\{\id,-f_1,-f_2,f_3\}$ which fixes $\v\ell$,
\med

\bul
and five containing $-\id$ hence not \gaf:
\noi
$\{\id,-\id\}$, $\{\id,f_n,-\id,-f_n\}$ with $1\le n\le3$,
and~$G_0$. 
\med

As a consequence $G_0$ contains no eccentric subgroup, a contradiction.
\ep
\sub{5.4}{The projective case}\index{projective space}
Usually $\RP_n$ is the set of vector lines of $\R^{n + 1}$.
In this article we identify $\RP_n$ with the quotient of $\S_n$
by the equivalence relation
$$
x\sim y\ssi x=y~\mbox{ or }~x=-y.
$$
For $x\in\S_n$, we write $\dot x = \{x,-x\}$ the corresponding class
in $\RP_n$.
Given $ \dot x = \{x, -x \} $ and $ \dot y = \{y, -y \} $ in $ \RP_n $,
the distance between $ \dot x $ and $ \dot y $ is then given by
$ d (\dot x,\dot y) = \arccos|\scal xy|$ where $\scal xy$ is the
scalar product between $x$ and $y$.

Given an isometry $f$ of $\S_n$ and $x\in\S_n$, the class of $f(x)$ in
$\RP_n$ is the same as $f(-x)=-f(x)$, so we can define a
function from $\RP_n$ to $\RP_n$, denoted by $\psi(f)$, which maps
$\dot x=\{x,-x\}$ to the class of $f(x)$.
It is known that the mapping
$$
\psi:\isom\S_n\to\isom\RP_n
$$
defined in this manner is a surjective morphism of kernel $\{\pm\id\}$,
cf.~\cite b, Chap.19.
Thus the group $\isom\RP_n$ is the image of $\isom\S_n$
by $\psi$ and similarly $\isomp\RP_n$ is the image of $\isomp\S_n$
by $\psi$. When $n$ is even, $-\id:\S_n\to\S_n$ reverses the orientation,
so $\isomp\RP_n$ is equal to $\isom\RP_n$.
\th{t5.12}{
The group $ \isom\RP_n $ is fixating if and only if $ n = 1 $.
\\
The group $ \isomp\RP_n $ is fixating if $ n = 1 $ and nonfixating if
$ n $ is odd and greater than or equal to $ 5 $.
}
{\noi\sl Remarks.}
\sm\\
1.
We do not know whether $ \isomp \RP_3 $ is fixating or not.
\sm\\
2.
Notice that the rotations of $ \RP_2 $ have a single fixed point.
Thus the rotation group $ \RP_2 $ is another example of a nonfixating
group which is not Abelian and such that all elements
have a single fixed point.
\med\\
\pr
It is easy to check that $ \isom \RP_1 $ is fixating.

If $ n $ is even and $ \dot f = \psi (f) \in \isom \RP_n $ then the matrix
of $ f $, seen as the vector isometry of $ \R^{n + 1} $,
is odd-sized, so always admits $ 1 $ or $ -1 $
as eigenvalue, hence there exists $ x \in \S_n $ such that $ f (x) \in \{x, -x \} $.
Thus any isometry of $ \RP_n $ has at least one fixed point $ \dot x $
but no point of $ \RP_n $ can be fixed by all the elements of
$ \isom\RP_n$. The group $ \isom \RP_n $ is therefore eccentric itself
hence nonfixating.

We now study the case when $n$ is odd. For $n=5$, set
$$
I=\left(\begin{matrix}1&0\\0&1\end{matrix}\right)~\mbox{ and }~
R=\left(\begin{matrix}0&-1\\1&0\end{matrix}\right),
$$
and consider the isometries $ f, g, h $  of $ \S_5 $ whose matrices
are the following block diagonal matrices
$$
\mat f=\diag(I,R,R),\;\mat g=\diag(R,I,R)~\mbox{ and }~\mat h=\diag(R,R,I).
$$
These are positive isometries.
Let $ G_5 $ denote the subgroup of $ \isomp \S_5 $ generated by $ f, g $ and $ h $.
One finds for $ G_5 $ the following group of order $ 32 $:
$$
G_5=\{\diag(\pm I,\pm I,\pm I),\diag(\pm I,\pm R,\pm R),
\diag(\pm R,\pm I,\pm R),\diag(\pm R,\pm R,\pm I)\}
$$
and ones checks that its image $ \psi (G_5) $ in $ \isomp \RP_5 $ is eccentric.
The verification is somewhat tedious but without difficulty.
It follows that $ \isomp \RP_5 $ is nonfixating hence so is $\isom\RP_5$.

For $ n $ odd, $ n \ge7 $, we consider the positive isometries
$ f, g, h $ of $ \S_n $ whose matrices are the block diagonal matrices:
$ \mat f = \diag (I, R,\dots, R) $,
$ \mat g = \diag (R, I, R,\dots, R) $
and $ \mat h = \diag (R, R, I, R, \dots, R) $.
We verify similarly that the  image $ \psi (G_n) $ of the subgroup $ G_n $ of
$ \isomp \S_n $ generated by $ f, g $, and $ h $ is eccentric,
so neither $ \isomp \RP_n $
nor $ \isom \RP_n $ are fixating.

It remains to treat the case $n=3$.
Let $f,g\in\isom\S_3$ of matrices
$$
\mat f=\left(\begin{matrix}
0&0&0&1\smallskip\\
1&0&0&0\smallskip\\
0&1&0&0\smallskip\\
0&0&1&0\end{matrix}\right)
~\mbox{ and }~
\mat g=\left(\begin{matrix}
0&0&0&-1\smallskip\\
1&0&0&0\smallskip\\
0&-1&0&0\smallskip\\
0&0&1&0\end{matrix}\right).
$$
For simplicity, we write $ f $ and $ g $ as signed permutations:
$ f = (1 \, 2 \, 3 \, 4) $ and $ g = (1 \, 2 \, - 3 \, - 4) $.
The computation gives $ fg = (1 \, 3) (2 \, - 4) = - gf $ and
$ f^2 = (1 \, 3) (2 \, 4) = - g^2 $.
Therefore the group generated by $ \pm f $ and $ \pm g $ is
$$
G_3 = \{\pm \id, \pm f, \pm f^2, \pm f^3, \pm g, \pm fg, \pm f^2g, \pm f^3g \}.
$$
By checking one by one all the elements of $G_3$, we easily verify that
$\psi(G_3)$ is eccentric in $\isom\RP_3$ hence $\isom\RP_3$ is nonfixating.
We did not find any eccentric subgroup of
$ \isomp\RP_3 $, nor were we able to adapt to the projective case the proof that
$ \isomp\S_3 $ is fixating (Proposition~\reff {p5.9}).
\ep
%
%
\sec{5.}{Groups acting on discrete sets}
\sub{5.2}{Symmetric groups}\index{symmetric group} 
The aim of this section is the following result.
\th{t6.1}{
\be\item
The symmetric group $\sy_n$ acting on $\{1,\dots,n\}$ is fixating
if and only if $n\le4$.
\item
The alternating\index{alternating group}  group $\A_n$ of even permutations 
acting on $\{1,\dots,n\}$ is fixating if and only if $n\le4$.
\ee
}
\pr
{\bf a.}
We split the proof for the symmetric group into five steps.
We use the decomposition in cycles with disjoint supports.
A cycle of order $n$ is called an {\sl$n$-cycle}.
\med\\
{\sl Step 1. $\sy_n$ is fixating when $n\le3$.} 
\sm

For $n=1$ and $n=2$, there is nothing to prove since there is no
nontrivial subgroup.
For $ n = 3 $, the nontrivial subgroups of $ \sy_3 $ are
$\mathcal{A}_3 = \langle (123) \rangle $, which is not a \gaf,
and the three subgroups of order 2 generated by each transposition,
which are \gag, hence $\sy_3$ has no eccentric subgroup.
\med\\
{\noi\sl Step 2. $\sy_4$ is fixating.}
\sm

Let $ G \le \sy_4 $ be a \gaf. So $ G $  contains neither a
double-transposition nor a 4-cycle
(since these permutations have no fixed points), so $G\setminus\{\id\}$
contains only transpositions and/or 3-cycles.
Two cases occur.
\begin{itemize}
\item[$\triangleright$]
If $ G $  contains no 3-cycle then $ G = \{\id, \tau \} $ where
$ \tau $ is a transposition and we are done.
\item[$\triangleright$]
If $ G $ contains a 3-cycle $ \g $, say $ \g = (123) $,
let us show by contradiction that
$ \fix G = \{4 \} $.
Otherwise there exists $ g \in G $ such that $ g(4) \ne 4 $, say $ g(4) = 1 $.
If $ g = (14) $ then $ g \g = (1234) $ which has no fixed points.
If $ g $ is a 3-cycle, say $ g = (124) $, then $ g \g = (14) (23) $
which has no fixed points either. 
\end{itemize}

{\noi\sl Step 3. $\sy_5$ is nonfixating.}
\sm

Let $ G = \langle f, g \rangle $ with $ f = (123) $ and $ g = (12) (45) $.
We verify that $ G = \{\id, f, f^2, g, h, k \} $ with $ h = fg = (13) (45) $
and $ k = gf = (23) (45) $, and that $ G $ is eccentric.
\sm\\
{\noi\sl Remark.}
We can see $ G $ as $ (\sy_3 \times \sy_2)^+ $, the set of even permutations
acting separately on $ \{1,2,3 \} $ and on $ \{4,5 \} $.
We also can interpret $ G $ as the group of isometries of the
``double tetrahedron'', i.e. the hexahedron obtained by gluing two regular
isometric tetrahedra on one of their faces.
\med\\
{\sl Step 4. $\sy_6$ is nonfixating.}
\sm

Let $ G = \langle f, g \rangle $ with $ f = (12) (34) $ and $ g = (12) (56) $.
We verify that $ G = \{\id, f, g, h \} $ with $ h = (34) (56) $ and
that $ G $ is eccentric.
\eep\\
{\noi\sl Remark.}
Again $ G $ can be interpreted as a set of even permutations:
those acting separately on $ \{1,2 \} $, on $ \{3,4 \} $, and on $ \{5,6 \} $,
and also as a group of isometries: the half-turns of axes the coordinate
axes in $ \R^ 3 $.
\med\\
{\sl Step 5. $\sy_n$ is nonfixating when $n\ge7$.} 
\sm

Let $G=\langle f,g\rangle$ with $f=(123)(6\dots n)$ and
$g=(12)(45)(6\dots n)$. We check that $ G $ is eccentric.
\med\\
{\bf b.}
For $n\le4$, $\A_n$ is fixating, as a subgroup of a fixating group.
For $n=5$ and $n=6$, both eccentric subgroups built in steps 3 and 4 
are precisely subgroups of $\A_5$, resp. $\A_6$,
which shows that $\A_5$ and $\A_6$ are nonfixating.
For $n\ge8$, $n$ even, the cycle $(6\dots n)$ is even, so the eccentric
subgroup of Step 5 is still in $\A_n$.
For $n\ge9$, $n$ odd, the cycle $(7\dots n)$ is even, so we
complete as in Step 5, but with the group of Step 4~:
We choose $G=\langle f,g\rangle$ with $f=(12)(34)(7\dots n)$ and
$g=(12)(56)(7\dots n)$. We verify that $G$ is an eccentric subgroup of
$\A_n$.
We leave to the reader the pleasure to explore the most interesting case
$\A_7$, see Exercise~\reff{e15}.
\ep

{\noi\sl Remark.}
Theorem~\ref{t6.1} shows that the groups $\sy_n$ and $\A_n$ are
nonfixating when $n>4$. Some of their subgroups however may be fixating. 
For example, given $d\ge1$ and $q=p^n$ integers with $p$ prime,
let $G=\GL(d,\F_q)$ denote the group of isomorphisms
of the $d$-dimensional vector space on the $q$-element finite field $\F_q$.
Then $G$ acts naturally on $X=\F_q^d\setminus\{\v0\}$ and can be identified
with a subgroup of $\sy_X$, the group of permutations of $X$. 
One can then ask in which cases $G$ is fixating.
The answer is that only $\GL(1,\F_q)$ and $\GL(2,\F_2)$ are fixating,
see Exercises~\reff{e16} and~\reff{e17}.
We did not study in detail the case of the affine group $G=\GA(d,\F_q)$
acting on $\F_q^d$.
\sub{6.2}{Isometries of $\Z^n$}
Let $ n \ge 1 $ be an integer. We equip $ \Z^n $ with the Euclidean norm
$ \| \;\| $. Any $f\in\isom \Z^n$ can be extended
into an isometry of $ \R^n $, denoted by the same notation~$f$,
cf.~Exercise~\reff{e54} and its solution.
\th{t6.2}{
The group $\isom\Z^n$ is fixating.
}
\pr
By contradiction, let $G$ be an eccentric subgroup of $\isom\Z^n$
with $n$ minimal.
\vspace{0.2cm}\\
{\sl Step 1. $G$ is finite.}
\sm\\
Let $\mathcal{C}_n=\{x=(x_1,\dots,x_n)\tq x_i=-1,\,0\mbox{ or }1\}$ be
the unit hypercube in $ \Z^n $. Since the only translation of $G$ is $\id$,
the morphism from $ G $ to $ \isom \mathcal{C}_n $ which maps any isometry $f$
to its linear part $ \v f$ is injective. Now $ \isom \mathcal{C}_n $
is finite, so $ G $ is finite too.
\eep

Denote by $N$ the cardinality of $G$ and by $\w$ the centroid of
the orbit of $\v0$ by $G$, i.e. $\w=\frac1N\sum_{f \in G}f(\v0)$.
\med\\
{\sl Step 2. All coordinates of $\w$ are
congruent to $\frac12$ modulo $1$.}
\sm\\
Indeed, $\w$ is fixed by all the elements of $G$.
Since $G$ is eccentric, $\w\notin \Z^n$.
Let
$$
\W=\{x\in\Z^n\tq \|\w-x\|=d(\w,\Z^n)\}.
$$
Let $ I = \big \{i \in \{1, \dots, n \} \tq \exists x, y \in \W, \; x_i \ne y_i \big \} $.
We have $ \w_i \equiv \frac12 \modd1 \ssi i \in I $.
Indeed, if $ x, y \in \W $ are such that $ x_i \ne y_i $ then $ | x_i-y_i | = 1 $
and $ | x_i- \w_i | = | y_i- \w_i | = \frac12 $.
Conversely, if $ \w_i \equiv \tfrac12 \modd1 $ and $ x \in \W $ then
$ x_i = \w_i \pm \tfrac12 $ and the  point $y$ with the same
coordinates as $ x $ except the $ i $-th  coordinate
equal to $ \w_i \mp \tfrac12 $ is also in $\W$.
In summary we have
\[
\W=\big\{(x_1,\dots,x_n)\in \Z^n\tq x_i=\w_i\pm\tfrac12\text{ if }i\in I
\text{ and } x_i=\lfloor\w_i\rceil \text{ if } i\notin I\big\},
\]
where $\lfloor\w_i\rceil$ is the unique integer such that 
$|\,\w_i-\lfloor\w_i\rceil|<\frac12$ (recall that
$\w_i\not\equiv \frac12\modd1$ when $i\notin I$).
Let $ E = \Z^n \cap \aff \W $.
It is a `` lattice'' isometric to $ \Z^k $, where $ k $ is the cardinality of $ I $.
Precisely, let us set $ E_i = \Z $ if $ i \in I $ and $ E_i = \{\lfloor\w_i\rceil \} $ otherwise;
then we have $ E = E_1 \times \cdots \times E_n $.
For all $ f \in G $ we have $ f (\W) = \W $, so $ f(\aff \W) = \aff \W $, and in addition
$ f(\Z^n) = \Z^n $, so $ f(E) = E $. This allows to define
$ G_E $, the set of all restrictions to $ E $ of the elements in $ G $.
These are  isometries of $E$. Since $ G $ is not a \gag, $ G_E $ is not a \gag\;either.
Moreover, for all $ f \in G $ and all $ x \in \fix f $, the orthogonal projection
$ \pi_Ex $ is also in $ \fix f $, so that $ G_E $ is a \gaf.
By minimality of $ n $ we deduce that $ k = n $ hence $ I = \{1, \dots, n \} $.
\med\\
{\sl Step 3. One is reduced to $\w=\v0$ and one changes
$\Z^n$ into $(2\Z+1)^n$.}
\sm\\
Let $ \f: \Z^n \to (2 \Z + 1)^n$ be the map defined by $x \mapsto 2x-2 \w $ and let $ f \in G $. Then the isometry
$ \~f = \f f \f^{- 1} $ fixes $ \v0 $ and maps $ (2 \Z + 1)^n $ into $ (2 \Z + 1)^n $.
In addition, a small calculation shows that $ \~f = \v f $, the linear part
associated with $ f $, so $ \f $ also globally fixes the lattice $ \Z^n $,
and the entries of its matrix are only $ 0 $, $ 1 $ or~$-1$.
\eep

To avoid multiple notations, we still denote by $G$ the conjugate
of $G$ by $\f$.
For each $ f\in G $, the matrix of $f$, 
denoted $\mat f=(a_{i,j})_{1\le i,j\le n}$, is thus a matrix of
a signed permutation: There is
one and only one nonzero entry on each row and each column
and this entry is equal to $1$ or~$-1$.
\med\\
{\sl Step 4. The diagonal coefficients of $\mat f$ are never $-1$.}
\sm\\
Since $G$ is a \gaf, $\fix f$ is a nonempty subset of $(2\Z+1)^n$.
Let $x=(x_1,\dots,x_n)\in\fix f$.
It is an element of $(2\Z+1)^n$, hence $x_i\ne0$ for all $i=1,\dots,n$.
If $a_{i,i}\ne 0$, the $i$-th coordinate of the equality $f(x)=x$ gives
$a_{i,i}x_i=x_i$ hence $a_{i,i}=1$.
\eep

Denote by $\sigma:G\to\sy_n,f\mapsto\sigma_f$ the function which maps $f$ of
matrix $(a_{i,j})_{1\le i,j\le n}$ to the permutation matrix
$ (|a_{i, j}|)_{1\le i,j\le1}$.
This is clearly a group homomorphism which is injective by Step 4
(its kernel is reduced to $\id$).
\med\\
{\sl Step 5. Towards  the construction of a global fixed point.}
\sm\\
Let $(e_1,\dots,e_n)$ denote the canonical basis of $\R^n$ and 
consider  the two relations on $\{1,\dots,n\}$:

\noi $i\sim j$ if there exists $f\in G$ such that $f(e_i)\in\{e_j,-e_j\}$,

\noi $i\approx j$ if there exists  $f\in G$ such that $f(e_i)=e_j$.

It is easy to check that they are equivalence relations.
Let $ p $ be the number of classes for the relation $ \sim $.
For each $ k = 1, \dots, p $, the class $ C_k $  for $ \sim $ is
partitioned into two classes for $ \approx $
(possibly, one of the classes is empty).
Let $ C_k^+ $ be one of these classes (arbitrarily chosen)
and $ C_k^- = C_k \setminus C_k^+ $. Denote
$$
C^+= C_1^+\cup\dots\cup C_p^+\quad\mbox{and}\quad C^-= C_1^-\cup\dots\cup C_p^-.
$$
According to Step 4, given $ i \sim j $, the equality $ f (e_i) = e_j $
is not possible for one  $ f \in G $ and $ g (e_i) = - e_j $ for another,
so we have $i\approx\sigma_f(i)\ssi f(e_i)=e_{\sigma_f(i)}$.

Let $f\in G$ be fixed and denote
\begin{align*}
C^{++}_f=C^+\cap\sigma_f^{-1}(C^+),&\quad C^{+-}_f=C^+\cap\sigma_f^{-1}(C^-),\\
C^{-+}_f=C^-\cap\sigma_f^{-1}(C^+),&\quad C^{--}_f=C^-\cap\sigma_f^{-1}(C^-).
\end{align*}
Thus we have $f(e_i)=e_{\sigma_f(i)}$ when $i\in C^{++}_f\cup C^{--}_f$ and
$f(e_i)=-e_{\sigma_f(i)}$ when $i\in C^{+-}_f\cup C^{-+}_f$.
We also have 
$$
C^+=C^{++}_f\cup C^{+-}_f,\quad C^-=C^{-+}_f\cup C^{--}_f
$$
and 
$$
\sigma_f(i)\in C^+\ssi i\in C^{++}_f\cup C^{-+}_f,\quad
\sigma_f(i)\in C^-\ssi i\in C^{+-}_f\cup C^{--}_f.
$$
{\sl Step 6. The point  $x=(x_1,\dots,x_n)$ defined by $x_i=1$ if 
$i\in C^+$ and $x_i=-1$ if $i\in C^-$ is fixed by all elements of $G$.}
\sm\\
Indeed, for any $f\in G$, we have 
\begin{align*}
f(x)&=f\Big(\sum_{i\in C^+}e_i\Big)-f\Big(\sum_{i\in C^-}e_i\Big)\\
&=
\bigg(\sum_{i\in C^{++}_f}e_{\sigma_f(i)}-\sum_{i\in C^{+-}_f}e_{\sigma_f(i)}\bigg)
-\bigg(\sum_{i\in C^{--}_f}e_{\sigma_f(i)}-\sum_{i\in C^{-+}_f}e_{\sigma_f(i)}\bigg)\\
&=\sum_{\sigma_f(i)\in C^+}e_{\sigma_f(i)}-\sum_{\sigma_f(i)\in C^-}e_{\sigma_f(i)}\\
&=\sum_{j\in C^+}e_j-\sum_{j\in C^-}e_j=x.
\end{align*}
It follows that $\fix G$ is nonempty, a contradiction.
\ep
\sub{5.1}{Isometries of trees\index{tree}}
All the results of this section are taken from 
J.-P.~Serre's book~\cite s\index{Serre, Jean-Pierre}.
A {\sl combinatorial tree} $X$ is a simple connected undirected graph
without cycles. The following fundamental  property is easy to check: 
Given two vertices $P$ and $Q$, there exists a unique injective path
joining them. We denote by $[P,Q]$ this path and by $d(P,Q)$ its length,
i.e. the number of its edges.
The map $d:X\times X\to\N$ is the {\sl combinatorial distance} on $X$
and $(X,d)$ is a discrete metric space.
Two vertices are joined by an edge if and only if their distance 
is $1$. We recall that an isometry of $X$ is a bijection of $X$ that
preserves the distances.
\lm{l7.1}{
Let $s$ be an isometry with at least one fixed point.
Then, for every vertex $x$ of $X$, the	distance $d(x,s(x))$ is even and 
the midpoint $z$ of the path $[x,s(x)]$ is a	fixed point of $s$.
Moreover, $z$ is the unique element of $\fix s$ such that
$d(x,z)=d(x,\fix s)$.
}
\pr
The uniqueness of the geodesic path between two vertices and the fact
that the image
of an injective path  by the isometry $s$ is an injective path with the
same length ensure that $\fix s$ is connected hence is a subtree of $X$.

If $x\in \fix s$ then the statement is clear.
Suppose that $x\notin\fix s$. Then there exists $b\in\fix s$
such that $n=d(x,\fix s)=d(x,b)\ge1$. If $b$ were not unique,
one could construct a nontrivial cycle thanks to the  connectedness of $\fix s$.
Let $[x_{0},x_{1,}x_{2},\dots,x_{n}]$ be the geodesic path  joining $b=x_0$
to $x=x_n$. 
By definition of $b$, $s(x_i)\neq x_i$ for all $i>0$.
The two paths $[x,b]=[x_n,x_{n-1},\dots,x_0]$ and
$[b,s(x)]=[x_0,s(x_1),s(x_2),\dots,s(x_n)]$ have no common vertex except
$x_0=b$.
Indeed, if $x_j=s(x_i)$ with $i\neq j$, one could construct a path
of length  strictly smaller than $n$ joining $b$ either to $x$ or to 
$s(x)$. It follows that the geodesic path  $[x,s(x)]$ is  
the concatenation of the geodesic paths  $[x,b]=[x_{n},x_{n-1},\dots,x_{0}]$
and $[b,s(x)]=[x_{0},s(x_{1}),s(x_{2}),\dots,s(x_{n})]$.
Therefore, the distance $d(x,s(x))$ is even and $b$ is the midpoint of
$[x,s(x)]$.
\ep
\th{t5.1}{{\rm~\cite s} \
If $G$ is a finitely generated group of isometries of $X$ which is
a \gaf\;then $G$ is a \gag. Moreover, $\fix G$ is a nonempty subtree of $X$.
}
\pr
We proceed by induction on the number of generators of $ G $.
Suppose $ G $ is generated by an isometry $ s $ that fixes a vertex, and
a subgroup $ G_0 $ having a global fixed point $ x_0 $.
If $s(x_0)\neq x_0$ then, according to Lemma~\reff{l7.1},
the midpoint $z$ of the geodesic path $[x_0,s(x_0)]$ is a fixed point
of $s$.
Similarly, for all $ t\in G_0 $, $ st $ has a fixed point and
$[x_0, s(x_0)] = [x_0, s(t(x_0))]$. Therefore, $ z $ is fixed
by all the isometries belonging to the set $ sG_0 $.
Since $ s $ and $ sG_0 $ generate $ G $, the  point $z$ is fixed by $G$.
\ep

{\noi\sl Remark.}
The assumption ``$G$ is finitely generated'' is necessary.
Indeed, let $ [x_0,x_1,\dots] $ be an 
infinite geodesic path in a tree $ X $. Let $ k $ be an integer and let
$$
N(k)=\{s\in\isom X\tq\forall\ell\ge k,\,s(x_\ell)= x_\ell\}.
$$
The sequence $ (N (k))_{k \in \N} $ is a nondecreasing sequence of 
subgroups of $ \isom X $
whose intersection is a subgroup of $ G $. Any element $s$ in $G$ admits
fixed points but in general there is no global fixed point.
Especially, when $X$ is a homogeneous tree of degree $\ge3$,
the  group $G$ has no global fixed point and so $\isom X$
is nonfixating.
In fact the group $G$ admits a kind of {\sl fixed point at infinity}:
the equivalence class of geodesic paths ending as $[x_0, x_1,\dots]$.
\medskip
	
We now give an explicit condition to obtain a \gaf\,operating on a tree.
\prop{p7.7}{
Let $f$ and $g$ be two isometries of a tree. Assume that $f$, $g$
and $h=fg$ have fixed points. Then any element of the group
generated by $f$ and $g$ has fixed points.
}
\pr
It is enough to show that $\fix f $ meets $\fix g $.
If $ \fix f\cap \fix g = \emptyset$, let $ [P,Q]$ be the geodesic
joining $ \fix f $ to $ \fix g $, with $P\in\fix f$ and $Q\in\fix g$.
According to Lemma~\ref{l7.1}, $P$ is the midpoint of the geodesic
$ [Q, f (Q)] $. Also $f(Q) =f(g(Q)) = h(Q) $ and, as
$\fix h $ is nonempty, Lemma~\ref {l7.1} also implies that $P$ is
fixed by $ h $. We deduce $(fg)(P) = P $ hence $ g(P)=f^{-1}(P) = P $,
contradicting $ P \notin \fix g$.
\ep
\smallskip

{\noi\bf Bounded orbits in a tree}
\smallskip

The  median inequality~\rf6 almost holds for
combinatorial trees but an edge has no midpoint!
This is why isometries whose fixed points should be midpoints of edges
play a special role.
\lm{lem:arbrefini}{
If $X$ is a finite tree then there exists either a vertex fixed by every
element of $\isom X$ or an edge stable by every element of $\isom X$.
}
\pr
We proceed by induction on the number of vertices.
If there are one or two vertices, the result holds.
If there are strictly more than two vertices, the set of the vertices of $ X $
having at least two neighbors is a nonempty subtree $X'$
stable by the elements of $\isom X $.
Since there is at least one vertex that has a single neighbor in a finite tree,
the induction assumption can be applied to $X'$.
\ep
\df{5.l}{
An isometry $g$ of $X$ is called an {\sl inversion} if there exists an edge
$\{a,b\}$ such that $g(a)=b$ and $g(b)=a$.
}
\prop{p7.8}{
Let $X$ be a tree whose all vertices have a finite degree.
Let $G$ be a subgroup of $\isom X$ without inversion.
If $G$ has a bounded orbit then $G$ has a global fixed point.
}
\pr
Suppose $G$ has a bounded orbit $\Delta $.
Consider the set $ T(\Delta)$ of vertices in $ X $ that are in a geodesic
joining two vertices of $\Delta $.
Since $\Delta $ is finite, $ T(\Delta) $ is finite too.
In addition $T(\Delta)$ is the smallest subtree containing $\Delta$,
the {\sl convex hull} of the orbit.
Since an isometry sends a geodesic segment to a geodesic segment,
the subtree $ T(\Delta) $ is stable by the elements of
$ \isom X $. Therefore, according to Lemma~\reff{lem:arbrefini} applied to
$T(\Delta)$, there is either a vertex or an edge invariant by all elements
of $G$. If an edge $\{a,b\}$ is stable by $G$ then necessarily $G$ fixes $a$
and $b$ because $G$ does not contain an inversion.
\ep
	
The following result is an easy consequence.
\cor{c7.9}{
An element of finite order in the group of  isometries
of a tree which is not an inversion has a fixed point.
}
{\noi\sl Application.} 
Consider the  Schwartz group $G$ defined by two generators $a$ and $b$
related by the relations $a^{A} = b^{B} = (ab)^{C} = 1 $, where $A$,
$B$, and $C$ are
integers greater than or equal to $2$.
Any action by isometries without inversion of $G$ on
a tree $X$ has a global fixed point. Indeed, according to
Proposition~\reff{p7.8}
about bounded orbits, each of the isometries of $X$ defined by the
actions of $a$, $b$ and $ab$ has fixed points.
Then thanks to Proposition~\ref {p7.7} we conclude that
the group generated by $a$ and $b$ has a global fixed point.

By a similar argument Serre\index{Serre, Jean-Pierre}
proves that the group $\SL(3,\Z)$
has the same property: Each action by isometries of $G$ on
a tree has a global fixed point.
\med

{\noi\sl Generalization.} \,
One can define a notion of $\Lambda$-tree where $\Lambda$ is a
totally ordered Abelian group. 
Then Theorem~\ref{t5.1} still holds~\cite{ms}.
\sub{qgraphe}{Questions about isometries of finite graphs}\index{graph}
In the subsection we call {\sl graph} a simple unoriented connected graph.
The set of the vertices is equipped with the distance $d$
defined by the minimum number of edges joining vertices.
Several simple and natural questions arise in this context:
\medskip\\
1.
Which are the finite graphs whose isometry group is fixating? 
\medskip\\
2.
Which are the finite groups admitting a generating system defining
a Cayley graph whose isometry group is fixating?
\medskip\\
3.
Find infinite families of finite graphs whose isometry group is fixating.
\medskip

The first two questions are ambitious and probably difficult.
On the other hand, the third question admits simple partial answers.

Consider the complete graph $ K_n $ on $ n $ vertices.
Its isometry group  is the symmetric group $ \sy_n $.
Therefore, according to Theorem~\ref{t6.1},
it is fixating if and only if $ n \le 4 $.

Let $ \mathcal C_n$ be the graph associated with the $ n $-dimensional
hypercube $\{0,1 \}^n$  of $ \R^n $.
The set of vertices of $ \mathcal C_n $ is $\{0,1 \}^n$ and the edges
are the pairs of vertices of which exactly one
coordinate differs.

\prop{t6.4}{
The group $\isom\mathcal C_n$  is fixating
for all integer $n\ge1$.
}
\pr
The proof is analogous, in simpler form, to that of Theorem~\reff{t6.2}
on $ \Z^n $.
It is easy to prove that an isometry of $ \mathcal {C}_n $ is the restriction to
 $ \mathcal {C}_n $ of a unique isometry of the Euclidean space
$\R^n$, cf. Exercise~\reff{e54}.
The center of the cube $\omega=(\tfrac12,\dots,\tfrac12)$ is fixed by
all the isometries of the cube.
Conjugating by the map $x\mapsto2x-2\w$, we reduce to the
case where the isometries are linear maps whose matrices
are matrices of signed permutations.
We find that the matrices of the isometries of a subgroup \gaf\
never have $ -1 $ on the diagonal, which implies
the existence of a global fixed point as in Steps 5 and 6
of the proof of Theorem~\ref{t6.2}.
\ep
\sub{graphe}{A result about infinite graphs}
Let $ X =(S,A) $ be a simple undirected connected graph whose edges
are colored\index{colored graph}.
The color of an edge is given by a map with values in a color set
$\mathcal C $, defined on the set $A$ of edges.
Recall that a cycle is called {\sl simple}\index{simple cycle} if no edge occurs more
than once, and  {\sl elementary}\index{elementary cycle} if no vertex occurs more than once
except the beginning and the end.
It is easy to see that an elementary cycle of length at least $ 3 $
is simple and that a cycle, simple or not, always contains at least one
elementary cycle.
We make the following assumptions:
\begin{itemize}
\item
The edges of an elementary cycle of the graph all have the same color,
\item
For each $ c \in \mathcal C $, the connected components of the partial graph
$ X_c $ obtained by keeping only the edges of color $c$
are complete graphs.
\end{itemize}
For every $c\in\mathcal C$, we will call any connected component of the
partial graph $X_c$ {\sl a cell of color $c$}.
Each edge of $X$ belongs to a single cell.

Notice that the Cayley graph of the free product of two groups
$G_1$ and $G_2$ satisfies the above assumptions if we choose
$(G_1\cup G_2)\setminus\{e\} $ as set of generators and if we
color in blue the edges $\{w,wg\}$, $g \in G_1$, and in red the others.
\med

Here is a generalization of Theorem~\reff {t5.1} on finitely generated groups 
of isometries of trees.
\th{t6.13}{
Suppose that each cell of the graph $ X $ has at most four vertices.
Let $G$ be a finitely generated subgroup of $ \isom X $.
If $G$ is a \gaf\;then $G$ is a \gag.
}
Observe that the hypothesis about the cardinality of the cells is necessary:
If $ X $ is a complete graph with at least five vertices then $\isom X$
contains an eccentric subgroup by Theorem~\ref{t6.1}.
It is the same for the Cayley graph of the free product of two
groups of which one at least has more than  five elements.
The proof of Theorem~\reff{t6.13} needs some
preliminary lemmas.
\lm{l6.9}{
Given two vertices $x$ and $y$ of the graph $ X $, there is a unique
path joining them such that two consecutive edges on this path never have
the same color. Moreover this path is the only geodesic from $x$ to $y$.
Therefore a path is a geodesic if and only if two
consecutive edges never have the same color.
}
\pr
Let $ x$ and $y$ be in $S$. Since $X$ is connected, there is a path
joining $ x $ to $ y $, and therefore at least one geodesic.
Since the cells of the graph are complete graphs, this geodesic
cannot have two
consecutive edges of the same color. This proves the existence.

For the uniqueness, let us first notice that all the vertices of a path
verifying the color change property of the edges, are distinct.
Indeed, if two vertices of the path coincided then an elementary cycle
could be extracted.
By assumption on $X$, all the edges of this cycle have the same color,
contradicting the property of color change.

If two distinct paths $ [a_0 = x, \dots, a_m = y] $ and $ [b_0 = x, \dots, b_n = y] $
join $ x $ to $ y $ and satisfy the  property
of color change then there is an integer $ i\ge 0 $ such that $ a_i = b_i $ and
$ a_{i + 1} \ne b_{i + 1} $. Consider the first vertex $ a_j $ in
the path $ [a_{i + 1}, \dots, a_m] $ which also belongs to the path
$ [b_{i + 1}, \dots, b_n] $. So we have $ a_j = b_k $ for some $ k> i $.
The integer $k$ is chosen minimal.
By choice of $ j $ and $ k $, the vertices of the path
$ [a_i, a_{i + 1}, \dots, a_j = b_k, b_{k-1}, \dots, b_i = a_i] $ are all
distinct, except the two ends.
This cycle has at least three edges (otherwise $ a_{i + 1} = b_{i + 1} $),
it is elementary and has at least two colors, a contradiction.
\ep
\lm{l6.11}{
Let $ s $ be an isometry of $ X $ having at least one fixed point and
let $ x $ be a vertex of
$ X $ that is not a  fixed point of $s$. Let $F$ be
the set of points $ y \in \fix s $ such that $ d(x,y) = d(x,\fix s) $.
\be
\item
Then $ F $ is included in a cell whose color is the one
of the last edge of the geodesic going from $ x $ to any
element of $ F $.
\item
Moreover:
\bee
\item If $d(x,s(x))$ is even then the midpoint $z$ of the geodesic
joining $x$ to $s(x)$ belongs to $F$.
\item If $d(x,s(x))$ is odd then $F$ is included in the cell
$Y$ containing the middle edge $[a,b]$ of  the geodesic going from $x$
to $s(x)$. In addition, $s(Y)=Y$ and $s(a)=b$.
\ee
\ee
}
\pr 
{\bf a}. 
Let $ u $ and $ v $ be two points in $ F $. Consider the geodesics
$ [u_0 = x, \dots, u_n = u] $ and $ [v_0 = x, \dots, v_n = v] $
(with $ n = d(x, \fix s) $).
Since $ s $ is an isometry that fixes $ u $ and $ v $,
it fixes each vertex of the geodesic $ [w_0 = u, \dots, w_m = v] $
joining $ u $ to $ v $. Let $ i $ be the smallest integer such that $ u_i = v_i $ and
$ u_{i + 1} \ne v_{i + 1} $.
We verify, as in the proof of Lemma~\reff {l6.9},
that the three branches of geodesics $ [u_i, \dots, u_n = u] $,
$ [u, \dots, v] $ and $ [v_n, \dots, v_i = u_i] $ form an elementary cycle.
Therefore, all the edges of this cycle have the same color
$ c $ and the vertices $ u $ and $ v $ are joined by an edge of color $ c $.
We conclude by noticing that $ c $ is the color of the edge $ [u_{n-1}, u_n] $.
\med\\
{\bf b}.
Let $ u $ be a point in $ F $ and let $ [u_0 = x, \dots, u_n = u] $ be the geodesic
joining $ x $ to $ u $.
According to item~{\bf a}, $F$ is included in a cell  $Y$ 
whose color $ c $ is that of the edge $ [u_{n-1}, u_n] $.
Consider the image $ [s (x), \dots, s(u_n) = u] $ of this geodesic
by $ s $. Let $ c'$ the color of the edge $ [s(u_{n-1}), u] $.
\med\\{\bf(i)}
If $c\ne c'$ then
$[u_0,\dots,u_n=s(u_n),s(u_{n-1}),\dots,s(u_0)=s(x)]$ is a geodesic
since the colors of two  consecutive edges of this path are never the
same. Hence it is the geodesic joining $x$ to $s(x)$.
The length of this geodesic is $2n$ and its midpoint is $u\in F$. 
\med\\{\bf(ii)}
If $ c=c'$ then $ s(u_{n-1}) $ is a vertex of $ Y $ which is complete,
so $ [u_{n-1}, s(u_{n-1})] $ is also an edge
$ Y $. Two consecutive edges of the path
$ [u_0, \dots, u_{n-1}, s(u_{n-1}), $ $ \dots, s(u_0)] $ never have the same color.
Hence this path is the geodesic joining $ x $ to $ s(x) $.
By construction, the length of this geodesic is odd and its middle edge
 $ [u_{n-1}, s(u_{n-1})] $ is of color $ c $.
It remains to prove that $ s(Y) = Y $. Let $ y $ be a vertex of $ Y $.
The image under $ s $ of the triangle $ u_{n-1} u y $ is a triangle that contains
the edge $ [s(u_{n-1}), u] $ which is an edge of $ Y $.
Therefore the image of this triangle is a triangle of $ Y $
hence $ s(y) $ is a vertex of $ Y $.
\ep
\med\\
{\noi\sl Proof of Theorem~\reff{t6.13}.}
We proceed by induction on the number of generators of $G$.
Suppose that $G$ is generated by $s$ and $G_0$, where $s$ is
an isometry which fixes at least one vertex and $G_0$ is a
subgroup with a global fixed point $x_0$. If $s(x_0)=x_0$ we are done.
Otherwise consider the geodesic going from $ x_0 $ to
$ s(x_0) $. Observe that it is also the geodesic going from
$ x_0 $ to $ st(x_0) = s(x_0) $ for any $ t \in G_0 $.

If this geodesic has an even length then, according to Lemma~\reff{l6.11},
its midpoint $ z $ is a fixed point of $ s $ but
also of all $ st $, $ t \in G_0 $. Since $ s $ is injective,
$ s (t (z)) = z = s (z) $ implies $ t (z) = z $, so $ z $ is a fixed point of $ t $.
Thus $ z $ is a global fixed point of $ G $.

If this geodesic has an odd length then, again according to 
Lemma~\reff{l6.11}, the cell $ Y $ containing the middle edge of
this geodesic is stable by all $ st $, $ t \in G_0 $.
It is therefore stable by $ G $. In addition, the cell $Y$ contains fixed points of
$ st $. Therefore the  group $H$ of  restrictions to $Y$ of
elements of $ G $ is a \gaf. Now $ \isom Y $ is the group of
permutations of the vertices of $ Y $  and by hypothesis the  cell
$ Y $ has at most $ 4 $ vertices so, according to Theorem~\ref{t6.1},
$ \isom Y $ is
fixating. Therefore $ H $ is a \gag\;hence $G$ is a \gag\;too.
\ep
\sec{7.}{Exercises} 
\exo{e0}{
Show that the group $\isom \R^2$ is fixating.
}
\exo{e1}{
A group is called {\sl superfixating}
if, for any set $X$ and any morphism $\rho:G\to\bij X$, the pair
$(X,\rho(G))$ is fixating.
By considering the action of a group on all of its
nontrivial parts, show that a group is superfixating
if and only if it is cyclic (finite or not).
}
\exo{e2}{
Show that the additive group $\Q$ is {\sl finitely superfixating}
in the following sense: If $X$ is a finite set and $\rho:\Q\to\bij X$
a morphism then $(X,\rho (\Q))$ is fixating.
}
\exo{e3}{
Prove that, in a metric space verifying the median inequality~\rf6,
the point $m$ is unique and satisfies
$d(x,m)=d(y,m)=\tfrac12\,d(x,y)$.
}
\exo{e4}{
Let $(E,\|\;\|)$ be a normed vector space verifying the median
inequality~\rf6.
\be
\item
Prove that $E$ satisfies the so-called {\sl parallelogram identity}
\eq4{
\forall x,y\in E\qquad\|x+y\|^2+\|x-y\|^2=2\big(\|x\|^2+\|y\|^2\big).
}
\item
Deduce that $E$ is a pre-Hilbert space (a result of M.~Fr\'echet,
P.~Jordan and J.~von
Neumann~\cite{f,jn})\index{Fr\'echet, Maurice}\index{Jordan, Pascual}\index{Neuman@von Neumann, John}.
\ee
}
\exo{e6}{
According to the comparison theorem of Rauch~\cite{ec}\index{Rauch, Harry},\index{Rauch comparison theorem}\index{Theorem!Rauch comparison}
the classical cosine law becomes an inequality in hyperbolic trigonometry: 
In a hyperbolic triangle of side lengths $a,b,c$
and angle $\g$ opposite to the side of length $c$, we have
\eq{cos}{
a^2+b^2-2ab\cos\g\le c^2.
}
From this inequality, show that $\H_n$ satisfies the median inequality~\rf6.
}
\exo{ex:homographie1}{
On the set $\hat{\R}=\R\cup\{\infty\}$, consider the {\sl M\"obius group}
\index{M\"obius group}
$$
M(\hat {\R})=\big\{\f:\hat\R\to\hat\R,\,x\mapsto\tfrac{ax+b}{cx+d}
\tq a,b,c,d\in \R,\, ad-bc=\pm 1\big\}.
$$
Show that this group is nonfixating.
Hint: Consider the matrices $A$ and $B$ of Example~\reff{e2.4}
and use Lemma~\ref{p}.
}
\exo{e12}{
Using only the fact that $\isom\H_3$ is fixating,
prove that $\isom\H_2$ is fixating.
}
\exo{e4b}{
Show that a map, a priori surjective or not, from $\R^n$ to $\R^n$
which preserves the Euclidean distance is an affine bijection of $\R^n$.
}
\exo{e5}{
Show that a function between two normed vector spaces
which is continuous and preserves the midpoints is affine.
}
\exo{exo:mazurulam}{
{\sl Mazur-Ulam theorem.\index{Mazur, Stanislaw}\index{Ulam, Stanislaw}} 
This exercise is inspired by~\cite v.
\smallskip

\noi
Let $E$ be a real normed vector space.
\be
\item
Let $a,b\in E$ and $m$ be the midpoint of $[a,b]$, let $W_{a,b}$ be the
set of isometries of $E$ fixing $a$ and $b$,
and let $\l=\sup\{\Vert g(m)-m\Vert\tq g\in W_{a,b}\}$.
\bee
\item
Prove that $\l\le\Vert a-b\Vert$.
\item
Let $s_m$ be the symmetry of center $m$, i.e. such that
$s_m(x)=2m-x$ for all $x\in E$.
For $g\in W_{a,b}$, we set $g^*=s_m\,g^{-1}s_m\,g$. Prove that
$\Vert g^*(m)-m\Vert=2\Vert g(m)-m\Vert$. 
\item
Show that any isometry that fixes $a$ and $b$ fixes $m$.
\ee
\smallskip

\item
Let $f$ be an isometry of $E$.
\bee
\item
Let $a,b\in E$. Denote by $m$ the midpoint of $[a,b]$ and by 
$m'$ the one of $[f(a),f(b)]$.
Prove that $h=s_m\,f^{-1}s_{m'}\,f\in W_{a,b}$ and deduce that $f(m)=m'$.
\item
Prove that $f$ is affine. Hint: Use Exercise~\ref{e4b}.
\ee
\ee
}
\exo{e9}{
{\sl Kakutani Theorem in finite dimension}\index{Kakutani, Shizuo}.
This exercise is inspired by R.~Antetomaso~\cite{a}.
Let $E$ be a normed vector space of finite dimension, 
$G$ a compact subgroup of $GL(E)$,
and $K$ a nonempty compact convex subset of $E$.
We assume that $g(K)\subseteq K$ for all $g\in G$ and we aim to show
that $(K, G)$ is a \gag.
%
\be
\item 
Let $\Vert\;\;\Vert_2$ be the Euclidean norm on $E$.
For any $x\in E$, we set
$\Vert x\Vert=\sup\{\Vert g(x)\Vert_2\tq g\in G\}$.
Show that this defines a strictly convex norm on $E$
for which every element of $G$ is an isometry.
\item
Let $f$ be an endomorphism of $E$ such that $f(K)\subseteq K$.
Let $x_1\in K$ and let $(x_n)_{n> 0}$ be the sequence defined by
$x_{n+1}=f(x_n)$. 
Considering the sequence $(\sigma_n)$ of Ces\`aro means of
the sequence $(x_n)$, show that $f$ has a fixed point in $K$.
\item
For $g\in G$, denote $V_g=\{x\in K\tq g(x)\neq x\}$.
By contradiction, assume that for every $x\in K$
there exists $g\in G$ satisfying $g(x)\neq x$.
\bee
\item
Show that there exist $g_1,\ldots,g_N\in G$ such that
$K\subset V_{g_1}\cup\cdots\cup V_{g_N}$.
\item
Show that there exists $a\in K$ such that $(g_1+\cdots+g_N)(a)=Na$.
\item
Show that $g_k(a)=a$ for all $k\in\{1,\ldots,N\}$. Conclude.
\ee
\ee
}
\exo{e5.3}{
{\sl The isometries of the sphere extend into isometries of
the Euclidean space.}\index{extension of an isometry}
\\
We equip the sphere $\S_n$ with the spherical distance
$d(x,y)=\arccos\scal xy$.
We want to prove that every isometry of  $\S_n$ is the restriction of
a unique isometry of $\R^{n+1} $ endowed with the Euclidean distance.
Let $ f $ be an isometry of $\S_n$.
\be
\item
Show that $f$ preserves the scalar product.
\item
We define $\~f$ on $\R^{n+1}$ by
$\~f(x)=\Vert x\Vert f\big(\tfrac{x}{\Vert x \Vert}\big)$ if $x\neq\v0$
and $\~f(\v0)=\v0$. Show that $\~f$ preserves the scalar product.
\item
Conclude.
\ee
}
\exo{e54}{ 
{\sl Extension of an isometry.\index{extension of an isometry}}
\\
Let $A$ be a subset of $\R^d$ and set $E=\aff A$, the affine subspace
generated by $A$. 
We endow $A$ with the Euclidean distance induced by that of $\R^{d}$.
We want to show that any isometry of $A$ is the restriction of a
unique affine isometry of $E$.

We know that there are $a_0,\dots,a_n\in A$
such that $E=\aff(a_0,\dots,a_n)$, where $n=\dim E$.
Thus every point $x\in E$ is written in a unique way
$x=\sum_{i=0}^n\l_i(x)a_i$ with $\sum_{i=0}^n\l_i(x)=1$.
The numbers $\l_i(x)$ are the {\sl barycentric coordinates} of $x$.
Let $f$ be an isometry of $A$.
\be
\item
Show that for all $i,j\in\{0,\dots,n\}$ we have
$\scal{f(a_i)-f(a_0)}{f(a_j)-f(a_0)}=\scal{a_i-a_0}{a_j-a_0}$.
\item
We define $\~f$ on $E$ by $\~f(x)=\sum_{i=0}^n\l_i(x)f(a_i)$.
Prove that $\~f$ is an isometry.
\item
Show that $\~f$ extends $f$. Conclude.
\ee
}
\exo{e15}{
{\sl The alternating group $\A_7$ is nonfixating.\index{alternating group}}
\\
In Step 3 of the proof of Theorem~\reff{t6.1},
we built an eccentric subgroup of $\sy_5$
using permutations which act separately on 
$\{1,2,3\}$ and on $\{4,5\}$, thanks to the following key point.
The group $\sy_3$ has a normal subgroup
(the alternating group $\A_3$) which has the following two properties:
\begin{itemize}
\item[$\triangleright$]
The quotient $\sy_3/\A_3$ is isomorphic to $\Z/2\Z$,
\item[$\triangleright$]
Any element of $\sy_3\setminus\A_3$ has at least one fixed point.
\end{itemize}
We then obtained an eccentric subgroup $G=\langle f,g\rangle$ of
$\sy_5$ by taking for $f$ an element of $\A_3$ on $\{1,2,3\}$
and the identity on $\{4,5\}$, and for $g$ an element of
$\sy_3\setminus\A_3$ on $\{1,2,3\} $ and an element of $\sy_2$
without fixed point on $\{4,5\}$.
Use a similar construction to build an eccentric subgroup of
$\A_7$.
}
\exo{e16} {
{\sl The action of the group $\GL(3,\F_2)$ on $X=\F_2^3\setminus\{\v0\}$
is nonfixating.}
\\
Let $(e_1,e_2,e_3)$ be the canonical basis of $\F_2^3 $ and,
for $A\subseteq\{1,2,3\}$, let us denote $e_A=\sum_{i\in A}e_i$.
In this manner we have $X=\{e_A\tq A\neq\emptyset\}$.
Let $f$ and $g$ be the elements of $\GL(3,\F_2)$ defined by
\[
f(e_1)=e_2,\,f(e_2)=e_3,\,f(e_3)=e_1\quad\text{and}\quad
g(e_1)=e_{123},\,g(e_2)=e_2,\,g(e_3)=e_3.
\]
Let $G=\langle f,g\rangle$. Finally denote  $Y=\{e_1,e_2,e_3,e_{123}\}$
and $Z=\{e_{23},e_{13},e_{12}\}$.
\be
\item
Show that $f$ and $g$ induce even permutations on $X$.
\item 
Show that the map which sends $h\in G$ to its restriction on $Y$
induces an injective morphism from $G$ to the group of permutations
$\sy_Y$ of $Y$.
\item
Compute $(fg)^2$, $(gf)^2$, and $(fgf)^2$.
Deduce that, if the restriction of $h\in G$ to $Y$ is a double transposition,
then its restriction to $Z$ is the identity.
\item
Deduce that the action of $\GL(3,\F_2) $ is nonfixating.
\item
Observe that this produces another eccentric subgroup of $\A_7$.
\ee
}
\exo{e17}{
{\sl The action of the group $\GL(d,\k)$ on $X=\k^d\setminus\{\v0\}$
is nonfixating, except in the obvious cases
$\GL(1,\k)$ and $\GL(2,\F_2)$.}
\\
Let $\k$ be a field, finite or not, let $d\ge1$ be an integer,
and let $X=\k^d\setminus\{\v0\}$. In this exercise we will say
that the group $\GL(d,\k)$ is fixating if its action on $X$ is.
\be
\item
Prove that, if $\k$ has at least three elements,
then $\GL(2,\k)$ is nonfixating.
\item
Prove that, if $\GL(d,\k)$ is nonfixating, then
$\GL(d+1,\k)$ is nonfixating either.
\item
Deduce that $\GL(d,\k)$ is fixating if and only if:

\bul either $d=1$, 

\bul or $d=2$ and $\k=\F_2=\{0,1\}$. 
\ee
}

%
%
\sec{8.}{Appendices} 
%
\sub{8.1}{
A short introduction to hyperbolic geometry
}\index{hyperbolic space}
Let $n\ge2$ be an integer. A model 
of the $n$-dimensional hyperbolic space is the upper half-space
\[
\H_n=\{(x_1,\dots,x_n)\in \R^n\tq x_n>0\}
\]
endowed with the Poincar\'e metric 
\[
ds^2=\frac1{x_n^2}\,\big(dx_1^2+\cdots+dx_n^2\big).
\]
A  calculation shows that the geodesic distance associated with the 
Riemannian metric is given by 
\eq c{
d(x,y)=\operatorname{argcosh}\left(1+\frac{\|x-y\|^2}{2x_ny_n}\right)
}
where $x=(x_1,\dots,x_n)$ and $y=(y_1,\dots,y_n)$. 
\th{geodesiques}{{\rm (\cite{cfkp}, Theorem 9.3)}
The geodesics in $\H_n$ are the half-lines (affine lines) and the 
half-circles (Euclidean circles), the endpoints of which are in the
horizontal hyperplane  
\[
\partial\H_{n}=\{(x_1,\dots,x_n)\in \R^n\tq x_n=0\},
\] 
and orthogonal to this hyperplane at their endpoints. 
}
\ssub{sec:groupehyp}{The isometry group of $\H_n$} 
The similarity transformations (similarities for short),
the inversions with respect to spheres,
and the reflections through affine hyperplanes, act on
$\widehat \R^n=\R^n\cup\{\infty\}$ and
form a subgroup of the group of homeomorphisms of $\widehat \R^n$.
This subgroup is call the {\sl M\"obius group\index{M\"obius group} of}
$\widehat \R^n$ and is denoted by $M(\widehat \R^n)$. 
The inversions and the reflections are enough to generate the
Möbius group.
The restrictions to  the upper half-space of some elements of the 
Möbius group give rise to the isometry group of $\H_n$:
\th{th:Mobius}{{\rm (\cite{r}, Theorem 4.6.2)}
The isometry group $\isom\H_n$ is the group of the restrictions to
$\H_n$ of the M\"obius transformations $\f$ such that $\f(\H_n)=\H_n$.
It is generated by the reflections through spheres centered
at points in $\partial\H_{n}$ and through vertical hyperplanes.
}
All results we need follow from this theorem,
from Theorem~\ref{geodesiques}, and from formula~\rf c.

A first consequence of Theorem~\reff{th:Mobius} is that a similarity
$\f:\R^n\rightarrow \R^n$ such that $\f(\H_n)=\H_n$ induces an isometry
of $\H_n$.
Conversely, since the only elements in $M(\widehat \R^n)$ that fix $\infty$
are the similarities (\cite{r}, Theorem 4.3.2), we have the following result.
\cor{cor:infinifixe}
{ If an isometry $f$ of $\H_n$ is the restriction of a map 
$\widehat f$ in $M(\widehat \R^n)$ which fixes $\infty$ then the restriction 
$\widehat f$ to $\R^n$ is a similarity in $\R^n$.
}
A second consequence of Theorem~\reff{th:Mobius}
is that  isometries of $\H_n$ are smooth. On the one hand it follows that
the Poincar\'e metric is invariant, on the other hand
it follows that the sign of the Jacobian of an isometry
is constant. This latter fact leads to the usual decomposition of the group
$\isom\H_n$:
It is the union of the subgroup $\isomp\H_n$ of isometries with positive
Jacobian and its complementary $\isomm\H_n$.
\ssub{sec:sousespacehyp}{Subspaces of $\H_n$}
A {\sl hyperbolic subspace\index{hyperbolic subspace} of} $\H_n$ is a subset
$X$ of $\H_n$ of one of the following forms:
\begin{itemize}
\item[$\triangleright$]
the empty set,  $\dim X=-1$,
\item[$\triangleright$]
a single point, $\dim X=0$,
\item[$\triangleright$]
the intersection of $\H_n$ with a vertical affine subspace $A$,
$\dim X=\dim A$,
\item[$\triangleright$]
the intersection of $\H_n$ with a vertical affine subspace $A$
and a sphere $S$ centered at a point of $\partial\H_{n}$,
$\dim X=\dim(A\cap S)$ ($A=\R^n$ is permitted).  
\end{itemize}
By abuse of language, subspaces of dimension 1 will be called
{\sl lines} and those of dimension 2 {\sl planes}.
Using that the image of a sphere by an inversion 
is an affine subspace provided that the pole of the inversion belongs 
to the sphere,
we see that a $p$-dimensional hyperbolic subspace is isometric to
\[
\{(x_1,\dots,x_n)\in \H_n\tq x_1=x_2=\dots=x_{n-p+1}=0\}
\] 
and therefore to $\H_p$.
By definition, the {\sl geodesics} are the lines of the hyperbolic space.
Moreover, a  hyperbolic subspace is {\sl totally geodesic}\index{totally geodesic},
which means that if it contains two points $x$ and $y$ then it contains
the geodesic joining $x$ to $y$.
\lm{lem:intersection}{
Every intersection of  hyperbolic subspaces in $\H_n$
is a  hyperbolic  subspace.
}
\pr
Every intersection of spheres and/or affine subspaces is a sphere $S$
(or an affine subspace if there is no sphere)
of the affine subspace $A$ generated by the intersection.
Observe that this subspace $A$ can be empty or reduced to a single point.
We need to check that, if the affine subspaces are all vertical,
if the spheres are centered in $\partial\H_{n}$, and if
$\dim A\ge1$, then $A$ is vertical and the center of $S$ belongs
to $\partial\H_{n}$.  If the intersection is defined by at least 
two hyperbolic subspaces  then it is included
in a vertical  hyperplane. Since such a hyperplane is isometric to
$\H_{n-1}$, the result follows by induction on $n$.
\ep	
\lm{lem:sousespaceengendre}{
Let $H$ be a subspace of $\H_n$ and $x\in\H_n\setminus H$.
The dimension of the smallest subspace containing $H$ and $x$ is $\dim H+1$.
}
\pr
Let $f$ be an isometry that sends $H$ onto the intersection of $\H_n$
with a vertical affine subspace  $A$ of dimension $p=\dim H$.
The smallest hyperbolic subspace containing $A\cap\H_n$ and $f(x)$ is
$A'\cap\H_n$ where $A'$ is the affine subspace generated by $A$ and $f(x)$,
which is of dimension $p+1$.
\ep
\lm{lem:med}{
Let $a$ and $b$ be two distinct points in $\H_n$.
Then the set $\md (a,b)$ of  points equidistant from $a$ and $b$ is a
hyperplane of $\H_n$.
}
\pr
By \rf c, a point $x$ is in  $\md(a,b)$ if and only if 
\[
d(x,a)=\operatorname{argcosh}\left(1+\frac{\|x-a\|^2}{2x_na_n}\right)=
d(x,b)=\operatorname{argcosh}\left(1+\frac{\|x-b\|^2}{2x_nb_n}\right),
\]
that is
\[
\frac{\|x-a\|^2}{2x_na_n}=\frac{\|x-b\|^2}{2x_nb_n}.
\]
In the case where $a_n=b_n$, we get a vertical hyperplane.
In the case where  $a_n\neq b_n$, we get a half-sphere the center of which
is in $\partial\H_{n}$. Indeed the coefficient of
$x_n$ vanishes in the cartesian equation of the sphere, given by
$b_n\|x-a\|^2-a_n\|x-b\|^2=0$.
\ep
\lm{lem:sousespacefixe}{
If $f$ is an isometry of $\H_n$ then $\fix f$ is a
hyperbolic subspace.
}
\pr
By  Lemma \ref{lem:intersection}, there exists a smallest hyperbolic subspace
 $H$ containing $\fix f$. Let $x\in H$ be arbitrary. 
By contradiction, if $f(x)\neq x$ then $\md(x,f(x))$ is a hyperplane which
contains $\fix f$ hence $\fix f\subset F=H\cap\md(x,f(x))$
but $F$ does not contain $x$. This contradicts the minimality of $H$.
Therefore $x=f(x)$; this proves $\fix f=H$.
\ep
\lm{lem:projdroite}{
Let $D$ be a line in $\H_n$ and $x\in\H_n$.
There exists a unique point $y\in D$ such that 
\[
d(x,y)=d(x,D):=\inf\{d(x,z)\tq z\in D\}.
\]
This point is denoted $\pi_Dx$.
Moreover, if $x\notin D$ then the geodesic through $x$ and $\pi_Dx$
is orthogonal to $D$.
}
Recall that the Poincar\'e metric is conformal to the Euclidean metric
hence the orthogonality relations are equivalent for these two metrics.
\medskip

\pr
By compactness, the distance from $x$ to $D$ is realized:
There exists $y\in D$ such that $d(x,D)=d(x,y)$.
The point $x$ and the line $D$ are included in some hyperbolic plane $P$.
On the one hand the geodesic from $x$ to $y$ is included in $P$.
On the other hand this plane is isometric to the hyperbolic plane $\H_2$.
Therefore we can suppose that $x$ and $D$ are in $\H_2$.
Thanks to another isometry, we can suppose that $D=\{ti\tq t>0\}$.
Denoting $x=a+bi$, we have 
\[
d(x,ti)=\operatorname{argch}\left(1+\frac{a^2+(b-t)^2}{2bt}\right).
\]  
A calculation shows that the function $t\rightarrow\frac{a^2+(b-t)^2}{2bt}$
attains its minimum at  $t=\sqrt{a^2+b^2}$ hence $\pi_Dx=i\sqrt{a^2+b^2}$.
The geodesic joining $x$ to its projection is the half-circle of center
$\v0$ and radius $\sqrt{a^2+b^2}$, which is indeed orthogonal to $D$.
\ep
\cor{cor:proj}{
If $F$ is a hyperplane of $\H_n$ and $x$ is in $\H_n$
then there exists a unique point $\pi_Fx\in F$ such that
\[
d(x,\pi_Fx)=\inf\{d(x,z)\tq z\in F\}=d(x,F).
\]
Moreover, if $x\notin F$ then the geodesic going through $x$ and
$\pi_Fx$ is orthogonal to $F$. 
}
\pr
By compactness, the minimal distance is attained in at least one point.
If two distinct points $y,z\in F$ give this minimal distance
then these two points also give the minimum distance from $x$ to $D$,
where $D$ is the geodesic going through $y$ and $z$, since $D\subset F$.
This contradicts Lemma~\reff{lem:projdroite}.
This lemma also implies that the geodesic from $x$ to $\pi_Fx$ is
orthogonal to all the lines going through $\pi_Fx$ and included in $F$
hence is orthogonal to $F$.
\ep

{\noi\sl Remark.}
As already said, the uniqueness of the orthogonal projection on hyperbolic
subspaces also comes from Proposition~\ref{p4.5} and Exercise~\ref{e6}.
\ssub{sec:h3}{Isometries of $\H_3$}
The next lemma asserts that an angle and an axis are associated with
a positive isometry of $\H_3$ which has fixed points.
\lm{lem:rotationH3}{
Let $f\in \isomp \H_3$. Assume $\fix f\ne\emptyset$ and $f\ne\id$.
\be
\item
Then $\fix f$ is a hyperbolic line.
\item
There exists an angle $\theta$ such that for all $x\in \fix f$
the differential $df(x)$ is a rotation  with angle $\pm\theta$
and with axis the affine line tangent to $\fix f$.
\ee
}
\pr
{\noi\bf a.}
Let $x\in \fix f$. The differential $df(x)$ is a positive isometry 
of $\R^3$ hence a rotation with an eigenvector $\v u$.
Consider the geodesic $D$ through $x$ of direction $\v u$. 
It is invariant by $f$, and $f$ has a fixed point in $D$,
hence $D\subset\fix f$.
If $D$ were strictly included in $\fix f$  then $\fix f$ would contain a
hyperbolic plane. The differential of $f$ at a point of this plane would be a
rotation with two independent eigenvectors; hence it would be the identity.
It follows that $f$ would be the identity, a contradiction, hence $\fix f=D$.
\med

{\noi\bf b.}
Conjugating by an isometry, we can restrict ourselves to the case where $D$ is
a vertical line. By Corollary \ref{cor:infinifixe}, $f$ is the restriction
of a similarity. 
Since $f$ has a half-line of fixed points, $f$ is the restriction of a
rotation $r$ of axis $D$. The differential of $f$ along this axis is then
always equal to $\v r$ or to its inverse.
\ep

With the same method one easily proves the following statement.
\lm{lem:symetrie}{
Let $f\in\isomm\H_3$. If $\fix f\ne\emptyset$ then $\fix f$
is a hyperbolic plane and $f$ is a reflection through this plane.
}
Our last statement concerning isometries of $\H_3$ is the following.
\lm{lem:rotationsymetrie}{
Let $f\in\isomp \H_3$. Assume that $\fix f\neq\emptyset$.
Denote by $D$  the axis of $f$ and by $\theta$ its angle.
Let $P_1$ and $P_2$ be two hyperbolic planes of intersection
$D$ and let $\sigma_{P_1}$ and $\sigma_{P_2}$
be the  reflections through the planes $P_1$ and $P_2$.
Denote $\partial P_1$ and $\partial P_2$ the circles or lines determined
by these planes in $\partial \H_3$.
If  $\angle(\partial P_1,\partial P_2)=\pm\tfrac\q2\mod\pi$ then
$f=\sigma_{P_2}\sigma_{P_1}\text{ or } \sigma_{P_1}\sigma_{P_2}$.
}
\pr
As in the previous proof, thanks to a conjugacy by an isometry,
we can assume that $D$ is a vertical half-line.
It follows that $P_1$ and $P_2$ are vertical planes 
and that $f$, $\sigma_{P_1}$, and $\sigma_{P_2}$ are similarities
hence isometries of $\R^3$.
Since the relation about angles remains true up to the sign, we have
$f=\sigma_{P_2}\sigma_{P_1}\text{ or } \sigma_{P_1}\sigma_{P_2}$. 
\ep
\ssub{sec:h2}{Isometries of $\H_2$}
The hyperbolic plane $\H_2$ is identified with the upper half-plane
$\{z\in\C\tq\im z>0\}$. 
With this identification, by Theorem \ref{th:Mobius},
the positive isometries of $\H_2$ are the conformal transformations
of the upper half-plane, i.e., the maps $h:\H_2\rightarrow\H_2$ of the form 
\[
h(z)=h_{a,b,c,d}(z)=\frac{az+b}{cz+d}\;\tq\hspace{10pt} a,b,c,d\in\R,\,ad-bc=1.
\]
The negative isometries are obtained by composing the positive isometries
with the reflection $z\mapsto-\overline z$.
\lm{lem:transitiviteh2}{
If $x,y,x',y'\in\H_2$ are such that $d(x,y)=d(x',y')$
then there exists $f\in\isomp \H_2$ such that $f(x)=x'$ and $f(y)=y'$.
}
\pr
It is enough to prove the statement in the case where $x'=i$ and $y'=it$
with $t>1$.
The homography $h\in\isomp H_2$ defined by $h(z)=\frac{1}{\im x}(z-\re x)$
sends $x$ to $i$.
Let $\v u$ be the unit vector tangent at $x$ to the geodesic going
from $x$ to $y$.
By composing with a homography of the form
$r(z)=\frac{\cos\q\,z-\sin\q}{\sin\q\,z+\cos\q}$ which fixes $i$, we can 
suppose that the differential of $rh$ at $x$ sends $\v u$ to  $i$.
The isometry $rh$ maps the geodesic segment joining $x$ to $y$ onto a
geodesic segment
of the same length starting at $i$ with initial speed $i$ hence
included in the imaginary axis. Therefore $rh(y)=it$ where $t$ is the
unique real number greater than $1$ such that $d(i,it)=d(x,y)=d(x',y')$
hence $rh(y)=y'$. 
\ep
\np
{\noi\sl Proof of Lemma~\reff{l6.3}.}
\med

{\noi\bf a.} \
The search of the fixed points of $h=h_{a,b,c,d}$ leads to the equation
\eq8{
cz^2+(d-a)z-b=0,
}
the discriminant of which is
$\Delta=(d-a)^2+4bc=(a+d)^2-4=(\tr h)^2-4$.

So, if $\tr h<2$ then $h$ has a unique fixed point $a\in\H_2$
(the other root of equation~\rf8, the conjugate, is not in $\H_2$).
The isometry $h$ is {\sl elliptic} and $a$ is the {\sl center of $h$}.

If $c=0$ or if $\tr h=2$ then equation~\rf8 has a unique root in
$\bar\R=\R\cup\{\infty\}$, hence $h$ has no fixed points. 
The isometry $h$ is {\sl parabolic}. This case contains the case $c=a-d=0$.

If $c\ne0$ and $\tr h>2$ then~\rf8 has two roots in $\bar\R$
hence $h$ has no fixed points. The isometry $h$ is {\sl hyperbolic}.
\med

{\noi\bf b.} \
The condition $h_{a,b,c,d}(i)=\frac{ai+b}{ci+d}=i$ implies 
$a=d$ and $c=-b$. The condition $ad-bc=1$ then implies that there exists
$\q\in\R$ such that $a=\cos\tfrac\q2$ and $c=\sin\tfrac\q2$.
\med

{\noi\bf c.} \
If $\f$ and $\psi\in\isomp \H_2$ send  $i$ to $z_0\in \fix h$
then, denoting $\f^{-1} h\f$ by $r_{\q}$, $\f^{-1}\psi$ by $r_\a$ 
(these are elements of $\isomp \H_2$ fixing $i$, hence  rotations),
we obtain 
$$
\psi^{-1}h\psi=\psi^{-1}\f r_\q\f^{-1}\psi=
r_\a^{-1}r_\q r_\a=r_\q=\f^{-1}h\f.
$$
{\noi\bf d.} \
The rotation of center $ix$ and angle $\q$ is $f=\f r_\q\f^{-1}$, where
$\f\in\isomp \H_2$ sends $i$ to $ix$. Let us choose $\f=h_{a,b,c,d}$ with 
$a=\sqrt x,$ $b=c=0$ and $d=\frac1{\sqrt x}$, i.e. the homography
associated with the matrix
$M(x)=\left(\begin{matrix}
\sqrt x&0\\
0&1/\sqrt x
\end{matrix}\right)$.
Then the matrices associated with $\f r_\q\f^{-1}$ are
$\pm M$ with 
$$
M=M(x)^{-1}R\big(\tfrac\q2\big)M(x)=\left(\begin{matrix}
\cos\tfrac\q2&-x\sin\tfrac\q2\smallskip\\
x^{-1}\sin\tfrac\q2&\cos\tfrac\q2
\end{matrix}\right).
$$
\ep
\sub{8.2}{Proof of Lemma~\reff{lem:libreSL2}}
Let $\T=\{M\in\SL(2,\Z)\tq\tr M=3\}$. 
We first prove the following preliminary result.
\lm{p8.2}{ 
Let $M\in\T$.
\be
\item One has $M^{-1}\in\T$.
\smallskip

\item 
\bee
\item
For any integer $n\ge1$, one has $M^n=\alpha_n\,M-\alpha_{n-1}\,I$,
with $\alpha_0=0$, $\alpha_1=1$, and $\alpha_{n}=3\alpha_{n-1}-\alpha_{n-2}$
for $n\ge2$.
\item
The sequence $(\alpha_n)_{n\in\N}$ takes its values in $\N$ and strictly
increases. Moreover, for any $n\in\N$ we have 
$\alpha_{n+1}-\alpha_{n}\ge2^n$ and $\,\alpha_n\ge2^n-1$.
\ee
\item For any integer $n\ge1$, one has $\tr M^n\ge2^{n+1}-1$.
\ee
}
\pr Item a results from the fact that
$$
M=\left(\begin{matrix}a&b\\c&d\end{matrix}\right)\;\Rightarrow
\;M^{-1}=\left(\begin{matrix}\hfill d&-b\\-c&\hfill a\end{matrix}\right).
$$
Item b~(i) is proved by induction. 
It is obvious for $n=1$ and also for $n=2$ (i.e. $M^2=3M-I$) by
the Cayley-Hamilton theorem.\index{Cayley-Hamilton theorem}\index{Theorem!Cayley-Hamilton}
Next, if we assume that $M^{n-1}=\alpha_{n-1}\,M-\alpha_{n-2}\,I$
for some $n\ge2$ then
$$
M^n=\alpha_{n-1}\,M^2-\alpha_{n-2}\,M=(3\alpha_{n-1}-
\alpha_{n-2})M-\alpha_{n-1}\,I=\alpha_n\,M-\alpha_{n-1}\,I.
$$
For item b~(ii), 
an easy induction shows that $(\a_n)$ is increasing and that
$\alpha_n\ge0$ for all $n\in\N$. Therefore we have
$$
\alpha_{n+1}-\alpha_{n}=2\alpha_{n}-\alpha_{n-1}
\ge2(\alpha_{n}-\alpha_{n-1})\ge2^n(\alpha_{1}-\alpha_{0})=2^n,
$$
hence
$$
\alpha_n=\alpha_n-\alpha_0=\sum_{k=0}^{n-1}(\alpha_{k+1}-\alpha_k)
\ge\sum_{k=0}^{n-1}2^k=2^n-1.
$$
For item c, considering the above, for any $n\ge1$, one has
$$
\tr M^n=\tr(\alpha_n\,M-\alpha_{n-1}\,I)=3\alpha_n-2\alpha_{n-1}
=\alpha_n+2(\alpha_n-\alpha_{n-1})\ge2^{n+1}-1.
$$
\ep

Let $G_0$ be the subgroup of $\SL(2,\Z) $ generated by the matrices
$$
A=\left(\begin{matrix}\hfill0&1\\-1&3\end{matrix}\right)\qquad
\mbox{and}\qquad B=\left(\begin{matrix}-1&-1\\
\hfill5&\hfill4\end{matrix}\right).
$$
The matrices $A$ and $B$ (as well as their inverses $A^{-1}$ and
$B^{-1}$) are elements of $\T$.
We will show that any matrix $M\in G_0\setminus\{I\}$ has a trace
different from $2$, which is equivalent to saying that $1$ is not
an eigenvalue of $M$ (since $\det M = 1$), i.e. that
$M-I$ is invertible.

\smallskip
Let $\e:\R\to\{-1,1\}$ denote the {\em sign function},
i.e. $\e(x)=1$ if $x\ge0$ and $\e(x)=-1$ if $x<0$.
Let $\gg$ denote the {\em product order} on $\SL(2,\Z)$:
We write $X\gg Y$ if $x_{ij}\ge y_{ij}$ for all $i,j\in\{1,2\}$. 
It is a partial order compatible with the addition and
the multiplication. Especially for $X,Y,X',Y'\in\SL(2,\Z)$ one has
$$
(X\gg Y\gg 0\;\mbox{ and }\;X'\gg Y'\gg 0)\;\Rightarrow\;XX'\gg YY'.
$$
We also have $X\gg Y\;\Rightarrow\;\tr X\ge\tr Y.$
\medskip

The following is inspired by \cite{n}, VIII 26,
pp.158--162, and proves Lemma~\reff{lem:libreSL2}.
\lm{l8.3}{
For $k,l\in\Z\setminus\{0\}$, one has
\begin{equation}
\e(kl)A^{\,k}B^{\,l}\gg \left(\begin{matrix}5&0\\
\hfill0&\hfill1\end{matrix}\right)\!.\label{ineg}
\end{equation}}
\pr Note first that 
$$
A^{-1}=\left(\begin{matrix}\hfill3&\hfill -1\\1&0\end{matrix}\right)\qquad
B^{-1}=\left(\begin{matrix}
\hfill4&\hfill 1\\
\hfill -5&-1\end{matrix}\right)\qquad
AB=\left(\begin{matrix}\hfill5&\hfill 4\\16&13\end{matrix}\right)
$$
and that
$$
AB^{-1}=\left(\begin{matrix}\hfill-5&\hfill-1\\
-19&-4\end{matrix}\right)\qquad
A^{-1}B=\left(\begin{matrix}\hfill-8&\hfill -7\\
-1&-1\end{matrix}\right)\;\;\mbox{ and }\;\;
A^{-1}B^{-1}=\left(\begin{matrix}\hfill17&\hfill4\\
\hfill 4&1\end{matrix}\right).
$$
Now let $k,l\ge1$.
With the notation of Lemma~\ref{p8.2}.b, 
one has
\begin{align*}
A^{\,k}B^{\,l}&=(\alpha_{k}\,A-\alpha_{k-1}\,I)(\alpha_l\,B-\alpha_{l-1}\,I)
\\
\noalign{\vskip2mm}
&=\alpha_k\alpha_l\,AB-\alpha_{k}\alpha_{l-1}\,A-
\alpha_{k-1}\alpha_{l}\,B+\alpha_{k-1}\alpha_{l-1}\,I
\\
\noalign{\vskip2mm}
&=\left(\begin{matrix}\hfill 5\alpha_{k}\alpha_{l}+
\alpha_{k-1}\alpha_{l}+\alpha_{k-1}\alpha_{l-1}
&4\alpha_{k}\alpha_{l}-\alpha_{k}\alpha_{l-1}+\alpha_{k-1}\alpha_{l}
\\
\noalign{\vskip2mm}
16\alpha_{k}\alpha_{l}+\alpha_{k}\alpha_{l-1}-5\alpha_{k-1}\alpha_{l}
&\;13\alpha_{k}\alpha_{l}-3\alpha_{k}\alpha_{l-1}-
4\alpha_{k-1}\alpha_{l}+\alpha_{k-1}\alpha_{l-1}
\end{matrix}\right).
\end{align*}
Therefore, since the sequence $(\alpha_n)$ strictly increases and
$\alpha_0=0$, we obtain
\begin{equation}
A^{\,k}B^{\,l}\gg \left(\begin{matrix}5&0\\
\hfill0&\hfill6\end{matrix}\right)\!.\label{in}
\end{equation}
Similarly, for $k,l\ge1$ we have
\begin{align*}
A^{\,k}B^{-l}&=\left(\begin{matrix}\hfill \alpha_{k-1}\alpha_{l-1}-
4\alpha_{k-1}\alpha_{l}-5\alpha_{k}\alpha_{l}
&-\alpha_{k}\alpha_{l-1}-\alpha_{k-1}\alpha_{l}-\alpha_{k}\alpha_{l}
\\
\noalign{\vskip2mm}
\alpha_{k}\alpha_{l-1}+5\alpha_{k-1}\alpha_{l}-19\alpha_{k}\alpha_{l}
&\alpha_{k-1}\alpha_{l-1}-3\alpha_{k}\alpha_{l-1}+\alpha_{k-1}\alpha_{l}-
4\alpha_{k}\alpha_{l}
\end{matrix}\right)
\\
\noalign{\vskip2mm}
A^{-k}B^{\,l}&=\left(\begin{matrix}\hfill \alpha_{k-1}\alpha_{l-1}-
3\alpha_{k}\alpha_{l-1}+\alpha_{k-1}\alpha_{l}-8\alpha_{k}\alpha_{l}
&\alpha_{k}\alpha_{l-1}+\alpha_{k-1}\alpha_{l}-7\alpha_{k}\alpha_{l}
\\
\noalign{\vskip2mm}
-\alpha_{k}\alpha_{l-1}-5\alpha_{k-1}\alpha_{l}-\alpha_{k}\alpha_{l}
&\alpha_{k-1}\alpha_{l-1}-4\alpha_{k-1}\alpha_{l}-\alpha_{k}\alpha_{l}
\end{matrix}\right)
\\
\noalign{\vskip2mm}
A^{-k}B^{-l}&=\left(\begin{matrix}\hfill \alpha_{k-1}\alpha_{l-1}-
3\alpha_{k}\alpha_{l-1}-4\alpha_{k-1}\alpha_{l}+17\alpha_{k}\alpha_{l}
&\alpha_{k}\alpha_{l-1}-\alpha_{k-1}\alpha_{l}+4\alpha_{k}\alpha_{l}
\\
\noalign{\vskip2mm}
-\alpha_{k}\alpha_{l-1}+5\alpha_{k-1}\alpha_{l}+4\alpha_{k}\alpha_{l}
&\alpha_{k-1}\alpha_{l-1}+\alpha_{k-1}\alpha_{l}+\alpha_{k}\alpha_{l}
\end{matrix}\right)
\end{align*}
from which one deduces the inequalities
$$-A^{\,k}B^{-l}\gg \left(\begin{matrix}5&0\\
\hfill0&\hfill3\end{matrix}\right)
\qquad  -A^{-k}B^{\,l}\gg \left(\begin{matrix}7&0\\
\hfill0&\hfill1\end{matrix}\right)
\quad\mbox{ and }\quad A^{-k}B^{-l}\gg \left(\begin{matrix}10&0\\
\hfill0&\hfill1\end{matrix}\right).$$
This and \rf{in} give the desired result.
\ep


\lm{p}{
Let $m\ge1$ be an integer and $k_1,\ldots k_m$,
$l_1,\ldots l_m$, be nonzero integers, except possibly for $k_1$ and $l_m$,
and such that $(k_1, l_1)\ne(0,0)$.
Then the matrix $\,M=A^{\,k_1}B^{\,l_1}\ldots A^{\,k_m}B^{\,l_m}\,$
satisfies $\vert\tr M\vert\ge 3$.
}
\pr
If $m=1$ then the result follows from
Lemma~\ref {p8.2}.c if $k_1l_1 = 0$, 
and from \rf{ineg} if $k_1$ and $l_1$
are both nonzero, because in this case
$$
\vert\tr (A^{\,k_1}B^{\,l_1})\vert=\e(k_1l_1)\tr (A^{\,k_1}B^{\,l_1})\ge6.
$$
Now suppose $m\ge2$. 
\begin{itemize}
\item[$\triangleright$]
If $k_1$ and $l_m$ are nonzero then
$$
\vert\tr M\vert=\tr\big(\e(k_1l_1)A^{\,k_1}B^{\,l_1}\times\cdots
\times\e(k_ml_m)A^{\,k_m}B^{\,l_m}\big)\ge
\tr\left(\begin{matrix}5^m&0\\0&\hfill1\end{matrix}\right)=5^m+1\ge3.
$$
\item[$\triangleright$]
If $k_1=l_m=0$ then
$$
\tr M=\tr(B^{\,l_1}A^{\,k_2}B^{\,l_2}\ldots 
A^{\,k_m})=\tr(A^{\,k_2}B^{\,l_2}\ldots A^{\,k_m}B^{\,l_1})
$$
and we fall in the previous case.
\item[$\triangleright$] 
If $k_1=0$ and $l_m\neq0$ then
$$
\tr M=\tr(B^{\,l_1}A^{\,k_2}B^{\,l_2}\ldots A^{\,k_m}B^{\,l_m})=
\tr(A^{\,k_2}B^{\,l_2}\ldots A^{\,k_m}B^{\,l_m+l_1}).
$$
If $l_m+l_1\neq0$ we fall into one of the previous cases
(depending on whether $m=2$ or $m>2$). If $l_m+l_1=0$ then
$$
\tr M=\tr(A^{\,k_2}B^{\,l_2}\ldots A^{\,k_m})=
\tr(A^{\,k_2+k_m}B^{\,l_2}\ldots A^{\,k_{m-1}}B^{\,l_{m-1}})
$$
and we iterate until we are reduced to an already considered case.
\item[$\triangleright$] 
If $k_1\neq0$ and $l_m=0$ then
$$
\tr M=\tr(A^{\,k_1}B^{\,l_1}\ldots A^{\,k_m})=
\tr(A^{\,k_1+k_m}B^{\,l_1}\ldots A^{\,k_{m-1}}B^{\,l_{m-1}})
$$
and we proceed as above.
\ep
\end{itemize}
\sub{8.5}{Proof of Lemma~\reff{l5.5}}
Let $ w $ be a nontrivial word of $ G_0 = \langle \sig, \tau \rangle $
of length $ \ell \ge1 $, i.e. an element $ w = a_1 \dots a_\ell $ with
$ a_i \in \{\sig, \sig^{- 1}, \tau, \tau^{- 1} \} $ and $ a_ia_{i + 1} \ne \id $
for all $ i \in \{1, \dots, \ell-1 \} $.
  The word $ w $ has one of the four forms
$ \sig^\pm \cdots \tau^\pm $, $ \tau^\pm \cdots \sig^\pm $, $ \sig^\pm \cdots \sig^\pm $
or $ \tau^\pm \cdots \tau^\pm $. The second form comes down to the first by
considering $ w^{- 1} $ instead of $ w $ and the last two forms
can be reduced to one of the first two by conjugation, unless
$ w $ is simply a power of $ \sig $ or $ \tau $, or a conjugate
of such a power, in which case the result comes from $ \q \notin \pi \Q $.
So we just have to show
that $ 1 $ is not an eigenvalue of $ w $ when $ w = \sig^\pm \cdots \tau^\pm $.
An easy induction shows that $ w $ has a matrix of the form
$$
\left(\begin{matrix}\hfill P&-Q&-R&-S\\Q&P&-S& R\\
R&S&P&-Q\\S&-R&Q&P\end{matrix}\right)
$$
where $ P, R $ are polynomials in $ \cos \q $ with integer coefficients
and $ Q, S $ are products of such polynomials with $ \sin \q $.

Since $ w $ is orthogonal, we have $ P^2 + Q^2 + R^2 + S^2 = 1 $.
A straightforward computation then shows that the characteristic polynomial of
$ w $ is $ \l^4-4P \l^3 + (4P^2 + 2) \l^2-4P \l + 1 $.
\med

We show below that the degree of $ P $ is equal to $ \ell $,
the length of $ w $.
If $ 1 $ were an eigenvalue of $ w $ we would have $ 4P^2-8P + 4 = 0 $,
a contradiction since $ \cos \q $ is transcendent.

In the sequel, the notation $ \as $ indicates that only the term of
highest degree in $ \cos \q $ has been retained.
Using $ \cos m \q \as2^{m-1} \cos^m \q $ and
$ \sin m \q \as2^{m-1} \cos^{m-1} \q \sin \q $, and denoting 
$ \cos \q $ by $ c $
and $ \sin \q $ by $ s $, we get, with $ \eps = \pm1 $ and $ \d = \pm1 $
$$
\sig^{\eps m}\as2^{m-1}c^{m-1}\left(\begin{matrix}\hfill
c&-\eps s&0&0\\
\eps s&c&0&0\\
0&0&c&-\eps s\\0&0&\eps s&c\end{matrix}\right)
\qquad
\tau^{\d k}\as2^{k-1}c^{k-1}\left(\begin{matrix}\hfill
c&0&0&-\d s\\0&c&-\d s&0\\0&\d s&c&0\\\d s&0&0&c
\end{matrix}\right).
$$
Multiplying both expressions and using $ s^2 = 1-c^2 \as-c^2 $ we get
\begin{equation*}
\sig^{\eps m}\tau^{\d k}\as A(m,k,\eps,\d)=2^{m+k-2}c^{m+k-1}
\left(\begin{matrix}\hfill
c&-\eps s&-\eps\d c&-\d s\\
\eps s&c&-\d s&\eps\d c\\
\eps\d c&\d s&c&-\eps s\\
\d s&-\eps\d c&\eps s&c\end{matrix}\right).
\end{equation*}
Let us denote $w=\sig^{\eps_1m_1}\tau^{\d_1k_1}\cdots\sig^{\eps_nm_n}\tau^{\d_nk_n}$
with $\eps_i,\d_i=\pm1$ and $m_i,k_i\ge1$; thus the length of $w$ is
$\ell=m_1+k_1+\cdots+m_n+k_n$. 

We assert that, for every integer $n\ge1$, there exist $ \xi_n $, $ \mu_n $,
$ \zeta_n $ and $ \nu_n $ equal to $ \pm1 $ such that
$\xi_{n}\zeta_{n}=\mu_{n}\nu_{n}$ and
$$
w\as A_n=2^{\ell-n-1}c^{\ell-1}
\left(\begin{matrix}
\xi_n c&\hfill-\mu_n s&-\zeta_n c&-\nu_n s\\
\mu_n s&\hfill\xi_n c&-\nu_n s&\hfill\zeta_n c\\
\zeta_n c&\hfill\nu_n s&\hfill\xi_n c&-\mu_n s\\
\nu_n s&\hfill-\zeta_n c&\hfill\mu_n s&\hfill\xi_n c\end{matrix}\right).
$$
By induction on $ n $, let us show that this is indeed the case.
For $ n = 1 $, this is because $ A_1 = A (m_1, k_1, \eps_1, \d_1) $
and because we have $ \xi_1 = 1 $, $ \mu_1 = \eps_1 $,
$ \zeta_1 = \eps_1 \d_1 $, and $ \nu_1 = \d_1 $.

Now assume that the property holds for an integer
$ n \ge1$ and let us check it for $ n + 1 $.
One sees that $A_nA(m_{n+1},k_{n+1},\eps_{n+1},\d_{n+1})\as A_{n+1}$ with
\begin{align*}
\xi_{n+1}&=\tfrac12(\xi_n+\eps_{n+1}\mu_n-\eps_{n+1}\d_{n+1}\zeta_n+
\d_{n+1}\nu_n)\\
\mu_{n+1}&=\tfrac12(\eps_{n+1}\xi_n+\mu_n+\d_{n+1}\zeta_n-
\eps_{n+1}\d_{n+1}\nu_n)\\
\zeta_{n+1}&=\tfrac12(\eps_{n+1}\d_{n+1}\xi_n+\d_{n+1}\mu_n+
\zeta_n-\eps_{n+1}\nu_n)\\
\nu_{n+1}&=\tfrac12(\d_{n+1}\xi_n+\eps_{n+1}\d_{n+1}\mu_n-
\eps_{n+1}\zeta_n+\nu_n).
\end{align*}
From these equalities, it is easily shown that
$\xi_{n+1}\zeta_{n+1}=\mu_{n+1}\nu_{n+1}$. 
It remains to check that $\xi_{n+1},\mu_{n+1},\zeta_{n+1},\nu_{n+1}$
are equal to $\pm1$. 
Noticing that $\zeta_n=\xi_{n}\mu_{n}\nu_{n}$,
one obtains
$$
\xi_{n+1}=\tfrac12(\xi_n+\eps_{n+1}\mu_n+\d_{n+1}\nu_n-
\eps_{n+1}\d_{n+1}\xi_{n}\mu_{n}\nu_{n})=\pm1
$$
since, in general, if $ a, b, c $ are equal to $ \pm1 $
then $ a + b + c-abc = \pm2 $.
In the same way we show that $ \mu_{n + 1} $, $ \zeta_{n + 1} $ and $ \nu_{n + 1}$
are equal to $ \pm1 $, which completes the proof.
\sub{8.4}{Solutions of the exercises}
{\noi\bf Exercise~\reff{e0}.} 
From Corollary~\reff{c4.5}, it is enough to show that $\isomp\R^2$
is fixating. Let $G\le\isomp\R^2$ be a \gaf.
If $G=\{\id\}$ then $G$ is obviously a \gag.
Otherwise let $f\in G\setminus\{\id\}$. Then $f$ is a rotation of center
$a\in\R^2$ and nonzero angle. Let $g\in G\setminus\{\id\}$ be arbitrary;
$g$ is a rotation of center $b\in\R^2$ and nonzero angle.
In addition we have $\v{[f,g]}=\id$, thus $[f,g]$ is the translation of vector 
$\v{cf(c)}$, with $c=g(a)$. Since $G$ is a \gaf, $\fix[f,g]$ is not
the empty set, hence $f(c)=c$, then $c=a$, $g(a)=a$, and $a=b$.
Therefore all the elements of $G\setminus\{\id\}$ are rotations of center
$a$. It follows that $\fix G=\{a\}\ne\emptyset$ and $G$ is a \gag.
\medskip

{\noi\bf Exercise~\reff{e1}.} 
Let $G=\langle g\rangle$, let $X$ be an arbitrary set and
$\rho:G\to\bij X$ be a morphism.
We have $\rho(G)=\langle\rho(g)\rangle$.
Let $H\le\rho(G)$. Then $H$ is cyclic generated by some $h$
hence $\fix H=\fix h$. If $H$ is a \gaf\;then $\fix h\neq\emptyset$ hence
$\fix H\neq\emptyset$ too, hence $H$ is a \gag. 
\\
Conversely, let $G$ be a noncyclic group. 
Let $G$ act by left multiplication on the set 
$X=\P(G)\setminus\{\emptyset,G\}$ of nontrivial subsets of $G$:
Precisely, let $\rho:G\to\bij X$ be defined by $\rho(g)(A)=gA=\{ga\tq a\in A\}$.
For every $g\in G$ the subset $\langle g\rangle$ is nontrivial since $G$ is
noncyclic and $\langle g\rangle$ is fixed by $\rho(g)$
hence $\rho(G)$ itself is a \gaf\;but $\rho(G)$ has no global fixed point:
Given $A\in X$, let $a\in A$ and $b\notin A$.
Since $b=(ba^{-1})a$, we have $b\in\rho(ba^{-1})(A)$ hence $\rho(ba^{-1})(A)\ne A$.
Thus $(X,G)$ is eccentric,  therefore nonfixating.
\med

{\noi\bf Exercise~\reff{e2}.} 
Denote $\rho(\Q)=\{g_1,\dots,g_n\}$.
Let $r_i=\tfrac{p_i}{q_i}\in\Q$ such that $\rho(r_i)=g_i$, with $p_i$ and $q_i$
relatively prime, and denote by $M$ the least common multiple of the $q_i$. 
Then $\rho(\Q)$ is cyclic, generated by $\rho\big(\tfrac1M\big)$,
hence fixating by Exercise~\reff{e1}.
\med

{\noi\bf Exercise~\reff{e3}.} 
Taking successively $z=x$  and $z=y$ in \rf6 
we obtain $d(x,m)$ and $d(y,m)\le\frac12\,d(x,y)$.
The triangular inequality  $d(x,y)\le d(x,m)+d(m,y)$ yields the
equalities.
If $m'$ is another point satisfying \rf6 then, applying \rf6 with
$z=m'$ and using $d(x,m')=d(y,m')=\frac12\,d(x,y)$,
we get $d(m,m')^2\le0$, hence $m=m'$.
\med

{\noi\bf Exercise~\reff{e4}.}
{\sl Theorem of M. Fr\'echet, P. Jordan and J. von
Neuman}~\cite{f,jn}\index{Fr\'echet, Maurice}\index{Jordan, Pascual}\index{Neuman@von Neumann, John}.
\be
\item
In a normed vector space, due to the fact that
$m=\tfrac12(x+y)$, the median inequality~\rf6 for $z=\v0$ reads as
$\|x+y\|^2\le2\big(\|x\|^2+\|y\|^2\big)-|\|x-y\|^2$.
Rewritten with $x=\frac12(a+b)$ and $y=\tfrac12(a-b)$,
this yields the opposite inequality.
\item
Set $\scal xy=\tfrac14\big(\|x+y\|^2-\|x-y\|^2\big)$.
It is enough to verify that $\scal{\;\,}{\;\,}$ is a positive
definite symmetric bilinear form.
Easily we have $\scal xy=\scal yx$,
$\scal xx\ge0$ and $\big(\scal xx=0\Rightarrow x=\v0\big)$.
It remains to prove that
$\scal{x+x'}y=\scal xy+\scal{x'}y$
and $\scal{\l x}y=\l\scal xy$.

First we have $\scal{\v0}y=0$.
By using \rf4\;we get
\begin{align*}
\scal{x+x'}y&+\scal{x-x'}y\\
&=\tfrac14\big(\|x+x'+y\|^2+\|x-x'+y\|^2-\|x+x'-y\|^2-\|x-x'-y\|^2\big)\\
&=\tfrac12\big(\|x+y\|^2+\|x'\|^2-\|x-y\|^2-\|x'\|^2\big)\\
&=2\scal xy.
\end{align*}
For $x'=x$ we obtain $\scal{2x}y=2\scal xy$.
\\
For $x=\tfrac12(a+b)$ and $x'=\tfrac12(a-b)$ we deduce
$\scal ay+\scal by=\scal{a+b}y$.
Therefore $\scal{kx}y=k\scal xy$, first
for integer $k$, then for rational $k$, finally for real $k$
by continuity of the norm and by density of $\Q$ in $\R$.
\ee
\med

{\noi\bf Exercise~\reff{e6}.}
Let $xyz$ be a hyperbolic triangle and $m$ the midpoint of the geodesic 
segment from $x$ to $y$.
Let $\a=\angle zmx$ and $\b=\angle zmy\in\,]0,\pi[$~;
one has $\a+\b=\pi$ hence $\cos\a+\cos\b=0$.
For short denote by $ab$ the geodesic distance between two points $a$ and $b$.
The cosine inequality for the triangle $myz$ is:
$mz^2+my^2-2mz.my\cos\a\le yz^2$. By using $my=\tfrac12xy$, we obtain
$mz^2+\tfrac14xy^2-mz.xy\cos\a\le yz^2$.
In the same way, one has $mz^2+\tfrac14xy^2-mz.xy\cos\b\le xz^2$.
By adding both inequalities, we get
$2mz^2+\tfrac12xy^2\le xz^2+yz^2$
which indeed corresponds to the median inequality~\rf6.
\med

{\noi\bf Exercise~\reff{ex:homographie1}.}
Let $f:x\mapsto\frac1{-x+3}$ and $g:x\mapsto\frac{-x-1}{5x+4}$ ;
$f$ and $g$ have distinct fixed points.
Let $k_1,l_1,\dots,k_m,l_m$ be relative integers all nonzero,
except possibly $k_1$ and $l_m$.
The product $h=f^{k_1}g^{l_1}\cdots f^{k_m}g^{l_m}$
is of the form $h:x\mapsto\frac{ax+b}{cx+d}$ with $ad-bc=1$. 
Looking for a fixed point $x$ of $h$ leads to the equation
$cx^2+(d-a)x-b=0$, whose discriminant is $\Delta=(a+d)^2-4$. 
From Lemma~\ref{p}, one has $|a+d|\ge3$ hence $\Delta\ge0$ hence $h$ has
always a fixed point. Therefore the group $\langle f,g\rangle$ is eccentric.
\med

{\noi\bf Exercise~\reff{e12}.}
Let $G$ be a \gaf\, of $\isom\H_2$. Let us identify $\H_3$ with $\R\times\H_2$.
For any $g\in G$, let $\f(g):\H_3\to\H_3,\,(x,y)\mapsto(x,g(y))$.
We easily check that $\f(g)\in\isom\H_3$ and that
$\fix\f(g)=\R\times\fix g\ne\emptyset$.
Thus $H:=\f(G)$ is a \gaf\, of $\isom\H_3$, hence a \gag\, since
$\isom\H_3$ is globalizing, and $\fix H$ is of the form $\R\times A$.
We obtain $\fix G=A\ne\emptyset$ hence $G$ is a \gag\, of $\isom\H_2$.
\med

{\noi\bf Exercise~\reff{e4b}.}
Let $f:\R^n\to\R^n$ preserving the Euclidean distance. The equality
$$
\scal xy=\tfrac12\big(\|x\|^2+\|y\|^2-\|x-y\|^2\big)
$$
shows that the function $g:x\mapsto f(x)-f(\v0)$ preserves the scalar product:
One has $\scal{g(x)}{g(y)}=\scal xy$ for every $x,y\in\R^n$.
Thus, for every $x,y\in\R^n$ and every $\l\in\R$, we have
\begin{align*}
\|g(x+\l y)-g(x)-\l g(y)\|^2=\;&
\|g(x+\l y)\|^2+\|g(x)\|^2+\l^2\|g(y)\|^2-\\
&2\scal{g(x+\l y)}{g(x)}
-2\l\scal{g(x+\l y)}{g(y)}+2\l\scal{g(x)}{g(y)}\\
=\;&\|x+\l y\|^2+\|x\|^2+\l^2\|y\|^2-\\
&2\scal{x+\l y}x-
2\l\scal{x+\l y}y+2\l\scal xy\\
=\;&\|x+\l y-x-\l y\|^2=0.
\end{align*}
Therefore $g$ is linear hence $f$ is affine.
Since $g$ preserves the Euclidean distance, it is injective.
Since we are in finite dimension, $g$ is surjective.
It is the same with $f$.
\med

{\noi\bf Exercise~\reff{e5}.}
Let $E$ be a normed vector space and $f:E\to E$  be a continuous function such that 
\eq{77}{
\forall a,b\in E,\quad f\big(\tfrac12(a+b)\big)=\tfrac12\big(f(a)+f(b)\big).
}
Let $g:E\to E,\;x\mapsto f(x)-f(\v0)$.
Then $g$ is continuous and satisfies~\rf{77} and $g(\v0)=\v0$.
Since $g(x)=g\big(\tfrac12(\v0+2x)\big)=\tfrac12\,g(2x)$,
we can rewrite~\rf{77} as
$$
\forall a,b\in E,\quad g(a+b)=g(a)+g(b).
$$
It follows that $g(ax)=ag(x)$ for all $a\in\N$, then for all 
$a\in\Q$, finally for all $a\in\R$ by continuity of $g$ and by
density of $\Q$ in $\R$, hence $g$ is linear, hence $f$ is affine.
\med

{\noi\bf Exercise~\reff{exo:mazurulam}.}
{\sl Mazur-Ulam Theorem}~\cite v\index{Mazur, Stanislaw}\index{Ulam, Stanislaw}.
\be
\item
\bee
\item
For every $g\in W_{a,b}$, one has 
$$
\Vert g(m)-m\Vert\le\Vert g(m)-g(a)\Vert+\Vert a-m\Vert
=2\Vert a-m\Vert=\Vert a-b\Vert,
$$
hence $\l\le\Vert a-b\Vert$.
\item
Let $g\in W_{a,b}$. Since $s_m$ is an isometry which fixes $m$, one has
\begin{align*}
\Vert g^*(m)-m\Vert&=\Vert s_m\,g^{-1}s_m\,g(m)-m\Vert
=\Vert g^{-1}s_m\,g(m)-m\Vert
\\&
=\Vert s_m(g(m))-g(m)\Vert
=2\Vert g(m)-m\Vert.
\end{align*}
\item
From the above, on the one hand $\l$ is finite, and on the other hand,
for all $g\in W_{a,b}$, since $s_m$ permutes $a$ and $b$, we have
$g^*\in W_{a,b}$, hence $2\Vert g(m)-m\Vert\le\l$, hence $2\l\le\l$,
hence $\l=0$. As a consequence we have $g(m)=m$ for all $g\in W_{a,b}$.
\ee
\smallskip

\item
\bee
\item One has 
$h(a)=s_m\,f^{-1}s_{m'}(f(a))=s_m\,f^{-1}(f(b))=s_m(b)=a$.
Similarly one proves $h(b)=b$, hence $h\in W_{a,b}$, therefore from {\bf a}
one has $h(m)=m$.

It follows that $f^{-1}s_{m'}f(m)=s_m(m)=m$, then $s_{m'}(f(m))=f(m)$.
Since $s_{m'}$ admits $m'$ as unique fixed point, one has $f(m)=m'$.
\item
From item~(i), for all $a,b$ in $E$ we have
$f(\frac12(a+b))=\frac12(f(a)+f(b))$.
Since $f$ is continuous, we deduce that $f$ is affine by Exercise~\reff{e5}.
\ee
\ee
\med

{\noi\bf Exercise~\reff{e9}.}
{\sl Kakutani Theorem in finite dimension}~\cite a\index{Kakutani, Shizuo}.
\be
\item 
One easily checks that $\|\;\;\|$ is a norm.
In order to prove that $\|\;\;\|$ is strictly convex on $E$ we
consider $x,y\in E$ such that $x\ne\v0$ and $y\notin\R^+x$.
Then we have, for all $g\in G$, 
$g(x)\ne\v0$ and $g(y)\notin\R^+g(x)$.
Put $\f(g)=\|g(x)\|_2+\|g(y)\|_2-\|g(x)+g(y)\|_2$.
The function $\f$ is continuous and takes positive values 
on the compact set $G$, hence $\f$ is bounded below by some $\d>0$.
Thus, for all $g\in G$,
$\|g(x)+g(y)\|_2\le\|g(x)\|_2+\|g(y)\|_2-\d\le\|x\|+\|y\|-\d$ hence,
taking the supremum: $\|x+y\|\le\|x\|+\|y\|-\d$.

Finally $\|g(x)-g(y)\|=\sup_{h\in G}\|hg(x)-hg(y)\|_2=
\sup_{k\in G}\|k(x)-k(y)\|_2=\|x-y\|$
since $G$ is a group, showing that every element of $G$ is an isometry.
\item
Since $K$ is convex, the Ces\`aro mean
$\sigma_n=\tfrac1n(x_1+\cdots+x_n)$ is in $K$.
By compactness, there exists a subsequence $(\sigma_{n_k})_{k\in\N}$
tending to some $a\in K$.
We have 
$$
f(\sigma_{n_k})=\tfrac1{n_k}(x_2+\cdots+x_{n_k+1})=\sigma_{n_k}-
\tfrac1{n_k}(x_{n_k+1}-x_1),
$$
hence the sequence $(f(\sigma_{n_k}))_{k\in\N}$ tends to $a$ too.
Since $f$ is an endomorphism on a finite dimensional space, $f$ is continuous
hence $(f(\sigma_{n_k}))_{k\in\N}$ tends to $f(a)$,
from which we deduce that $f(a)=a$.
\item
Given $g\in G$, set $V_g=\{x\in K\tq g(x)\neq x\}$ and assume by contradiction
that for all $x\in K$ there exists $g\in G$ such that $g(x)\neq x$.
\bee
\item
Since an isometry is continuous, each $V_g$ is an open subset of $K$.
By assumption, $K$ is the union of the $V_g$, $g\in G$, hence by compactness
there exist $g_1,\ldots,g_N\in G$ such that $K=V_{g_1}\cup\cdots\cup V_{g_N}$.
\item
Let $f=\tfrac1N(g_1+\cdots+g_N)$. By convexity, one has $f(K)\subseteq K$
hence there exists $a\in K$ such that $f(a)=a$ by item b.
\item
One has $\|Na\|=\|g_1(a)+\cdots+g_N(a)\|\le
\|g_1(a)\|+\cdots+\|g_n(a)\|=N\|a\|$ because the $g_k$ are isometries,
therefore the inequality is an equality, hence the $g_k(a)$ are positively
collinear by strict convexity of the norm, hence  all equal one to each other, 
hence  all equal to $a$. As a consequence the point $a$ would be in none of 
the $V_{g_k}$, a contradiction.
\ee
\ee
\med

{\noi\bf Exercise~\reff{e5.3}.}
{\sl Isometries of the sphere.}\index{extension of an isometry}
\be
\item 
For all $x,y\in\S_n$, one has $\scal{f(x)}{f(y)}=
\cos\,d\big(f(x),f(y)\big)
=\cos\,d(x,y)=\scal xy$.
\item
For $x$ or $y=\v0$, one has $\big\langle\~f(x)\,|\,\~f(y)\big\rangle=0=\langle x\,|\,y\rangle$.
For $x$ and $y\ne\v0$, one has
$$
\big\langle\~f(x)\,|\,\~f(y)\big\rangle=
\|x\|\,\|y\|\,\big\langle f\big(\tfrac x{\|x\|}\big)\,\big|\,
f\big(\tfrac y{\|y\|}\big)\big\rangle=
\|x\|\,\|y\|\,\big\langle\tfrac x{\|x\|}\,\big|\,\tfrac y{\|y\|}\big\rangle
=\scal xy.
$$
\item
Firstly we have
$$\|\~f(x)\|^2=\big\langle\~f(x)\,\big|\,\~f(x)\big\rangle=\langle x\,|\,x\rangle=\|x\|^2$$
for all $x\in\R^{n+1}$. 
Then, for all $x,y\in\R^{n+1}$
$$
\|\~f(x)-\~f(y)\|^2=
\|\~f(x)\|^2+\|\~f(y)\|^2-2\big\langle\~f(x)\,|\,\~f(y)\big\rangle=
\|x\|^2+\|y\|^2-2\langle x\,|\,y\rangle=\|x-y\|^2.
$$
By Exercise~\reff{e4b}, $\~f$ is also affine and surjective hence an isometry.
Since $\~f(\v0)=\v0$, $\~f$ is linear hence it is
determined by its values on a subset spanning $\R^{n+1}$,
for instance $\S_n$, hence the uniqueness follows.
\ee
\med

{\noi\bf Exercise~\reff{e54}.}\index{extension of an isometry}
{\sl Extension of an isometry.}
\be
\item\lb{e54a}
One has $\scal{f(a_i)-f(a_0)}{f(a_j)-f(a_0)}=
\|f(a_i)-f(a_0)\|^2+\|f(a_j)-f(a_0)\|^2
-\frac12\|f(a_i)-f(a_j)\|^2=\|a_i-a_0\|^2+\|a_j-a_0\|^2
-\frac12\|a_i-a_j\|^2=\scal{a_i-a_0}{a_j-a_0}$.
\item
One has
\begin{align*}
\|\~f(x)-\~f(y)\|^2
&=\Big\|\sum_{i=0}^n\big(\l_i(x)-\l_i(y)\big)f(a_i)\Big\|^2\\
&=\Big\|\sum_{i=0}^n\big(\l_i(x)-\l_i(y)\big)\big(f(a_i)-f(a_0)\big)\Big\|^2\\
&=\sum_{i=0}^n\sum_{j=0}^n\big(\l_i(x)-\l_i(y)\big)\big(\l_j(x)-\l_j(y)\big)
\scal{f(a_i)-f(a_0)}{f(a_j)-f(a_0)}\\
&=\sum_{i=0}^n\sum_{j=0}^n\big(\l_i(x)-\l_i(y)\big)\big(\l_j(x)-\l_j(y)\big)
\scal{a_i-a_0}{a_j-a_0}\\
&=\Big\|\sum_{i=0}^n\big(\l_i(x)-\l_i(y)\big)(a_i-a_0)\Big\|^2
=\|x-y\|^2.
\end{align*}
By Exercise~\reff{e4b}, $\~f$ is also surjective hence an isometry.
\item
Of course we have $\~f(a_i)=f(a_i)$ for all $i\in\{0,\dots,n\}$,
by uniqueness of the barycentric coordinates. 
Let $a\in A$; from above, one has $\|\~f(a)-f(a_i)\|=\|a-a_i\|$
for all $i\in\{0,\dots,n\}$ hence $\|\~f(a)-f(a_i)\|=\|f(a)-f(a_i)\|$
since $f$ is an isometry of $A$.
Thus $f(a_i)\in\md\big(f(a),\~f(a)\big)$ for all $i\in\{0,\dots,n\}$
but $\md\big(f(a),\~f(a)\big)$ is an affine subspace
hence $\~f(a)=\sum_{i=0}^n\l_i(a)f(a_i)\in\md\big(f(a),\~f(a)\big)$,
hence $\~f(a)=f(a)$, hence $\~f$ indeed extends $f$.
Since isometries of an affine space are affine maps by Exercise~\reff{e4b},
$\~f$ is affine and it is the only affine map of $\aff(a_0,\dots,a_n)$
taking the values $f(a_i)$ at points $a_i$ hence the uniqueness follows.
\ee
\med

{\noi\bf Exercise~\reff{e15}.}\index{alternating group}
{\sl The alternating group $\A_7$ is nonfixating.}
\\
The group $\A_4$ has a subgroup isomorphic to the Klein group,
\[
K=\{\id,(12)(34),(13)(24),(14)(23)\},
\]
such that:
\begin {itemize}
\item [$ \triangleright $]
The quotient $\A_4/K$ is isomorphic to $ \Z / 3 \Z $,
\item [$ \triangleright $]
Any element of $ \A_4 \setminus K $ has at least one fixed point 
(these are all cycles of order 3).
\end {itemize}
An eccentric subgroup $ G = \langle f, g \rangle $ of $ \A_7 $ is 
then obtained by taking two even permutations which act separately on 
$\{1,2,3,4\}$ and on $\{5,6,7\}$:
We choose for $ f $ an element of $ K $ on $ \{1,2,3,4 \} $
and the identity on $ \{5,6,7 \} $, say $ f = (12) (34) $,
and for $g$ an element of $\A_4\setminus K$ on $\{1,2,3,4\}$ 
and an element of $\A_3$ without fixed points on $\{5,6,7\}$, 
say $g=(123)(567)$.
We obtain for $G$ a 12-element semi-direct product of $K$ and $\A_3$.

Any element of $G$ whose restriction to $\{1,2,3,4\} $ has no fixed
points, i.e. is in $K\setminus\{\id\}$, is the identity on $\{5,6,7\} $,
so all the elements of $ G $ have a fixed point but $f$ and $g$ have
no common fixed points.
\med

{\noi \bf Exercise~\reff{e16}.}
{\sl The action of the group $\GL(3,\F_2)$  on $X=\F_2^3\setminus\{\v0\}$
is nonfixating.}\\
\be
\item
Identify $X$ with $\{1,\dots,7\}$ by
\[
e_i\mapsto i,\,e_{123}\mapsto4,\,e_{23}\mapsto5,
\,e_{13},\mapsto6,\,e_{12},\mapsto7.
\]
With this identification we obtain $f=(123)(567)$ and $g=(14)(67)$ which
are two even permutations.
\item
Since the elements of $G$ are linear and $Y$ generates the vector space
$\F_2^3$, if the restriction to $Y$ of $ h\in G$ is the identity map
then $h$ is the identity map. In addition, $f$ and $g$ send $Y$ to $Y$
hence the map which sends $h$ to its restriction to $Y$ is a one-to-one
morphism from $G$ to $\sy_Y$.
\item
We find
\[
(fg)^2=(12)(34),\qquad(gf)^2=(13)(24),\qquad(fgf)^2=(14)(23).
\]
By item b, there is no other element of $G$ whose restriction to
$\{1,2,3,4\} $ is a double transposition.
\item
As $ \fix f\cap\fix g=\emptyset$, in order to prove that the action of $G$
is eccentric, it is enough to check that $\fix h\neq\emptyset$ for all $h$
in $G$. The elements of $\sy_4$ without fixed point are the double
transpositions and the cycles of order 4.
From item c, if the restriction to $\{1,2,3,4\}$ of $h\in G$
is a double transposition then $\fix h=\{5,6,7\}$.
If the restriction of $h\in G$ to $\{1,2,3,4\}$ is a cycle of order 4,
since $h$ is an even permutation, its restriction to $\{5,6,7\}$ is a
transposition, so also has a fixed point.
Thus the action of $G$ is eccentric which shows that the action
$\GL(3,\F_2)$ on $\F_2^3\setminus\{\v0\}$ is nonfixating.
\item
It is easy to check that the above morphism is an isomorphism.
So the group generated by $(123)(567)$ and $(14)(67)$ is an
eccentric subgroup of order $24$ of $\A_7$ therefore is different
from the subgroup of Exercise~\reff{e15}.
\ee
\med
\np
{\noi\bf Exercise~\reff{e17}.}
{\sl The action of $\GL(d,\k)$ on $X=\k^d\setminus\{\v0\}$
is nonfixating, except in obvious cases.}
\be
\item
Let $G$ be the set of matrices of the form
$A=\left(\begin{matrix}a&b\\0&\hfill1\end{matrix}\right)$ with
$a\in\k^*=\k\setminus\{0\}$ and $b\in\k$. It is a group.
One finds that $\fix A$ is the straight line of equation
$x=0$ if $b=0$ and the straight line of equation $y=(1-a)x/b$ if
$b\ne0$ (more precisely the restrictions to $X$ of these lines).
Thus $G$ is a \gaf.
For $\card\k\ge3$, choose $a\in\k\setminus\{0,1\}$ and set
$A=\left(\begin{matrix}a&0\\0&\hfill1\end{matrix}\right)$ and
$B=\left(\begin{matrix}1&1\\0&\hfill1\end{matrix}\right)$.
Then $\fix A\cap\fix B$ is empty in $X$ hence $G$ is eccentric.
\item
Given an eccentric group $G$ of $\GL(d,\k)$, let
$\~G$ be the group of block-matrices of the form
$\~A=\left(\begin{matrix}A&b\\ \ta0&\hfill1\end{matrix}\right)$
where $A\in G$, $b\in\k^d$ is a column-vector, and $\ta0$ is the zero 
row-vector.
If $\~x=\left(\begin{matrix}x\\ \hfill y\end{matrix}\right)$
with $x\in\k^d$ and $y\in\k$ then
$\~A\~x=\left(\begin{matrix}Ax+yb\\y\end{matrix}\right)$
hence $\left(\begin{matrix}x\\ \hfill0\end{matrix}\right)$ 
is in $\fix\~A$ for all $x\in\fix A$, which shows
that $G$ is a \gaf.
On the other hand an element $\~x\in\fix\~G$ cannot be of the form
$\~x=\left(\begin{matrix}x\\ \hfill0\end{matrix}\right)$ (otherwise
$x$ would be in $\fix G$ which is empty by assumption) and,
if $\~x=\left(\begin{matrix}x\\ \hfill y\end{matrix}\right)$ 
were in $\fix\~G$ with $y\ne0$, we would have  $(I-A)x=yb$
for all $A\in G$ and all $b\in\k$, a contradiction since $\k$
has at least two elements.
\item
One easily checks that $\GL(1,\k)$ is fixating.
If $\k\ne\F_2$ then $\GL(2,\k)$ is nonfixating by item a
hence $\GL(d,\k)$  is nonfixating for all $d\ge2$ by item b.
If $\k=\F_2$ then $\GL(2,\k)$  is fixating seen as a
subgroup of $\sy_3$ but $\GL(3,\k)$ is nonfixating by
Exercise~\reff{e16} and $\GL(d,\k)$ is nonfixating for all 
$d\ge3$ by item b.
\ee
%
%
%

\printindex
\end{document}